\documentclass[a5paper,11pt]{book}

\usepackage{setspace} 

\usepackage{amsmath} 
\usepackage{amssymb}
\usepackage{theorem}

\usepackage{xspace} 
\usepackage{cite} 

\usepackage{endnotes} 

\renewcommand{\notesname}{Comments}
\makeatletter
\renewcommand{\enoteheading}%
{\chapter*{\notesname
  \@mkboth{\notesname}{\notesname}
  \addcontentsline{toc}{chapter}{Comments} }%
     \leavevmode\par\vskip-\baselineskip}

\renewcommand{\@makeenmark}{\hbox{$^{\@theenmark )}$}} 

\renewcommand{\enoteformat}{\rightskip\z@ \leftskip\z@ \parindent=1.8em
     \leavevmode\llap{\hbox{$^{\@theenmark )}$}}} 

\makeatother

\renewcommand{\thechapter}{\Roman{chapter}}

\renewcommand{\chaptermark}[1]{\markboth{\chaptername\ \thechapter. \ #1}{}}

\makeatletter
\renewcommand{\@makechapterhead}[1]{%
  \vspace*{20\p@}
  {\parindent \z@ \raggedright \normalfont
    \ifnum \c@secnumdepth >\m@ne
      \if@mainmatter
        \large \bfseries \@chapapp\space \thechapter
        \par\nobreak
        \vskip 10\p@
      \fi
    \fi
    \interlinepenalty\@M
    \LARGE \bfseries #1\par\nobreak
    \vskip 40\p@
  }}
\renewcommand{\@makeschapterhead}[1]{%
  \vspace*{20\p@}%
  {\parindent \z@ \raggedright
    \normalfont
    \interlinepenalty\@M
    \LARGE \bfseries  #1\par\nobreak
    \vskip 40\p@
  }}
\makeatother

\makeatletter
\def\@chapter[#1]#2{\ifnum \c@secnumdepth >\m@ne
                       \if@mainmatter
                         \refstepcounter{chapter}%
                         \typeout{\@chapapp\space\thechapter.}%
                         \addcontentsline{toc}{chapter}%
                                   {\protect \chaptername\ \thechapter.\ \ #1}
                       \else
                         \addcontentsline{toc}{chapter}{#1}%
                       \fi
                    \else
                      \addcontentsline{toc}{chapter}{#1}%
                    \fi
                    \chaptermark{#1}%
                    \addtocontents{lof}{\protect\addvspace{10\p@}}%
                    \addtocontents{lot}{\protect\addvspace{10\p@}}%
                    \if@twocolumn
                      \@topnewpage[\@makechapterhead{#2}]%
                    \else
                      \@makechapterhead{#2}%
                      \@afterheading
                    \fi}
\makeatother

\makeatletter
\def\@sect#1#2#3#4#5#6[#7]#8{%
  \ifnum #2>\c@secnumdepth
    \let\@svsec\@empty
  \else
    \refstepcounter{#1}%
    \protected@edef\@svsec{\@seccntformat{#1}\relax}%
  \fi
  \@tempskipa #5\relax
  \ifdim \@tempskipa>\z@
    \begingroup
      #6{%
        \@hangfrom{\hskip #3\relax\@svsec}%
          \interlinepenalty \@M #8\@@par}%
    \endgroup
    \csname #1mark\endcsname{#7}%
    \addcontentsline{toc}{#1}{%
      \ifnum #2>\c@secnumdepth \else
        \protect\numberline{\S \csname the#1\endcsname}
      \fi
      #7}%
  \else
    \def\@svsechd{%
      #6{\hskip #3\relax
      \@svsec #8}%
      \csname #1mark\endcsname{#7}%
      \addcontentsline{toc}{#1}{%
        \ifnum #2>\c@secnumdepth \else
          \protect\numberline{\csname the#1\endcsname}%
        \fi
        #7}}%
  \fi
  \@xsect{#5}}
\makeatother

\usepackage[russian,english]{babel}

\newtheorem{Pa}{Paper}[section]

\newtheorem{corollary}[Pa]{Corollary}

\newtheorem{definition}[Pa]{Definition}

\newtheorem{lemma}[Pa]{Lemma}

\begin{document}

\renewcommand{\contentsname}{Table of contents}
\renewcommand{\bibname}{References}

\thispagestyle{empty}

\begin{center}
{\large BELARUSIAN STATE UNIVERSITY}
\quad \\
\quad \\
\quad \\
\bigskip {\bf \LARGE Alexander Kiselev }\medskip\\
{\large \emph{e-mail}\ \ aakiselev@yahoo.com}
\quad \\
\quad \\
\quad \\
{\bf \Huge Inaccessibility}
\quad \\
\quad \\
{\bf \Huge and}
\quad \\
\quad \\
{\bf \Huge Subinaccessibility}\\
\quad \\
\quad \\
{\Large In two parts} \\
\smallskip
{\Large Part I} \\
\quad \\
\quad \\
{\em \large Second edition} \\
\smallskip
{\em \large enriched and improved} \\
\quad \\
\quad \\
\quad \\
\quad \\
\quad \\
\quad \\
\quad \\
\quad \\
\quad \\
{Minsk} \\
\smallskip
{``Publishing center of BSU''} \\
\smallskip
{2008}

\end{center}

\newpage

\thispagestyle{empty}

\noindent UDK 510.227
\\

{\small \textbf{Kiselev, Alexander.} Inaccessibility and
Subinaccessibility. In 2 pt. Pt. 1 / Alexander Kiselev. -- 2nd
ed., enrich. and improv. -- Minsk: Publ. center of BSU, 2008. --
\pageref{end} p. -- ISBN 978-985-476-596-9. }
\\

{\footnotesize The work presents the first part of second edition
of the previous edition of 2000 under the same title containing
the proof (in $ZF$) of the nonexistence of inaccessible cardinals,
now enriched and improved.

This part contains the apparatus of subinaccessible cardinals and
its basic tools -Ч theories of reduced formula spectra and
matrices, disseminators and others Ч- which are used in this proof
and are set forth now in their more transparent and simplified
form. Much attention is devoted to the explicit and substantial
development and cultivation of basic ideas, serving as grounds for
all main constructions and reasonings.

Ref. 26. }
\begin{center}
{R\;e\;f\;e\;r\;e\;e\;s:}
\\
Prof. {\em Petr P. Zabreiko;}
\\
Prof. {\em Andrei V. Lebedev}
\quad \\
\quad \\
\quad \\
\noindent Mathematics Subject Classification (1991): \\
03E05, 03E15, 03E35, 03E55, 03E60
\quad \\
\quad \\
\end{center}

\quad \\
\quad \\
\quad \\
\quad \\
\quad \\
\quad \\
\quad \\
\quad \\
\quad \\
\noindent {\scriptsize {\bf ISBN 978-985-476-596-9 (pt. 1)} \hfill
$\copyright$ Kiselev Alexander A., 2000 }

\noindent {\scriptsize {\bf ISBN 978-985-476-597-6}  \hfill
$\copyright$ Kiselev Alexander A., 2008, with modifications}

\newpage
\chapter*{Acknowledgements}

The author sends his first words of deep gratitude to
Hanna~R.~Ca\-lien\-do for understanding of the significance of the
theme, for hearty encouraging help and for the help in promotion of the work.

The very special appreciation is expressed to Prof. Sergei~R.
Kogalovskiy, who had taught the author the Hierarchy Theory, and
to Prof. Akihiro Kanamori for their valuable and strong
encouragement that gave the necessary impetus for the completing
the work.

The author would like to express his very special gratitude to
Prof. Andrei~V.~Lebedev and Prof. Petr~P.~Zabreiko who have
supported for many years the advancing of the work; the very
strong intellectual and moral debt is owed to both of them for the
spiritual and practical help.

The deep gratitude is expressed to Prof. Alexander~V.~Tuzikov and
Dr.~Yuri~Prokopchuk who gave the author extensive and expert help
in typing of the text.

Prof. Vasiliy~M.~Romanchack gave the author the financial and
moral help during the most complicated period of investigations of
the theme and the author sends to him many words of deep
thankfulness.

Many thanks are also to Ludmila Laptyonok who went through many
iterations of the difficult typing of previous works of the
author, preparing this work.

The deep hearty thankfulness goes to Nadezhda~P.~Zabrodina for
the support and encouragement for a long time, without which the
work would be considerably hampered.

The host of other people who over the years provided promotion of
the work or encouragement, is too long to enumerate and thanks are
expressed to all of them.

\newpage
{} \thispagestyle{empty}

\newpage

\begin{spacing}{1.5}

\tableofcontents

\end{spacing}

\newpage
{} \thispagestyle{empty}

\newpage
\chapter*{Introduction}
\addcontentsline{toc}{chapter}{Introduction}
\markboth{Introduction}{Introduction}

\setcounter{equation}{0}

The notion of inaccessibility in its various forms belongs to the
most important concepts of cultural traditions since ancientry and
consists in the representation about an {\it extraordinary} large
phenomenon that cannot be reached by means of ``lesser power''
than this phenomenon itself.
\\
In this meaning it appears to be one of the basic archetypes of
the Mankind.
\\

In the beginning of 20 century this concept have received its
sufficiently adequate expression in foundations of mathematics in
the following way. Cofinality of a cardinal \ $k$ \ is the minimal
ordinal \ $cf(k)$ \ that can be mapped by some function \ $f$ \
into \ $k$ \ such that \ sup$rng(f)=k$; \ the cardinal \ $k$ \ is
named regular \textit{iff} \ $cf(k)=k$, \ otherwise singular; it
is named weakly inaccessible \textit{iff} it is uncountable,
regular and for every \ $\alpha < k$ \ the next after \ $\alpha$ \
cardinal is still \ $ < k $. The notion of strongly inaccessible
cardinal is obtained by replacing the latter condition with the
stronger one: \ $\forall \alpha < k \ \ 2^{\alpha}< k $.
\\

Weakly inaccessible cardinals were introduced by Felix
Hausdorf~\cite{Hausdorf} in 1908; Paul Mahlo~\cite{Mahlo1, Mahlo2,
Mahlo3} had studied stronger limit points, Mahlo cardinals, in
1911-1913. Later on there appeared strongly inaccessible cardinals
introduced by Waclaw Sierpinski, Alfred Tarski~\cite{Sierp} and
Ernst Zermelo~\cite{Zermelo} in 1930. Since that starting period a
lot of inaccessible cardinals have appeared -- weakly compact,
indescribable, Erd\"{o}s cardinals, J\'{o}nsson, Rowbottom, Ramsey
cardinals, measurable, strong, Woodin, superstrong, strongly
compact, supercompact, extendible, Vo\-p\u{e}n\-ka's Principle,
almost huge, huge, superhuge cardinals etc., forming the hierarchy
corresponding to the substantial strengthening of an ``extent of
inaccessibility'' -- all of them have received the common name of
{\it ``large cardinals''.}
\\

So, the tradition of inaccessibility is firmly established in
foundations of mathematics serving as a confirmation to the faith
in the concept that the main trend of nowadays set theory
development lies in forming systems of large cardinal hypotheses,
natural models for these hypotheses and generic extensions of such
models. At any rate, there is firm belief, that:

``All of the large cardinal hypotheses ... can be given more or
less convincing justifications ... -- convincing enough, at any
rate, that almost no one expects them to be refused.

There are other, more problematic, large cardinal hypotheses. For
these there is as yet no convincing justification. Far from being
refused, however, these hypotheses have led to some extremely
beautiful mathematics.'' (Kunen~\cite{Kunen}).

In addition to a very beautiful justification of set theory as a
whole this faith is based on a lot of well-known results pointing
to close interrelations (by the relative consistency) between
large cardinal hypotheses, Axiom of Determinacy, regular set
properties, infinitary combinatorics, infinitary languages and
others. A big amount of results of this kind is expounded in
monographs of Frank Drake~\cite{Drake} and Akihiro
Kanamori~\cite{Kanamori}; the latter, remarkably exhaustive and
contemporary, contains also the outstanding demonstration of
diverse large cardinal ideas development.
\\
``The investigation of large cardinal hypotheses is indeed a
mainstream of modern set theory, and they have been found to play
a crucial role in the study of definable sets of reals, in
particular their Lebesgue measurability. Although formulated at
various stages in the development of set theory and with different
incentives, the hypotheses were found to form a $linear$ hierarchy
reaching up to an inconsistent extension of motivating concepts.
All known set-theoretic propositions have been gauged in this
hierarchy in terms of consistency strength, and the emerging
structure of implications provides a remarkably rich, detailed and
coherent picture of the strongest propositions of mathematics as
embedded in set theory.'' (Kanamori~\cite{Kanamori}).
\\

\noindent This work constitutes the first part of the second
edition of~\cite{Kiselev8}, which presents the proof of the
theorem:
\\

\noindent {\bf Main theorem } ($ZF$) \\
\quad \\
\hspace*{1cm} {\it There are no weakly inaccessible cardinals}
\\
\\

This proof is derived as a result of using the subinaccessible
cardinal apparatus which the author has developed since 1976. The
idea of this proof had arose in 1984 (though initial approaches to
it have been undertaken by the author since 1973) and in 1996 it
has accepted its present form. Its basic tools -- formula spectra
and subinaccessible cardinals -- were developed
in~\cite{Kiselev1}; reduced spectra and matrices
in~\cite{Kiselev2, Kiselev3, Kiselev4}; disseminators and \
$\delta$-matrices in~\cite{Kiselev4}; matrix informativeness was
developed by the author in~\cite{Kiselev5}; autoexorsizive
matrices in~\cite{Kiselev6}. The proof of main theorem was
presented for the first time in~\cite{Kiselev6} and its more
transparent and complete variant in~\cite{Kiselev7}; the more
systematic presentation of all this field is exposed
in~\cite{Kiselev8}.

However, one can say that this work~\cite{Kiselev8} exposes the
material which is too concise and overloaded by the technical side
of the matter to the detriment of explicit development of the
basic ideas which should advance their technical development.

These circumstances demand the certain preliminary exposition of the field
free of such reproaches.

So, the present work constitutes the first part of the second
edition of~\cite{Kiselev8}, called to overcome them.

The main theorem proof receives in the work its more transparent
and simplified form, though caused by some natural improvements in
final definition 8.2~\cite{Kiselev8} of \ $\alpha$-function, and
zero characteristic is subjected here to the more refined
consideration.

\noindent This circumstance makes possible simplification of the
discussion of \ $\alpha$-function and its properties. As the
result, the final part of the main theorem proof comes in its more
natural purified form. But, on the whole, the technical side of
the matter here repeats~\cite{Kiselev8} but in the more systematic
way.
\\

The special attention is given here to essential preliminary
developments and descriptions of ideas of all main constructions
and reasonings.
\\
It should be pointed out, that the main theorem proof itself occupies
only \S~11 in the second part of this work forthcoming. All
previous sections are devoted to the development of the basic and
special theories of subinaccessibles which, although used for the
creation of ideas and techniques of this proof, but have
independent value.
\\

Here, as before, we make a review of the material~\cite{Kiselev1,
Kiselev2, Kiselev3, Kiselev4, Kiselev5, Kiselev6, Kiselev7} which
is necessary for this proof. We shall omit technically simple
arguments. For set theory, symbols, concepts and other information
necessary for the further reasoning, the remarkable text
Jech~\cite{Jech} provides the basic development of the subject and
much more, so they are assumed to be known and shall be used
frequently without comments.
\\

The plan of this work is the following.
\\
To expose development of the idea of the proof, special notions
and terms are needed.
\\
Just with this end in view the first part of the work, Chapter~I,
is devoted to the basic theory developing, strongly necessary for
introducing the material required for creating the idea.
\\
The very special aim of investigations here makes it possible to
confine the theory to its the most reduced part.
\\
Thus, the existence of some cardinal \ $k$ \ inaccessible in $L$
is assumed and considerations and reasonings are conducted in the
structure $L_k$ or in its generic extensions providing the final
contradiction.
\\

So, in \S~1 the brief outline of the idea of the proof of main
theorem in its first approximation is presented without going in
any details. However, the precise and thorough implementation of
this idea requires the successive involving of new notions and
tools, which are demanded on each stage, since the theory goes on
as series of successive approximations caused by some
insufficiency of the idea in its previous forms. Therefore it
receives its more precise form sequentially involving the
necessary techniques
coming on.\\
Some other way of representing the idea, that is right from the
start in its \textit{final form}, even without details and rigor,
would be opaque and unnatural -- or it would be shallow and vague.
\\
Also \S~1 contains some classical information.
\\

In \S~2 \textit{formula spectra} are introduced (definitions 2.3,
2.4). Briefly saying here, the spectrum of a given formula \
$\varphi$ \ of some level \ $n$ \ is the function containing {\it
all information\/} about its  essential truth properties in {\it
all\/} generic extensions of the structure \ $L_k$ \ of certain
kind, by means of generic ultrafilters on \ $(\omega_0,
k)$-algebra $B$ L\'{e}vy. Such spectrum has range consisting of
\textit{all essential} Boolean values in $B$ of this formula, that
is its \textit{Boolean spectrum}, and domain consisting of
corresponding ordinals, named \textit{jump ordinals} of this
formula, providing these values, that is its \textit{ordinal
spectrum}.
\\

Among spectra the so called \textit{universal spectra} are of
most interest (definition 2.6), contracting \textit{all} formula
ordinal spectra of a given level \ $n$ \ into one
(lemma~\ref{2.7.}).
\\
These notions are considered also in their relativized form; in
this case the bounding ordinal \ $\alpha < k$ \ causes all truth
properties of formula \ $\varphi$ \ in all extensions of \ $L_k$ \
under consideration bounded by \ $\alpha$, \ therefore \ $\alpha$
\ is named the \textit{carrier} of these spectra.
\\

After this is done, with formula spectra in hand, we introduce in
\S~3 the central notion of basic theory -- the notion of
\textit{subinaccessibility} (definition 3.1). The principal
property of subinaccessible ordinal \ $\alpha$ \ of level \ $n$ \
lies in the \textit{impossibility of reaching} this ordinal ``by
means of language'' of level \ $n$. Since formulas truth
properties contains in their spectra, it means that this ordinal \
$\alpha$ \ includes all formula spectra of level \ $n$ \ with
constants \ $ < \alpha$ \ and so can not be reached by such
spectra.
\\
One can see here certain parallel with the notion of
\textit{inaccessibility} itself (see the note after definition
3.1).
\\
Some simple properties of these notions make possible to
introduce the so called \textit{subinaccessibly universal}
formula and its spectra, posessing only subinaccessible
ordinals of smaller level (definition 3.9). In what follows these
spectra provide the basic theory with its favourable tools.\\

In \S~4 the problem of spectrum ``complexity'' springs up. It is
almost evident that this spectrum characteristic must increase
unboundedly while its carrier increases up to \ $k$. \ So the more
significant aspect is investigated for \textit{reduced spectra}
(definition 4.1), received by reducing their Boolean values to
some given cardinal \ $\chi$.
\\
The main result of this section, lemma 4.6 about spectrum type,
shows, that under certain natural and productive conditions the
order type of subinaccessibly universal spectrum reduced to \
$\chi$ \ on a carrier \ $\alpha$ \ exceeds every ordinal \ $<
\chi^{+}$ \ defined below \ $\alpha$ \ (more precisely, defined
below its some jump ordinals).
\\

Still, such spectra on different carriers can be hardly compared
with each other because their domains, ordinal spectra, consist
of ordinals increasing unboundedly, while their carriers are
increasing up to \ $k$.
\\
Therefore in \S~5 we turn to \textit{reduced matrices} (definition
5.1) obtained from reduced subinaccessibly universal spectra by
isomorphic substituting their domains for corresponding ordinals.
The term ``matrix'' is pertinent here because it is possible to
use two-di\-men\-sio\-nal, three-dimensional matrices, etc. for
more fine analysis of the formula truth properties (just on this
way the author conducted his initial investigations of the problem
during long time).
\\
The author have proved in~1977 the main result of this section,
lemma 5.11, which makes possible to see, that such matrix can
contain \textit{all information} about \textit{every} initial part
of constructive universe and, moreover, about its \textit{every}
generic extension, bounded by its carrier (and again, more
precisely, bounded by jump ordinals of corresponding spectrum on
this carrier). So, it preserves this information when passing from
one carrier to another.
\\

Now, as the required notions and terms are formed, the instruments
of the proof of main theorem -- \textit{matrix functions} -- come
into play. Such matrix functions are certain sequences of matrices
reduced to the so called \textit{complete cardinal} \
$\chi^{\ast}$ \ (definition 5.4), the supremum of universal
ordinal spectrum. First the simplest version of such function
(definitions 5.7, 5.14) is exposed at the end of \S~5. After that
it is sharpened by means of supplying its matrices with the so
called \textit{disseminators} (definition 6.1), the special
cardinals extending information about lower levels of universe up
to their carriers (more precisely, up to jump ordinals of
corresponding spectra on these carriers).
\\
Therefore \S~6 is devoted to the investigation of this
disseminator notion in its most poor version in conformity with
the main aim; here we expose some methods of obtaining matrices
and matrix functions supplied by such disseminators.
\\

This stage concludes the developing of the basic theory. Here the
second part of this work forthcoming, Chapters~2 and~3, should be
exposed briefly. With the apparatus in hand, developed in the
first part, the special theory of matrix \
\mbox{$\delta$-}func\-tions starts (definitions 7.1, 7.2),
providing certain tools for correction of necessary technical
sources.

In \S~8 they are developed up to \ $\alpha$-functions,
representing the instrument of the proof in its final form.

In \S\S~9, 10 the properties of \ $\alpha$-functions are
investigated.

This information provides the main theorem proof exposed in
\S~11. After that some easy consequences of this theorem and of
some well-known results are presented in \S~12.

Also one should pay special attention on the comment
in the end of this work.

\newpage
\chapter{Basic Theory: Subinaccessibility, Formula Spectra and Matrices,
Disseminators}

\markboth{\chaptername\ \thechapter. \ Basic Theory}%
{\chaptername\ \thechapter. \ Basic Theory} 


\section{Preliminaries}
\setcounter{equation}{0}

We shall prove that the system
\[
    ZF+\exists k \hspace{2mm} (k \; \mbox{\it is weakly inaccessible
    cardinal})
\]
is inconsistent.
\\
In what follows all the reasoning will be carried out in this
system.
\\

As it was noted above, the idea of the main theorem proof is
rather complicated and untransparent and therefore requires some
special analysis.

So, it is best possible to introduce it into three stages
sequentially, approaching to it more precisely on each stage by
means of developing the corresponding new notions and techniques,
called to overcome certain impeding deficiencies of its previous
version.
\\
Thus it becomes complicated in essence on each stage and receives
its final form in \S~8.

The brief outline of the idea is the following.

First the idea springs up at the end of \S~5 and rests on the
formation of matrix functions that are sequences of matrices,
reduced to a fixed cardinal.
\\
One can obtain sufficiently adequate conception of such function
if beforehand acquaint himself with formulations of the following
notions in outline: the notions of formula spectra (definitions
2.3, 2.4), universal spectrum (definition 2.6, lemma 2.7),
subinaccessibility (definition 3.1), subinaccessibly universal
spectrum (definition 3.9), reduced spectra and matrices
(definitions 4.1, 5.1) and singular matrices (definition 5.7).
\\
On this foundation the simplest matrix function
\[
    S_{\chi f}=( S_{\chi \tau } )_{\tau }
\]
is introduced (definitions 5.13, 5.14) as the sequence of such
matrices of special kind.
\\
This function has range consisting of singular matrices reduced to
some fixed cardinal \ $\chi$ \ and defined as minimal, in the
sense of G\"{o}del function \ $Od$, \ on corresponding carriers;
this property evidently provides its monotonicity, also in the
same sense (lemma 5.17~~1)~). The domain of this function is
extraordinarily long and is cofinal to the inaccessible cardinal \
$k$ \ (lemma \ref{5.18.})
\\
The role of reducing cardinal \ $\chi$ \ is played further mainly
by the complete cardinal \ $\chi^{\ast}$ \ (definition 5.4).

Now the idea of the main theorem proof comes out in its
\textit{initial form}:
\\
\\
\textit{The required contradiction can be attained by means of
creation of some matrix function which should possess inconsistent
properties: it should be monotone and at the same time it should
be deprived its monotonicity.}
\\
\\
However, on this stage the direct proof of this function
nonmonotonicity is impeded by the following phenomenon: the
properties of universe under consideration change after its
bounding by carriers of this function values, that is of reduced
matrices \ $S_{\chi \tau}$ \ (one can see it from the discussion
in the end of \S~5 after lemma 5.18).

Just in order to get over this obstacle the process of matrix
function transformation starts.
\\
The special cardinals, \textit{disseminators}, are introduced.
Though presenting original and valuable phenomenon, this notion
is used here only for the transformation of the matrix function in
view of our aims.
\\
To this end its values, matrices \ $S_{\chi \tau}$, \ are
supplied with certain disseminators extending required universum
properties from lower levels up to their carriers (definitions
6.1, 5.9, 6.9).
\\
After that, in the second part of this work forthcoming, the second
approach to the idea of main theorem proof is undertaken and the matrix
function transforms to \ $\delta$\textit{-function,} also
defined on the set which is cofinal to \ $k$ \ (definitions 7.1,
7.2, lemma 7.6).
\\

But as the result now, vice versa, this new function loses its
monotonicity property (one can see it from the discussion of this
new situation in the end of \S~7 after lemma 7.7).
\\
The way out of this new stage of things lies in the third approach
to the idea of the proof, that is in transformation of this last
function into its more complicated recursive form --
the so called \ $\alpha$\textit{-function}
(definitions 8.1-8.3) which is also defined on the set cofinal to
$k$ (lemma 8.9). This recursive definition is formed in such a
way that cases of monotonicity are demanded in the
\textit{first} turn (so they are provided by \textit{``unit
characteristic''}), while cases of monotonicity violation are
allowed in the last turn only for want of anything better (and
they are of \textit{``zero characteristic''}).
\\
Thus, the priority belongs to matrix function values of unit
characteristic and just cause of that the situations of
monotonicity breaking are avoided.

As the result, \ $\alpha$-function, at last, delivers the required
contradiction: it cannot be monotone (theorem 1) and at the same
time it possesses monotonicity (theorem 2).
\\

Such is the idea of the main theorem proof in its outline form;
more detailed description draws into consideration too much
technical details. Therefore instead of that one should follow
more suitable way of assimilating notions mentioned above (up to
the end of \S~6) without going into details and getting along, may
be, with the help of some images. Using another approach (the
representation directly in its final form right from the start)
one should receive the exposition even more complicated than the
foregoing one, or noninformative and hazy.
\\

After that let us turn to realization of the programme sketched
above.
\\

Weakly inaccessible cardinals become strongly inaccessible in
G\"{o}del constructive class \ $L$; \ let us remind that it is the
class of values of G\"{o}del constructive function \ $F$ \ defined
on the class of all ordinals. Every set \ $a\in L$ \ receives its
ordinal number
\[
    Od ( a ) =\min \{ \alpha :F (\alpha  ) =a\}.
\]
If \ $\alpha $ \ is an ordinal then \ $L_{\alpha }$ \ denotes the
initial segment
\[
    \{ a\in L:Od ( a )<\alpha \}
\]
of this class. The starting structure in the further reasoning is
the countable initial segment
\[
    \mathfrak{M}=(L_{\chi ^{0}},\; \in , \; =)
\]
of the class \ $L $ \ serving as the standard model of the theory
\[
    ZF+V=L+\exists k \hspace{2mm} (k \; \mbox{\it is weakly
    inaccessible cardinal})
\]
Actually, only the finite part of this system will be used here
because we shall consider only formulas of limited length, as it
will be clear from what follows. Moreover, the countability of
this structure is required only for some technical convenience
(see below) and it is possible to get along without it.

Further \ $k$ \ is the \textit{smallest inaccessible cardinal in}
\ $\mathfrak{M}$. We shall investigate it ``from inside'',
considering the hierarchy of subinaccessible cardinals; the latter
are ``inaccessible'' by means of formulas of certain elementary
language. To receive this hierarchy rich enough it is natural to
use some rich truth algebra \ $B$. It is well-known
(Kripke~\cite{Kripke}) that every Boolean algebra is embedded in
an appropriate collapsing \ $ ( \omega _{0},\mu  ) $-algebra and
therefore it is natural to use as \ $B$ \ the sum of the set of
such algebras of power \ $k$ \ , that is L\'{e}vy \ $ ( \omega
_{0},k ) $-algebra \ $B$.

Namely, let us apply the set \ $P\in \mathfrak{M}$ \ of forcing
conditions that are finite functions \ $p\subset k\times k$ \ such
that for every limit \ $\alpha <k$ \ and \ $n\in \omega_{0}$
\vspace{-6pt}
\[
    \alpha +n\in dom ( p ) \longrightarrow p (\alpha +n ) <\alpha ;
\]
also let \ $p(n)\leq n $ \ for \ $\alpha=0 $. The relation \ $\leq
$ \ of partial order is introduced on \ $P$: \vspace{-6pt}
\[
    p_{1} \leq p_{2} \longleftrightarrow p_{2} \subseteq p_{1} .
\]
After that \ $P$ \ is embedded densely in the Boolean algebra \
$B\in \mathfrak{M} $, consisting of regular sections \ $\subseteq
P$, which is complete in \ $\mathfrak{M}$. \ The relation of
partial order \ $\leq $ \ is defined on \ $B$:
\[
    A_{1} \leq A_{2} \longleftrightarrow A_{1} \subseteq A_{2} ,
\]
and also Boolean operations \ $\cdot$, $+$, $\prod$, $\sum $ \
(see~\cite{Jech}) are defined on \ $B$. \ Every condition \
\mbox{$p\in P$} \ is identified with the section
\[
    \left[ p\right] =\{ p_{1}\in P:p_{1}\leq p\}
\]
and that is why \ $P$ \ is isomorphically embedded in \ $B$. \
Hereafter we recall the well-known results of Cohen~\cite{Cohen},
L\'{e}vy~\cite{Levy} (see also Jech~\cite{Jech}).

\begin{lemma} \label{1.1.} \quad \\
\hspace*{1em} Algebra \ $B$ \ satisfies the \ $k$-chain condition,
that is every set \  $X\subseteq B$, \ $X\in \mathfrak{M}$  \
consisting of pairwise disjoint Boolean values has the power \
$<k$ \ in \ $\mathfrak{M\;}$:
\[
    \left ( \forall A_{1},A_{2}\in X\quad A_{1}\cdot A_{2}=0 \right)
    \longrightarrow \left| X\right| <k.
\]
\end{lemma}
According to this lemma it is possible to consider instead of
values \  $A\in B$ \ only sets
\[
    P_{A}=\{ p\in A:dom (p) \subseteq \chi \}
\]
where
    \[\chi =\min \{ \chi^{\prime }:\forall p\in A\quad p\left|
    \chi ^{\prime }\right. \leq A\}
\]
(here \ $p\left| \chi ^{\prime }\right.$ \ is the restriction of \
$p$ \  to \  $\chi ^{\prime }$). Since \  $A=\sum P_{A}$, \  we
shall always identify \  $A$ \ and \  $P_{A}$, \ that is we shall
always consider \ $P_{A}$ \ instead of \ $A$ \ itself.
\\
Just due to this convention all Boolean values \ $A \in B$ \ are
sets in \ $L_k$, \ not classes, and this phenomenon will make
possible all further reasoning as a whole.
\\

We shall investigate the hierarchy of subinaccessible cardinals
with the help of Boolean values in \  $B$ \  of some propositions
about their properties, that is working inside the Boolean-valued
universe \  $\mathfrak{M}^{B}$. The countability of the structure
\ $\mathfrak{M}$ \  is needed here only to shorten the reasoning
when using its generic extensions by means of \
$\mathfrak{M}$-generic ultrafilters on \ $B$. \ It is possible to
get along without it developing the corresponding reasoning in the
Boolean-valued universe \ $L^{B}$ \ (see, for
example,~\cite{Jech}).

It will be more suitable to produce generic extensions of \
$\mathfrak{M}$ \
 not by means of ultrafilters but by means of functions. Namely,
as \  \mbox{$\mathfrak{M}$-generic} or L\'{e}vy function on \ $k $
\ we shall name every function \  $l\in {}^{k}k$ \  such that every set \
$X \in \mathfrak{M} , X\subseteq B$, which is dense in \ $P$, contains
some  \ $p\subset l$. All functions of this kind will be denoted
by the common symbol \  $l$. \ Obviously, \ $\mathfrak{M}$-generic
ultrafilters \  $G$ \  on \  $B$ \  and these functions mutually
define each other:
\[
    \quad l=\cup ( P\cap G) \quad , \quad G= \{ A\in B:\exists
    p \in P( p\leq A\wedge p\subset l) \}.
\]
In this case the interpretation \  $i_{G}$ \  of the universe \
$\mathfrak{M}^{B}$ \  is denoted  by \  $i_{l}$.  \  As usual, if
\ $\underline{a}\in \mathfrak{M}^{B}$, $a\in \mathfrak{M}\left[
l\right] $ \  and \  $i_{l}(\underline{a}) =a$, then \
$\underline{a}$ \ is named the label or the name of \ $a$. \
By \ $\left\| \varphi \right\| $
\  as usual is denoted the Boolean value of the proposition \
$\varphi $ \  with constants from \ $\mathfrak{M}^{B}$  \  in
algebra \ $B$.
\\
For some convenience we introduce the relation \
$\stackrel{\ast}{\in }~$: \ for every \ $l\in {}^{k}k$, \ $A\in B$
\[
    l\stackrel{\ast }{\in }A\longleftrightarrow \exists p\in P
    ( p\subset l\wedge p\leq A ).
\]

\begin{lemma} \label{1.2.} \hfill {} \\
\hspace*{1em} Let \ $l$ \ be an \ $\mathfrak{M}$-generic function
on \ $k$ \ and \ $\varphi  ( a_{1}...,a_{n} ) $ \ be a proposition
containing constants \ $a_{1}...,a_{n}\in \mathfrak{M}\left[ l
\right] $ \ with names \
$\underline{a}_{1}...,\underline{a}_{n}$, \ then
\[
    \mathfrak{M}[l] \vDash \varphi ( a_{1}...,a_{n} )
    \longleftrightarrow l\ \stackrel{\ast }{\in }\left\| \varphi
    ( \underline{a}_{1}...,\underline{a}_{n} ) \right\| .
\]
\end{lemma}

\vspace{-6pt}

\begin{lemma} \label{1.3.} \hfill {} \\
\hspace*{1em} Let \ $l$ \ be an \ $\mathfrak{M}$-generic function
on $k$, \ then:\quad \\ \quad
\\ {\em 1)} \quad \quad \quad \quad \quad \quad $
\mathfrak{M}\left[ l\right]
\vDash ZF+V=L\left[ l\right] +GCH+k=\omega _{1} ; $ \\ \quad \\
{\em 2)} \quad for every \ $\chi<k$ \ let \ $\chi_{1}=\chi$ \ iff
\ $\chi$ \ is regular and \ $\chi_{1}=\left( \chi^{+} \right)^
\mathfrak{M} $ \ iff \ $\chi$ \ is singular cardinal in \
$\mathfrak{M}$, then
\\
\\
\[
    \mathfrak{M}\left[ l|\chi \right]\models \forall \alpha < \chi_{1}
    \ \ |\alpha |\leq \omega_{0} \wedge \forall \alpha \geq
    \chi_{1}|\alpha |=|\alpha |^\mathfrak{M}\ .
\]

\end{lemma}

\vspace{-0pt}

\begin{lemma} \label{1.4.} \hfill {} \\
\hspace*{1em} Suppose that
\[
    t\in \mathfrak{M}^{B}, \: dom(t)\subseteq \left\{ \check{a} :
    a\in \mathfrak{M} \right\} \:\: \mbox{\it and} \:\: |rng(t)|<k.
\]
Let \ $B_{t}$ \ be the subalgebra of \ $B$ \ generated by \
$rng(t)$. \\
Then for every formula \ $\varphi (t)$
\[
    \parallel \varphi (t) \parallel \in B_{t}.
\]
\end{lemma}

Such is the preliminary information required to start the
developing the theory of subinaccessibles. \label{c1}
\endnote{
\ p. \pageref{c1}. \ Here the light should be shed
upon the following points.

First, one can ask, whether in this investigation inaccessible
cardinal hypothesis is used in the essential way.
\\
This hypothesis is used through all the text. The Power Set Axiom
is used in \ $L_k$ \ at almost every step; in particular, when
working with Skolem functions while formulas and their spectra
are transformed, and in every case providing some new notion.
\\
Besides, the inaccessibility of the cardinal \ $k$ \ is used
strongly in the $k$-chain property of algebra \ $B$ \ (lemma~1.1),
without which Boolean values, spectra and matrices are classes
in \ $L_k$, not sets, and they can not be compared with each
other in a required way. As a result basic theory and special
theory both become impossible in this situation.

One can ask also, whether the entire argument can be carried out
strictly in $ZFC$, with inaccessible $k$ replaced by the proper
class ordinal \ $k$ \ in $ZFC$ (strictly, of course, meaning this
proper class in Neumann-G\"{o}del-Bernays Set theory). One can
imagine that there is no problem with introducing the L\'{e}vy
forcing in this context, and, so, that it provides the
inconsistency of $ZF$.
\\
However, the direct transference of the whole reasoning from \
$L_k$ \ onto \ $L$ \ is considerably hampered. The point is that
several important parts of the proof are conducted by using
generic extensions of \ $\mathfrak{M}$ \ produced by  \
$\mathfrak{M}$-generic ultrafilters (or functions) over \
$(\omega_0,k)$-algebra Levy \ $B$. This algebra has the properties
significant for the reasoning, the \ $k$-chain property and
others, provided by the inaccessibility of \ $k$; in there turn
they provide the basic properties of spectra, subinnaccessibles,
reduced spectra, matrices, matrix functions, etc. It is
important that in this extensions every infinite ordinal \ $<k$ \
is countable due to first component \ $\omega_0$. The cardinal \
$\omega_0$ \ can be replaced here by some cardinal \ $\lambda < k$
\ ($k$-chain property will remain, etc.). But assume that Levy
forcing as a proper class forcing can be used in {\it this
context}, that is \ $k$ \ is replaced by the class of all
ordinals, -- then it provides the generic extension of \ $L$ \ in
which {\it every infinite cardinal is still countable}. It is
possible to get along without these extensions and work within
Boolean-valued universe $L^B$ but still the analogous obstacle
should arise -- there is the value
\[
    \| \mathrm{every\ infinite\ ordinal\ is\ countable} \|=1
\]
(or \ $\|$every ordinal has power $\le \lambda \|=1$ if we use \
$\lambda$ \ instead of \ $\omega_0$) in \ $B$. The author intends
to clarify this aspect more later.
\\
\quad \\
} 

\newpage

\section{Formula Spectra}
\setcounter{equation}{0}

The main instrument of the further reasoning is the notion of a
formula spectrum. In this section the basic spectrum theory is
exposed, containing the discussion of the most simple spectrum
properties. It is possible to introduce this notion in a more
general version (for an arbitrary Boolean algebra \  $B$ \  and a
partially ordered structure \ $\mathfrak{M\ }$); however, here it
is enough to use its the most poor variant.
\\

Among names from \  $\mathfrak{M}^{B}$ \  the canonical names are
distinguished, that are those which give analogous results under
any interpretation. For example, such are names \ $\check{a}$ \
of  sets \ $ a\in \mathfrak{M}$ \ which we shall identify with
these sets:
\[
    dom ( \check{a} ) =\{ \check{b}:b\in a\}, \quad
    rng ( \check{a} ) =\{ 1\}.
\]
The canonical name of every \ $\mathfrak{M}$-generic function \
$l$ \ on \  $k$ \ is the function
\[
    \underline{l}=\{ ( ( \alpha, \beta ) ,\{ ( \alpha ,\beta ) \}):
    \{( \alpha ,\beta ) \} \in P\}.
\]
It is easy to see that always \quad  $i_{l}( \check{a}) = a$, \
$i_{l}( \underline{l}) =l$.

Let us introduce the following elementary language \
$\mathcal{L}$ \ over the standard structure
\[
    ( L_{k}\left[ l\right] ,\; \in, \; =,\; l ).
\]
Its alphabet consists of usual logic symbols : quantors \ $\forall
$, $\exists $, \  connectives \  $ \wedge $, $\vee $, $\neg $,
$\longrightarrow $, $\longleftrightarrow $ , \ brackets $( \; , \;
)$, individual variables \ $x,y,z$... \ (with indices or without
them), all names from Boolean-valued universe \ $L_{k}^{B} $ \
serving as individual constants, and symbols \ $\in $,
$=$, $\underline{l}$.\\
When interpreting this language in the generic extension \
$L_{k}\left[ l\right] $ \ variables run through \ $L_{k}\left[
l\right] $, \ individual constants \  $a\in L_{k}^{B}$ \  denote
\  $i_{l} ( a ) $ \  and constants \  $\in $, $=$,
$\underline{l}$ \  denote respectively the standard relations of
membership, equality and function \  $l$. \  If this extension is
fixed, then constants and their interpretations will be
identified as usual.

The formulas of language \  $\mathcal{L}$ \  are defined in a
usual recursive way beginning with atomic formulas of the form \ $
t_{1}=t_{2}$, $\quad t_{1}\in t_{2}$, \ where \ $t_{1}$, $t_{2}$ \
are any terms that are also recursively formed of variables,
constants \  $\in L_{k}^{B}$ \ and \ $\underline{l}$ \ by
sequential superposition. However, following the tradition, some
conventional notation, relations and terms will be used in
writings of formulas if it will not cause difficulties. For
example, the train of \ $m$ \ variables or constants \
$x_{1},...,x_{m}$  \ is denoted by \ $( x_{1},...,x_{m} ) $ or, in
short, by \ $\overrightarrow{a}$; the order relation \ $ x_{1}\in
x_{2} \  $on the set of ordinals is denoted by \ $x_{1}<x_{2}$  \
and so on.

Further formulas will be considered as formulas of the language
\  $\mathcal{L}$ \  which we shall denote by small letters from
the end of Greek alphabet (unless otherwise specified).

Formula \ $\varphi $ \ which has free variables and individual
constants forming a train \  $ \overrightarrow{a}$ \ will be
denoted by \ $\varphi  ( \overrightarrow{a} ) $. \  If in
addition \  $\varphi $ \ contains a symbol to which it is
necessary to pay attention it should also be noted specifically;
for example, the notation \ $\varphi (
\overrightarrow{a},\underline{l} ) $ \  points out that \  $
\varphi $ \  contains an occurrence of \ $\underline{l}$.

As usual, occurrences of quantors \ $\exists x,\forall x$ \ in
formula \ $\varphi$ \ are named bounded by a term \ $t$ \
{\it iff} they have the form
\[
    \exists x~ ( x\in t\wedge ... )
    ,\forall x~ ( x\in t\longrightarrow ... );
\]
a formula is named bounded
{\it iff} all its quantors \ are bounded by some terms; it is
named prenex formula {\it iff} all occurrences of its \textit{unbounded}
quantors are disposed to the left from occurrences of other quantors and
connectives; this train of its unbounded quantors is named its
\textit{quantor prefix}.
\\
Formulas of the language \  $\mathcal{L}$  \  will be interpreted
in generic extensions  \  $L_{k}[l] $ \  and therefore we shall
name formulas \ $\varphi,\psi $ \  \textit{equivalent} and write \
$\varphi \longleftrightarrow \psi $ \ \textit{iff} they are
equivalent in $\mathrm{ZFC}^{-}+V=L[$\underline{l}$\ ]$, where
$\mathrm{ZFC}^{-}$   is ZFC with the Power Set Axiom deleted.
\\
Besides that, when interpreting formulas  \  $ \varphi, \psi $ \
in  \  $L_{k}$ \  we shall name them \textit{constructively
equivalent} and use the same notation  \  $\varphi
\longleftrightarrow \psi $ \ \textit{iff} they are equivalent in
$\mathrm{ZFC}+V=L$.
\\
We shall name them \textit{generically equivalent} and write \
$\varphi \dashv \vdash \psi $, \ \textit{iff}
\[
    \| \varphi \longleftrightarrow \psi \| =1
\]
for any values from \ $L_{k}^{B}$  \  of their free variables.
\\
Further the equivalence in  $\mathrm{ZFC}^{-}+V=L[$\underline{l}$\
]$  will be considered, if the context will not point to some
another case; at any rate, the meaning of the notion of
equivalence and of the symbol \ ``$\longleftrightarrow$'' \ will
be always obviously specified by the context.

\begin{definition}
\label{2.1.} \ \\
{\em 1)}\quad The class of all formulas \  $\varphi
(\overrightarrow{a},\underline{l}) $ \  in the prenex form that
have quantor prefix consisting of \  $n$ \  maximal blocks of
like quantors and begin with \  $\exists $ \  and also of all
formulas that are equivalent to such formulas is denoted by
$\Sigma _{n}( \overrightarrow{a}) $.
\\
The dual class is denoted by \  $\Pi _{n}( \overrightarrow{a}) $
\  and the class  \ $\Sigma _{n}( \overrightarrow{a} ) \cap \Pi
_{n} ( \overrightarrow{a})$ \ -- by \ $\triangle_{n} \left(
\vec{a}\right)$.
\\ As a result the elementary Levy hierarchy comes out:
\[
    \{ \Sigma _{n}( \overrightarrow{a} ) ;\ \Pi _{n}
    ( \overrightarrow{a}) \} _{n\in \omega _{0}} .
\]
{\em 2)}\quad  We denote by \  $\Sigma _{n}^{\dashv \vdash }(
\overrightarrow{a}) $ \  the class of all formulas generically
equivalent to formulas from \  $\Sigma _{n}(
\overrightarrow{a})$.
\\
The dual class is denoted by \ $\Pi_{n}^{\dashv \vdash }(
\overrightarrow{a}) $ \ and the class \ $\Sigma _{n}^{\dashv
\vdash } (\overrightarrow{a}) \cap \Pi _{n}^{\dashv \vdash }
(\overrightarrow{a}) $ \ -- \ by \  $\Delta _{n}^{\dashv \vdash
}( \overrightarrow{a}) $.
\\
\quad \\
The symbols \  $Q_{n} (\overrightarrow{a})$, $ Q_{n}^{\dashv
\vdash} (\overrightarrow{a}) $ \ serve as the common notation
respectively of the classes
\[
    \Sigma _{n}( \overrightarrow{a}), ~
    \Pi _{n}( \overrightarrow{a}), ~~~
    \Sigma _{n}^{\dashv \vdash}( \overrightarrow{a}), ~
    \Pi _{n}^{\dashv \vdash }(\overrightarrow{a}).
\]
The natural index  \ $n$ \ is named their level and the level of
their formulas.
\\
Considering fixed extension \ $L_{k}[ l] $, \ sets, that are
defined in it by formulas from \  $Q_{n}^{\dashv \vdash }
(\overrightarrow{a}) $, will be named \  $Q_{n}^{\dashv
\vdash }( \overrightarrow{a}) $-sets.
\\
As a result the hierarchy
\[
    \{ \Sigma _{n}^{\dashv \vdash } ( \overrightarrow{a} ) ;\
    \Pi_{n}^{\dashv \vdash } ( \overrightarrow{a} ) \}_{n \in \omega_0}
\]
and the corresponding hierarchy of relations come out.\\
In all the notation the train \  $\overrightarrow{a}$ \  will be
omitted if its value is arbitrary or meant by the context.
\end{definition}

\vspace{12pt}

Following the tradition, let us assume that every formula under
consideration is considered as transformed to the equivalent
prenex form of the \textit{minimal level} \ --- \ equivalent in
$ZFC^{-}$ or constructibly or generically depending on structures
where these formulas are interpreted \ --- \ of course, if the
context does not mean some another situation.
\\
Let us also assume that if in the context formulas interpreted in
the structure \  $ ( L_{k},\; \in ,\; = ) $ \ are considered,
then we use all these terms and notation, but without the index \
$\dashv \vdash $.
\\

Further the classes \  $Q_{n}$, $Q_{n}^{\dashv \vdash }$ \  of
some fixed level \  $n>3$ \  are considered (if the context does
not imply the opposite). This agreement is taken to have in hand
further sufficiently large subinaccessible tools (see lemma 3.5
below for example) and also to use some auxiliary formulas, terms,
relations and sets defined in \  $L_{k}$ \  directly as additional
constants in formulas notation without raising their level.
Obviously, in this way can be considered the sets \ $P$, \ $B$,
relations and operations on them mentioned above, and also the
following:
\\
\quad \\
1)\quad $On ( x ) $ \  - the  bounded formula meaning that \  $x$
\ is an ordinal:
\[
    \quad \forall y\in x ( \forall z\in x ( y\in
    z\vee z\in y ) \wedge \forall z_{1}\in y ( z_{1}\in x )  ) .
\]
Variables and constants bounded by this formula will be also
denoted by small letters from the beginning of Greek alphabet: \
$\alpha $, $\beta $, $\gamma $, ... \  omitting this formula
itself. With the help of this formula it is easy to define all
the natural numbers, the ordinals  \  $\omega_{0}, \ \
\omega_{0}+1 .... $ \  and so on by the corresponding bounded
formulas, so we shall use these
ordinals as the additional individual constants of the language.\\
2)\quad $F ( x,y ) $ \  -- the  \  $\Delta _{1}$-formula providing
the well-known recursive definition in \  $L_{k}$ \  of G\"{o}del
function \  $F$~~\cite{Godel}, that is for every \ $\alpha \in k,
\: a\in L_{k} \: :$
\[
    a=F\mid  ( \alpha +1 ) \longleftrightarrow L_{k}\vDash F
    (\alpha, a).
\]
3)\quad $\lessdot $, \ $\underline{\lessdot }$ \quad -- the
relations of well-ordering on\quad $L_{k}$:
$$ \quad a\lessdot b\longleftrightarrow Od ( a ) <Od ( b ) \quad ;\quad a\mathbf{\
}\underline{\lessdot }b\longleftrightarrow a\lessdot b\vee a=b .
$$
4)\quad $\vartriangleleft $, \ $\trianglelefteq $ \  -- the
corresponding relations on \ $L_{k}\times k$:
\[
    a \vartriangleleft \beta \longleftrightarrow Od ( a ) < \beta
    \wedge On ( \beta ) \quad ; \quad a \trianglelefteq \beta
    \longleftrightarrow a\vartriangleleft \beta \vee Od ( a ) = \beta.
\]
5) \quad It is not hard to use similarly the analogous  \  $
\Delta _{1} $ -formula  \  $F ( x,y,\underline{l} )$  \ providing
the recursive definition in \  $L_{k}\left[ l\right] $ \  of
G\"{o}del constructive function \  $F^{l}$ \   relalively to \
$l$ \  and to receive the function of ordinal number \ $Od^{l}$
\  and also the relations \  $\lessdot^{l}$, $\vartriangleleft
^{l}$:

\[
    ~Od^{l} ( a ) =\min \{ \alpha :F^{l} ( \alpha ) =a\} \ \ ;
\]
\[
    a \lessdot ^{l}b\longleftrightarrow Od^{l} ( a ) < Od^{l} ( b )~;
    \quad \quad\quad a\underline{\lessdot }^{l} b\longleftrightarrow
    a \lessdot ^{l}b\vee a=b ~;
\]
\[
    a\vartriangleleft ^{l}\beta\longleftrightarrow Od^{l} ( a )
    <\beta\wedge On ( \beta ) ~; \quad a\trianglelefteq^{l}\beta
    \longleftrightarrow a\vartriangleleft ^{l}\beta\vee
    Od^{l} ( a )=\beta ~.
\]

\vspace{12pt}

It is easy to define all these   functions and relations mentioned
above by \ $\Delta _{1}$-formulas respectively in \ $L_{k}$,
$L_{k}\left[ l\right] $ \ absolutely relatively to these
structures using formulas \ $F( x,y )$, $F(x,y,\underline{l} )$.
We shall denote these formulas by the same symbols as functions
and relations defined by them, but replacing the function \  $l$ \
by its name \  $\underline{l}$. \  The relation \  $\lessdot $ \
was used by Addison~\cite{Addison} over continuum and by
Kogalovskiy~\cite{Kogalovski} over arbitrary infinite structures
of arbitrary levels.
\\
With the help of ordinary Skolem functions techniques it is easy
to prove

\begin{lemma} \label{2.2.} \hfill {} \\
\hspace*{1em} Suppose a formula \  $\varphi $ \  contains in the
class \ $ Q_{n}$, then the formulas
\[
    \forall x\in y~\varphi, \quad \exists x\in y~\varphi
\]
are contained in the same class.
\\
Analogously for the class  \ $Q_{n}^{\dashv \vdash }$, \  $n\geq
1$, \ replacing the bounding formula \  $ x\in y$ \  by the
bounding formulas \  $x\lessdot ^{\underline{l}}y$, \
$x\vartriangleleft ^{\underline{l}}y$.
\end{lemma}

Let us turn to the notion of spectrum. In order to make it more
transparent we shall introduce it only for propositions \ $\varphi
 ( \overrightarrow{a},\underline{l} ) $ \  having the train of
individual constants \ $\overrightarrow{a}= ( a_{1},...a_{m} )$ \
consisting of ordinal constants (if the context does not point to
another case). It is possible to manage without this convention
replacing occurrences of each \  $a_{i}$ \ by occurrences of the
term \ $F^{\underline{l}} ( \alpha _{i} ) $ \ for the
corresponding ordinal constant \  $\alpha _{i}$.
\\
Let us also assume that every train \  $\overrightarrow{a}= (
\alpha _{1},...,\alpha _{m} ) $ \  of ordinals \ $ < k$ \ is
identified with the ordinal which is its image under the canonical
order isomorphism of \  $ ^{m}k$ \  onto \  $k$. \  The
isomorphism of \  $^{2}k$ \ onto \  $k$  \  of this kind will be
named the \textit{pair
function}.
\\
The next notion plays an important role in what follows:

\begin{definition} \label{2.3.} \hfill {} \\
\hspace*{1em} For every formula \ $\varphi $ \ and ordinal \
$\alpha _{1}\leq k$ \ by \ $\varphi ^{\vartriangleleft \alpha
_{1}}$  \ is denoted the formula obtained from \  $\varphi$ \ by \
$\vartriangleleft ^{\underline{l}}$~-bounding all its quantors by
the ordinal \ $\alpha _{1}$, \  that is by replacing all
occurrences of such quantors \  $\exists x$, $\forall x$ \ by the
corresponding occurrences of
\[
    \exists x~ ( x\vartriangleleft ^{\underline{l}} \alpha _{1}
    \wedge ... ) , ~~~ \forall x~ ( x\vartriangleleft ^{\underline{l}}
    \alpha _{1}\longrightarrow ... ).
\]
In addition, if \  $\alpha _{1}<k$, \  then we say that \ $\varphi
$ \  is restricted to \  $\alpha _{1}$ \  or relativized to \
$\alpha _{1}$; \  if, in addition, the proposition \ $\varphi
^{\vartriangleleft \alpha _{1}}$ \  holds, then we say that \
$\varphi $  \  holds below \  $\alpha _{1}$ \  or that \ $\varphi
$ \  is preserved under restriction or relativization to \ $\alpha
_{1}$.

The same terminology is carried over to all reasoning and
constructions having variables and individual constants \
$\vartriangleleft ^{\underline{l}}$- bounded by the ordinal \
$\alpha _{1}$.
\\
We shall consider  \  $\vartriangleleft$-restriction instead of
 \  $\vartriangleleft ^{\underline{l}} $~-restriction in all these
notation and notions iff formulas, constructions and reasoning
are interpreted in  \  $L_{k}$.
\\
In all such cases  \ $\alpha_{1}$ \ is named respectively the \
$\vartriangleleft ^{\underline{l}}$~-bounding or \
$\vartriangleleft $~-bounding ordinal.
\\
If \  $\alpha _{1}=k$, \  then the upper index \ $\vartriangleleft
\alpha _{1}$ \  is omitted and such formulas, reasoning and
constructions are named unrestricted or unrelativized.
\end{definition}

\vspace{-6pt}

\begin{definition}
\label{2.4.}\ \\
{\em 1)}\quad Let \  $\varphi (\overrightarrow{a},\underline{l} )
$ \  be a proposition \  $\exists x~\varphi _{1}
(x,\overrightarrow{a}, \underline{l} ) $ and \ $\alpha _{1}\leq
k$. For every \ $\alpha <\alpha _{1}$ \  let us introduce the
following Boolean values:

\vspace{12pt}
\[
    A_{\varphi }^{\vartriangleleft \alpha _{1}} ( \alpha ,
    \overrightarrow{a}) = \left\| \exists x\trianglelefteq ^{\underline{l}}
    \alpha ~\varphi_{1}^{\vartriangleleft \alpha _{1}}
    ( x,\overrightarrow{a},\underline{l} ) \right\| ;
\]

\[
    \Delta _{\varphi }^{\vartriangleleft \alpha _{1}} ( \alpha ,
    \overrightarrow{a} ) =A_{\varphi }^{\vartriangleleft \alpha _{1}}
    ( \alpha ,\overrightarrow{a} ) - \sum_{\alpha ^{\prime }<\alpha}
    A_{\varphi }^{\vartriangleleft \alpha _{1}} ( \alpha ^{\prime },
    \overrightarrow{a} ) .
\]
\vspace{12pt}

\noindent {\em 2)}\quad We name the following function \
$\mathbf{S}_{\varphi }^{\vartriangleleft \alpha _{1}} (
\overrightarrow{a} ) $ the spectrum of \  $\varphi $ \  on the
point \  $\overrightarrow{a}$ \  below \  $\alpha _{1}$:
\vspace{12pt}
\[
    \mathbf{S}_{\varphi }^{\vartriangleleft \alpha _{1}}
    ( \overrightarrow{a}) =\{  ( \alpha ,
    \Delta_{\varphi }^{\vartriangleleft \alpha_{1}}
    ( \alpha ,\overrightarrow{a} )  ) :\alpha <\alpha_{1}
    \wedge \Delta _{\varphi }^{\vartriangleleft \alpha _{1}}
    ( \alpha , \overrightarrow{a} ) >0\}.
\]
\vspace{0pt}

\noindent Projections
\[
    dom \left( \mathbf{S}_{\varphi }^{\vartriangleleft \alpha _{1}}
    ( \overrightarrow{a} ) \right), \quad  rng \left (
    \mathbf{S}_{\varphi }^{\vartriangleleft \alpha _{1}}
    ( \overrightarrow{a} ) \right )
\]
\vspace{-3pt}

\noindent are named respectively the ordinal and the Boolean spectra of \  $\varphi $
\ on the point \ $\overrightarrow{a}$ \  below \ $\alpha _{1}$.

\noindent {\em 3)}\quad  If \ $ ( \alpha ,\Delta  ) \in
\mathbf{S}_{\varphi }^{\vartriangleleft \alpha _{1}} (
\overrightarrow{a} ) $, \ then \  $\alpha $ \  is named
the jump ordinal of this
formula and spectra, while \ $\Delta $ \  is named
its Boolean value on the
point \  $\overrightarrow{a}$ \ below \ $\alpha _{1}$.
\\
\noindent {\em 4)}\quad The ordinal \ $\alpha _{1}$ \ itself is
named the carrier of these spectra.
\end{definition}
If a train \  $\overrightarrow{a}$ \  is empty, then we omit it
in notations and omit other mentionings about it.
\\

To analyze propositions it is natural to use their spectra. We
can develop more fine analysis using their two-dimensional,
three-dimensional spectra and so on. \label{c2}
\endnote{
\ p. \pageref{c2}. \ The two-dimensional spectrum can be
constructed in the following way. Let us consider the proposition
\[
    \varphi  ( \overrightarrow{a},\underline{l} ) =\exists x\ \varphi_{1}
    ( x,\overrightarrow{a},\underline{l} ) =\exists x\forall y\
    \varphi_{2} ( x,y,\overrightarrow{a},\underline{l} ) .
\]

To each jump ordinal \ $\alpha $ \ of \ $\varphi $\  on the point
\ $\overrightarrow{a}$\  below \ $\alpha _{1}$\  we adjoin the
spectrum of the proposition \  $\varphi _{1}^{\prime }=\exists
y\neg \varphi _{2} ( F^{\underline{l}} ( \alpha  ) ,y,
\overrightarrow{a},\underline{l} ) $ \  but on the point \ $ (
\alpha ,\overrightarrow{a} ) $.\  As a result we obtain the
function that is two-dimensional spectrum:
\\
\vspace{6pt}
\[
    \mathbf{S}_{\varphi ,2}^{\vartriangleleft \alpha _{1}}
    ( \overrightarrow{a} ) = \{  ( \alpha ,\beta ,
    \Delta _{\varphi}^{\vartriangleleft \alpha _{1}}
    ( \alpha ,\overrightarrow{a} ) ,
    \Delta _{\varphi _{1}^{\prime }}^{\vartriangleleft \alpha _{1}}
    (\alpha ,\beta ,\overrightarrow{a} )  ) :\alpha ,
    \beta<\alpha_{1}\wedge \\
\]
\[
    \qquad\qquad \qquad\qquad\wedge
    \Delta _{\varphi}^{\vartriangleleft \alpha _{1}}
    ( \alpha ,\overrightarrow{a} )>0\wedge
    \Delta _{\varphi _{1}^{\prime }}^{\vartriangleleft \alpha _{1}}
    ( \alpha ,\beta ,\overrightarrow{a} ) >0 \}.
\]
\vspace{6pt}

\noindent Its first and  second projection constitutes the
two-dimensional ordinal spectrum, while third and fourth
projection constitutes the two-dimensional Boolean spectrum of \
$\varphi $ \  on the point \  $ \overrightarrow{a}$ \  below \
$\alpha _{1}$. \  The three-dimensional spectrum are introduced
in a similar way by consequent adjoining the spectra of \
$\varphi _{2}$ \  and so on. However the investigation of
multi-dimensional spectra lies outside the limits of this work.
\\
\quad \\
} 

All spectra introduced possess the following simple properties:

\begin{lemma} \label{2.5.} \quad \\
\hspace*{1em} Let \ $\varphi $ \ be a proposition
\[
    \exists x\ \varphi _{1} ( x,\overrightarrow{a},l ), \quad
    \quad \varphi_{1}\in \Pi _{n-1}^{\dashv \vdash }, \quad
    \alpha _{1}\leq k,
\]
then:
\\
{\em 1)}\quad $\sup dom \left (
\mathbf{S}_{\varphi}^{\vartriangleleft \alpha_{1}} (
\overrightarrow{a} ) \right ) <k $ ;
\hfill {} \\
\hfill {} \\
{\em 2)}\quad $\mathbf{S}_{\varphi }^{\vartriangleleft \alpha
_{1}} ( \overrightarrow{a} ) $, \ $dom \left (
\mathbf{S}_{\varphi }^{\vartriangleleft \alpha _{1}} (
\overrightarrow{a} ) \right) $ \ are \ $\Delta _{n}$-definable,
while
\\
$rng \left (\mathbf{S} _{\varphi }^{\vartriangleleft
\alpha _{1}} ( \overrightarrow{a} ) \right )$ \ is \ $\Sigma
_{n}$-definable in \ $L_{k}$ \ for \ $\alpha _{1}=k$.
\\
For \ $\alpha _{1}<k $ \ all these spectra are \
$\Delta_{1}$-definable;
\hfill {} \\
\hfill {} \\
{\em 3)}\quad $\alpha \in dom \left (\mathbf{S}_{\varphi
}^{\vartriangleleft \alpha _{1}} ( \overrightarrow{a} ) \right)$ \
iff there exists some \ $\mathfrak{M}$-generic function
\[
    l\stackrel{\ast }{\in }\Delta _{\varphi }^{
    \vartriangleleft \alpha _{1}} ( \alpha , \overrightarrow{a} ).
\]
\end{lemma}
Here the important statement 1) comes directly from lemma
\ref{1.1.}. This statement will be frequently used in what
follows.
\\

The so called \textit{universal spectrum} is distinguished among
all other spectra.
\\
It is well known that the class \ $\Sigma_{n} (
\overrightarrow{a})$ \ for \  $n>0$ \  contains the formula which
is universal for this class (Tarski~\cite{Tarski}, see Addison
\cite{Addison}); let us denote it by \ $U_{n}^{\Sigma } (
\mathfrak{n}, \overrightarrow{a},\underline{l} ) $. \ Hence, it is
universal for the class  \
$\Sigma_{n}^{\dashv\vdash}(\overrightarrow{a}) $ \ also. Its
universality means that for any \  $\Sigma _{n}^{\dashv\vdash} (
\overrightarrow{a} ) $-formula \  $ \varphi  (
\overrightarrow{a},\underline{l} ) $ \  there is a natural \
$\mathfrak{n}$ \ such that
\[
    \varphi  ( \overrightarrow{a},\underline{l} ) \dashv \vdash
    U_{n}^{\Sigma } ( \mathfrak{n},\overrightarrow{a},\underline{l} ) ;
\]
this \ $\mathfrak{n}$ \ is named the G\"{o}del number of \
$\varphi $. \ The dual formula universal for
$\Pi_{n}^{\dashv\vdash} ( \overrightarrow{a} ) $  \  is denoted by
$ U_{n}^{\Pi } ( \mathfrak{n},\overrightarrow{a},\underline{ l} )
$. \ For some convenience we shall use \  $U_{n}^{\Sigma }$ \ in
the form \  $\exists x~U_{n-1}^{\Pi } (
\mathfrak{n},x,\overrightarrow{a}, \underline{l} ) $. \ In this
notation the upper indices \ $^{\Sigma} $,~$^{\Pi} $ \ will be
omitted in the case when they can be  restored from the context or
arbitrary.
\\
One should recall that universal formulas come out on the
following way. Let us consider a formula  \ $ \varphi\left(
\overrightarrow{{a}}, \underline{l}\right)$ \ and its equivalent
prenex form of the \textit{minimal} level $e$:
\[
    Q_{1}x_{1}Q_{2}x_{2}...,Q_{i}x_{i}\varphi_{1}\left(
    x_{1},x_{2}...,x_{i},\overrightarrow{{a}},\underline{l}\right),
\]
where \ $Q_{1}x_{1}Q_{2}x_{2}...,Q_{i}x_{i}$ \ is its quantor
prefix and \ $ \varphi_{1}$ \ contains only bounded quantors;
suppose \  $Q_{i}=\exists, \quad e>0$. \ Removing these last
bounded quantors to the left we can receive the equivalent
formula of the kind
\begin{eqnarray*}
    Q_{1}x_{1}Q_{2}x_{2}...,Q_{i}x_{i}
    \quad \qquad \qquad \qquad \qquad \qquad \qquad \qquad \qquad
    \qquad \qquad
\\
    \quad \forall y_{1}\in z_{1} \exists v_{1}\in w_{1}\forall y_{2}
    \in z_{2}\exists v_{2} \in w_{2}...,\forall y_{j}\in z_{j}
    \exists v_{j}\in w_{j}\qquad \qquad
\\
    \varphi_{2}\left(x_{1},x_{2}...,x_{i},y_{1},z_{1},
    v_{1},w_{1},y_{2},z_{2},v_{2},w_{2}...,y_{j},z_{j},v_{j},
    w_{j},\overrightarrow{{a}},\underline{l}\right),
\end{eqnarray*}
\vspace{0pt}

\noindent where each of the variables \
$z_{1},w_{1}...,z_{j},w_{j}$ \ is bounded by the others or by the
term from the train \ $x_{1}...,x_{i},\overrightarrow{{a}}$ \ and
\ $\varphi_{2}$ \ contains no quantors. Using Skolem functions
techniques one can transform it to the equivalent formula$ \
U_{\varphi}\left( \overrightarrow{{a}},\underline{l}\right)$ \ of
the same level \ $e$:
\[
    Q_{1}x_{1}Q_{2}x_{2}...,Q_{i}x_{i}\exists f_{1} \exists
    f_{2}...,\exists f_{j}\quad
    \varphi_{3}\left(x_{1},x_{2}...,x_{i},f_{1},f_{2}...,f_{j},
    \overrightarrow{{a}},\underline{l}\right),
\]
where \ $f_{1}...,f_{j}$ \ - variable Skolem functions and \
$\varphi_{3}$ \ contains only bounded quantors of the standard
kind and disposition.
\\
The case \ \mbox{$Q_{i}=\forall$} \ one should consider in a dual
way.

We shall name this formula \
$U_{\varphi}\left( \overrightarrow{{a}},\underline{l}\right)$ the
\textit{preuniversal form} of the formula  \ $\varphi
\left(\overrightarrow{{a}}, \underline{l}\right)$ under
consideration. After that it is easy to receive the equivalent
formula with blocks of like quantors contracted into one through
application of the Pairing Axiom. It remains to apply to formulas
like this one the recursive enumaration of their subformulas
following after the quantor prefix.

Let us assume in what follows that formulas we introduce are
considered in their \textit{preuniversal form} and use for them
the same notation (unless the context does not mean another
situation). This agreement preserves the equivalence of a
formula  \ $\varphi (a,\underline{l}) \in
\Sigma_{e}^{\dashv\vdash} $ \ with a G\"{o}del number \
$\mathfrak{n}$ \ to the universal formula under \
$\vartriangleleft^{ \underline{l}}$~-restriction:
\[
    \varphi^{\triangleleft \alpha_{1}}(\overrightarrow{a},
    \underline{l}) \longleftrightarrow U_{e}^{\Sigma \triangleleft
    \alpha_{1}}( \mathfrak{n},\overrightarrow{a},\underline{l})
\]
for ordinals  \  $\alpha_{1} $ \  of many kinds (for example, for
cardinals). Analogously for  \ $\Pi_{e}^{\dashv\vdash}$ \ and for
generic equivalence.

\begin{definition}
\label{2.6.}
\quad \\
{\em 1)}\quad We name as the
spectral universal for the class \ $\Sigma_{n}^{\dashv \vdash }$ \
formula of level \  $n$  \ the formula \  $u_{n}^{\Sigma }
(\overrightarrow{a}, \underline{l} ) $ \  obtained from the
universal formula \ $U_{n}^{\Sigma } (
\mathfrak{n},\overrightarrow{a},\underline{l} ) $ \  by replacing
all occurrences of the variable \ $\mathfrak{n}$ \ by occurrences
of the term \  $\underline{l} ( \omega _{0} ) $.
\\
The spectral universal for the class $\Pi_{n}^{\dashv \vdash }$
formula \ $u_{n}^{\Pi } ( \overrightarrow{a},\underline{l} )$ is
introduced in the dual way. Thus, we take
\[
     u_{n}^{\Sigma } ( \overrightarrow{a}, \underline{l} ) =
     \exists x~u_{n-1}^{\Pi } ( x,\overrightarrow{a},\underline{l}
     ),
\]
where \ $u_{n-1}^{\Pi} ( x, \overrightarrow{a}, \underline{l} )$
\ is the spectral universal for the class \ $\Pi_{n-1}^{\dashv
\vdash }$ \ formula.
\quad \\
{\em 2)}\quad The Boolean values
\[
    A_{\varphi}^{\vartriangleleft \alpha _{1}} ( \alpha ,
    \overrightarrow{a} ), \Delta _{\varphi }^{\vartriangleleft
    \alpha _{1}} ( \alpha,\overrightarrow{a} )
    \mathrm{\it \ \ and\ the\ spectrum\ \ } \mathbf{S}_{\varphi}^{\vartriangleleft
    \alpha _{1}} ( \overrightarrow{a} )
\]
of the formula \ $\varphi =u_{n}^{\Sigma } (\overrightarrow{a},
\underline{l} )$ \ and its projections (definition 2.4 where \
$\trianglelefteq^{\underline{l}}$ \ should be replaced with \
$\leq$)
 are named the
universal Boolean values and spectra of the level \ $n$ \  on the
point \ $\overrightarrow{a}$ \  below \ $\alpha _{1}$ \ and in
their notation the index \ $u_{n}^{\Sigma }$  \  is replaced by \
$n$, \  that is they are denoted by
\[
    A_{n}^{\vartriangleleft \alpha _{1}} ( \alpha ,
    \overrightarrow{a} ), \quad \Delta _{n}^{\vartriangleleft
    \alpha_{1}} ( \alpha ,\overrightarrow{a} ) ,\quad
    \mathbf{S}_{n}^{\vartriangleleft \alpha _{1}}
    ( \overrightarrow{a} ).
\]
\hfill {} \\
{\em 3)}\quad If \ $( \alpha ,\Delta  ) \in
\mathbf{S}_{n}^{\vartriangleleft \alpha _{1}} ( \overrightarrow{a}
) $, \ then \ $\alpha $ \ is named the jump ordinal of these
formula and spectra, while \ $\Delta $ \ is named its Boolean value on the
point \ $\overrightarrow{a}$ \ below \ $\alpha _{1}$.
\\
{\em 4)}\quad The ordinal \  $\alpha _{1}$ \ is named
the carrier of these spectra.
\end{definition}

Everywhere further  \  $\vartriangleleft ^{\underline{l}}$~- \
or  \  $\vartriangleleft $-bounding ordinals  \ $\alpha_1 $  \
are limit cardinals  \ $<k $ \ (in \ $L_{k} $) or \ $\alpha_1 =k $
\ (if the context does not point to the opposite).
\\
\hfill {} \\
Here the term ``universal'' is justified by the following

\begin{lemma} \label{2.7.} \hfill {} \\
\hspace*{1em} For every proposition \  $\varphi =\exists x\
\varphi _{1} ( x, \overrightarrow{a},\underline{l} ) ,\quad
\varphi _{1}\in \Pi _{n-1}^{\dashv \vdash }$:
\[
    dom \left (
    \mathbf{S}_{\varphi }^{\vartriangleleft \alpha _{1}} (
    \overrightarrow{a} )\right ) \subseteq dom \left (
    \mathbf{S}_{n}^{\vartriangleleft \alpha _{1}}
    (\overrightarrow{a} ) \right) .
\]
\end{lemma}
\vspace{-12pt}

\noindent \textit{Proof.} \ For a function  \  $l\in {}^{k}k $  \
let \ $l_0 $  \  denote every function  $\in {}^{k}k $  \ taking
values:
\\
\begin{equation} \label{e2.1}
l_0 \left( \alpha\right) = \left\{
\begin{array}{ll}
l(\alpha) & ,  \ \alpha< \omega_{0} \vee \alpha > \omega_{0}+1; \\
(l(\omega_{0}), l(\omega_{0}+1)) & , \ \alpha = \omega_{0}+1;
\end{array} \right.
\end{equation}
\hfill {} \\
the value  \ $l_0(\omega_{0}) $ \ is arbitrary here. Evidently \
$l $  \ is \ $\mathfrak{M}$-generic \textit{iff} \ $l_0$  \
possesses the same property. Let  \ $\varphi_{0} \left
(\overrightarrow{a}, l \right )$  \  denote the formula received
by the successive replacing every subformula \ $\underline{l}
\left (t_{1} \right ) = t_{2} $ \ of the formula \ $\varphi $  \
by the bounded subformula:
\\
\begin{eqnarray} \label{e2.2}
& \nonumber \exists y_{1},y_{2} < \omega_{0} ( \underline{l}
\left( \omega_{0}+1 \right) = (y_{1},y_{2})\wedge
\\
& \wedge \left ( \left ( t_{1} < \omega_{0} \vee t_{1}
> \omega_{0} + 1 \right) \longrightarrow \underline{l} (t_{1}) =
t_{2}  \right) \wedge
\\
& \nonumber \wedge \left ( t_{1} = \omega_{0} \longrightarrow
y_{1} = t_{2} \right )\wedge \left ( t_{1} = \omega_{0} +1
\longrightarrow y_{2} = t_{2} \right ) ).
\end{eqnarray}
\hfill {} \\
The subformulas of the form  \  $ \underline{l} (t_{1}) \in
t_{2}, \quad t_{2} \in \underline{l} (t_{1}) $  \  are treated
analogously. Let  \  $ \mathfrak{n}_{0} $  \  be the G\"{o}del
number of \ $\varphi_{0}$. \ For an arbitrary  \
$\mathfrak{M}$-generic function \ $l $  \  let us assume that \
$l_0(\omega_{0}) = \mathfrak{n}_{0} $. Obviously, for every  \
$\alpha < \alpha_{1} $

\[
    l \stackrel{\ast }{\in }~ A_{\varphi}^{\triangleleft \alpha_{1}}
    (\alpha, \overrightarrow{a}) \longleftrightarrow l_0
    \stackrel{\ast}{\in }~ A_{\varphi_{0}}^{\triangleleft \alpha_{1}}
    (\alpha, \overrightarrow{a}) \longleftrightarrow l_0
    \stackrel{\ast }{\in }~ A_{n}^{\triangleleft \alpha_{1}}
    (\alpha, \overrightarrow{a})\quad.
\]
\vspace{0pt}

It remains to apply now lemma 2.5 \ 3). \hfill $\dashv$ \\

From here one can imply that universal spectra accumulate while
their ordinal constants are increasing. For \ $\overrightarrow{a}=
( \alpha_{1}...,\alpha _{m} )$ \ let \ $\max \overrightarrow{a} =\max \{ \alpha _{1}...,\alpha _{m}\} $.
\\

\begin{lemma}\label{2.8.} \hfill {} \\
\hspace*{1em} For every \ $\overrightarrow{a}_{1},
\overrightarrow{a}_{2}<\alpha _{1}$:
\hfill {} \\
\[
    \max \overrightarrow{a}_{1}<\max \overrightarrow{a}_{2}
    \longrightarrow dom \left ( \mathbf{S}_{n}^{\vartriangleleft
    \alpha _{1}} (\overrightarrow{a}_{1} )\right )
    \subseteq dom \left ( \mathbf{S}_{n}^{\vartriangleleft
    \alpha _{1}} (\overrightarrow{a}_{2} )\right) .
\]
\end{lemma}
The proof can be carried out by the so called splitting method
(see the proof of lemma 4.6 below for example), but this lemma,
although clarifying spectrum properties, is not used further and
so its proof is omitted here.

\newpage
\quad 

\newpage

\section{Subinaccessible Cardinals}
\setcounter{equation}{0}

Here is developed the theory of subinaccessibility in its basic
aspect.

The further reasoning is conducted in \  $L_{k}$ \  (or in  \
$\mathfrak{M} $ \  if the context does not mean some another
situation).
\\
Let us introduce the central notion of subinaccessibility -- the
inaccessibility by means of our language. The ``meaning'' of
propositions is contained in their spectra and therefore it is
natural to define this inaccessibility by means of the spectra
of all propositions of a given level:

\begin{definition} \label{3.1.} \hfill {} \\
\hspace*{1em} Let \  $\alpha _{1} \leq k $.  \\
We name an ordinal  \ $\alpha <\alpha _{1}$ \  subinaccessible of
\  $a $ \  level \ $n$ \ below \  $\alpha _{1}$ \  iff it fulfills
the following formula denoted by \  $SIN_{n}^{<\alpha _{1}} ( \alpha )$:
\hfill {} \\
\[
    \forall \overrightarrow{a} <\alpha \quad dom
    \left ( \mathbf{S}_{n}^{\vartriangleleft \alpha_{1}}
    ( \overrightarrow{a} ) \right ) \subseteq \alpha .
\]
The set \
\[
    \{ \alpha <\alpha _{1}:SIN_{n}^{<\alpha _{1}} ( \alpha ) \}
\]
of all these ordinals is denoted by \ $SIN_{n}^{<\alpha _{1}}$ \
and they are named \ $SIN_{n}^{<\alpha _{1}}$-ordinals.
\\
As usual, for \ $\alpha _{1}<k$ \  we say that subinaccessibility
of \  $\alpha $ \  is restricted by \  $\alpha _{1}$ \ or
relativized to \  $\alpha _{1}$; \ for \ $\alpha_1 = k$ \ the
upper indices \ $< \alpha_1$, $\vartriangleleft \alpha_1$ \ are
dropped. \hfill $\dashv$
\end{definition}

Obviously, the cardinal \ $k$ \ is \textit{subinaccessible} itself
of any level, if we define this notion for \ $\alpha = \alpha_1 =
k$.
\\
So, the comparison of the notions of inaccessibility and
subinaccessibility naturally arises in a following way:
\\
The cardinal \ $k$ \ is weakly inaccessible, since it is
uncountable and cannot be reached by means of smaller powers in
sense that it possess two properties: 1) it is regular and 2) it
is closed under operation of passage to next power: \ $\forall
\alpha < k ~~~ \alpha^+ < k$.
\\
Turning to \textit{subinaccessibility} of an ordinal \ $\alpha <
k$ \ of the level \ $n$ \ (in its unrelativized form for some
brevity), one can see that the property of regularity is dropped
now, but \ $\alpha$ \ still can not be reached, but by another
more powerful means: the second condition is strengthened  and now
\ $\alpha$ \ is closed under more powerful operations of passage
to jump ordinals of universal spectrum:
\[
    \forall \overrightarrow{a} <\alpha \quad \forall \alpha^{\prime}
    \in dom \mathbf{S}_{n} ( \overrightarrow{a} )
    \quad \alpha^{\prime} < \alpha,
\]
that is by means of ordinal spectra \textit{of all}
propositions of level \ $n$ \ (see lemma 2.7 above).
\\
It implies the closure of \ $\alpha$ \ \textit{under all \
$\Pi_{n-1}^{\dashv\vdash}$-functions in all generic extensions}
of \ $L_k$, \ not only under operation of power successor in \
$L_k$ \ (see lemma 3.5 below).
\\

It is evident that working in \ $L_{k} $ \ one should treat the
formula \ $SIN_{n}^{< \alpha_{1}}(\alpha) $ \ actually as two
formulas: one of them without the constant \ $\alpha_{1} $ \ when
\ $\alpha_{1}=k$, \ and another containing \ $\alpha_{1} $ \
when  \  $\alpha_{1} < k$; \ the same remark concerns all
formulas, constructions and reasoning containing some parameter \
$\alpha_{1} \leq k$.\\
\quad \\

From definition 3.1 and lemma 2.7 obviously comes

\begin{lemma}
\label{3.2.} {\em (About restriction)} \hfill {} \\
\hspace*{1em} Let \  $\alpha <\alpha _{1}\leq k$, \ $\alpha \in
SIN_{n}^{<\alpha _{1}}$  \ and a proposition \  $ \exists
x~\varphi  ( x,\overrightarrow{a}, \underline{l} ) $ \ has \  $
\overrightarrow{a}<\alpha $, \ $\varphi \in \Pi _{n-1}^{\dashv
\vdash }$, \ then for every \ $\mathfrak{M}$-generic \  $l$
\\

\[
    L_{k}\left[ l\right] \vDash \left( \exists x\vartriangleleft ^{l}
    \alpha _{1}\mathit{\ }\varphi ^{\vartriangleleft \alpha _{1}}
    ( x,\overrightarrow{a}, l ) \longrightarrow \exists x
    \vartriangleleft ^{l}\alpha \ \varphi ^{\vartriangleleft
    \alpha _{1}} ( x,\overrightarrow{a},l) \right) \ ,
\]
\hfill {} \\
In this case we say that below \  $\alpha _{1}$ \  the ordinal \
$\alpha $ \  restricts or relativizes the proposition \ $\exists
x~\varphi $.
\\
Considering the same in the inverted form for \ $\varphi \in
\Sigma _{n-1}^{\dashv \vdash }$:
\\

\[
    L_{k}\left[ l\right] \vDash \left( \forall x\vartriangleleft ^{l}
    \alpha \ \varphi ^{\vartriangleleft \alpha _{1}} ( x,
    \overrightarrow{a},l ) \longrightarrow \forall x
    \vartriangleleft ^{l}\alpha _{1}\varphi ^{\vartriangleleft
    \alpha _{1}} ( x,\overrightarrow{a},l ) \right) \ ,
\]
\hfill {} \\
we say that below \  $\alpha _{1}$ \  the ordinal \  $\alpha $ \
extends or prolongs the proposition \  $\forall x\ \varphi $ \ up
to \ $\alpha _{1}$.
\end{lemma}

Of course, lemma~3.2 presents the stronger statement, the
criterion of \ $SIN_{n}^{<\alpha _{1}}$-subinaccessibility.
\label{c3}
\endnote{
\ p. \pageref{c3}. \ Obviously, the subinaccessibility notion can be
reformulated in terms of elementary equivalence but in the following
substantially more complicated and artificial form:

$\alpha$ \ is subinaccessible of level \ $n$ \ below \ $\alpha_1$
\ \textit{iff} $L_{\alpha}[l | \alpha_1]$ \ is \
$\Sigma_n$-elementary substructure of \ $L_{\alpha_1}[l |
\alpha_1]$ \ \textit{for every \ \mbox{$\mathfrak{M}$-generic}
function} \ $l$.
\\
Precisely this is done in lemma 3.2. Here few words must be said
in connection with this aspect.
\\
This notion can not be reduced naturally to the elementary
equivalence of constructive segments \ $L_{\alpha}$, \
$L_{\alpha_1}$ \ only, but requires involving all their generic
extensions of this kind.
\\
So, using these terms, when parameters \ $\alpha$, $\alpha_1$ \
are varying simultaneously along with many other cardinals that
are subinaccessible also, one receives some multilayer and
cumbersome description of the notion which is especially
unnatural, since in fact this notion works within \ $L$ \ only.
\\
Moreover, \ the notion of disseminator (\S 6), which is simple
generalization of the subinaccessibility notion, becomes
unnaturally complicated in terms of elementary substructures.

Therefore some more convenient description is required, pointing
directly to the very essence of the phenomenon, simply definable
and therefore more suitable for the investigation of the problem
-- that is subinaccessibility introduced above.
\\
\quad \\
} 
\\

Now the following lemmas \ref{3.3.} -- \ref{3.8.} can be easily
deduced from definition 3.1 and lemma \ref{2.5.} :

\begin{lemma} \label{3.3.} \hfill {} \\
\hspace*{1em} The formula \  $SIN_{n}^{<\alpha _{1}} ( \alpha
 ) $ \  belongs to the class \ $\Pi _{n} $ \
for \ $\alpha _{1}=k$ \ and to the class \ $\Delta _{1} $ \ for \
$\alpha _{1}<k$.
\end{lemma}

\vspace{-6pt}

\begin{lemma} \label{3.4.} \hfill {} \\
\hspace*{1em} For every \ $n>0$:\\
\quad \\
{\em 1)}\quad the set \  $SIN_{n}^{<\alpha _{1}}$ \  is closed in
\  $\alpha _{1}$, \ that is for any \  $\alpha <\alpha _{1}$
\[
    \sup  ( \alpha \cap SIN_{n}^{<\alpha _{1}} ) \in
    SIN_{n}^{<\alpha_{1}};
\]
{\em 2)}\quad the set \  $SIN_{n}$ \  is unbounded in \ $k$,
that is \ $\sup SIN_{n}=k$;\\
\quad \\
{\em 3)}\quad \quad \quad \quad \quad \quad \quad \quad \quad
$SIN_{n}^{<\alpha _{1}} ( \alpha  ) \longleftrightarrow
SIN_{n}^{\vartriangleleft \alpha _{1}} ( \alpha )$.
\end{lemma}

\vspace{-6pt}

\begin{lemma} \label{3.5.} \hfill {} \\
\hspace*{1em} Let \  $\alpha \in SIN_{n}^{<\alpha _{1}}$ \ and a
function \ $f\subset \alpha _{1}\times \alpha _{1}$ \ be defined
in  \ $L_{k}[l] $  \  by a formula \  $\varphi ^{\vartriangleleft
\alpha _{1}} ( \beta ,\gamma ,l ) $ \  where \  $\varphi \in \Pi
_{n-1}^{\dashv \vdash}$,
\\
then \  $\alpha $ \  is closed under \  $f$.
\\
In particular, for every \  $n \geq 2$
\[
    \mathrm{\quad if \quad} \alpha \in SIN_{n} \mathrm{\quad then \quad}
    \alpha =\omega _{\alpha } \mathrm{\quad (in\ } L_{k} ).
\]
\end{lemma}
The following lemmas 3.6-3.8 represents the important technical
tools of subinaccessibility investigations.

\begin{lemma} \label{3.6.} \hfill {} \\
\hspace*{1em} For every \ $m<n$
\\
1) \hspace*{3.5cm} $SIN_{n}^{<\alpha _{1}} \subset
SIN_{m}^{<\alpha _{1}}$;
\\
\\
2) moreover, every \ $\alpha \in SIN_{n}^{<\alpha _{1}}$ \ is a
limit ordinal in \ $SIN_{m}^{<\alpha _{1}}$~:
\[
    sup \left ( \alpha  \cap SIN_{m}^{<\alpha_{1}} \right ) = \alpha
\]
\end{lemma}
\textit{Proof.} \ Statement 1) is obvious, because every \
$\Sigma_m$-formula is at the same time \ $\Sigma_n$-formula.
\\
Turning to 2) let us consider the \ $\Sigma_n$-formula
\[
    \exists\gamma(\beta<\gamma\wedge SIN_m(\gamma))
\]
with arbitrary constant \ $\beta<\alpha$. \ This formula is true
below \ $\alpha_1$, because due to 1) the ordinal \ $\alpha$ \
itself can be used as \ $\gamma$. \ After that \
$SIN_{n}^{<\alpha_{1}}$-ordinal \ $\alpha$ \ restricts this
formula and some \ $SIN_{m}^{<\alpha_{1}}$-ordinal \
$\gamma>\beta$ appears below \ $\alpha$. \hfill $\dashv$
\\

\noindent It is obvious that the converse statement is false:
subinaccessible cardinals sometimes lose this property on the next
level -- for example, all the successors in a given class \
$SIN_{m}^{<\alpha _{1}}$.

\begin{lemma} \label{3.7.} \hfill {} \\
\hspace*{1em} Let
\[
    \overrightarrow{a}<\alpha _{2}<\alpha_{1}\leq k ,\quad
    \alpha _{2}\in SIN_{n}^{<\alpha _{1}}
\]
then for any \ $ Q_{n}^{\dashv\vdash}$-proposition \ $\varphi  (
\overrightarrow{a},\underline{l} )$
\\

\[
     \varphi ^{\vartriangleleft \alpha _{1}} ( \overrightarrow{a},
     \underline{l} ) \dashv \vdash \varphi ^{\vartriangleleft
     \alpha_{2}} ( \overrightarrow{a},\underline{l} ).
\]
\end{lemma}
\textit{Proof.} \ Let us consider the proposition  \  $\varphi =
\exists x \ \ \varphi_{1} \left (x,\overrightarrow{a}
,\underline{l} \right ) $,  \  $\varphi_{1}\in \Pi_{n-1} $, \ and
the ordinal
\[
    \alpha_{0} = min \left \{ \alpha : L_{k} [l] \vDash
    \varphi_{1}^{\triangleleft \alpha_{1}} \left (F^{l}(\alpha),
    \overrightarrow{a}, l \right ) \right \}
\]
for an  \  $\mathfrak{M} $-generic  \  $l $. \ By lemmas 2.5 \ 3),
2.7 \ $\alpha_{0}\in dom \left (\mathbf{S}_{n}^{\triangleleft
\alpha_{1}} (\overrightarrow{a}) \right ) $. \ Since  \
$\alpha_{2}\in SIN_{n}^{< \alpha_{1} }$ \  it implies  \
$\alpha_{0} < \alpha_{2} $ \  and therefore the  \
$\vartriangleleft^{l}$-restriction \  $\alpha_{1}$ \  in the
proposition  \  $ \varphi_{1}^{\triangleleft \alpha_{1}} \left
(F^{l}(\alpha_{0}), \overrightarrow{a}, l \right ) $ \  can be
replaced by the  \ $\vartriangleleft^{l}$-restriction \ $\alpha_{2}$.
Hence
\[
    L_{k}[l]\vDash \left( \varphi^{\triangleleft \alpha_{1}}
    (\overrightarrow{a},l) \longrightarrow
    \varphi^{\triangleleft \alpha_{2}}(\overrightarrow{a},l) \right) \ .
\]
It remains to convert this argument. \hfill $\dashv$

\begin{lemma} \label{3.8.} \hfill {} \\
\hspace*{1em} Let \ $\alpha _{2}<\alpha _{1}\leq k$, \ then:\\
\quad \\
{\em 1)}\quad If \  $\alpha _{2}\in SIN_{n-1}^{<\alpha _{1}}$,
then the set \  $SIN_{n}^{<\alpha _{1}}\cap \alpha _{2}$ \
constitutes the initial segment of the set \
$SIN_{n}^{<\alpha_{2}}$, \  that is:
\\

 \quad \

\qquad (i)\quad $SIN_{n}^{<\alpha _{1}}\cap \alpha
_{2}\subseteq SIN_{n}^{<\alpha _{2}}$;\\

\qquad (ii)\quad $SIN_{n}^{<\alpha_{2}} \cap sup \left (
SIN_{n}^{<\alpha_{1}} \cap \alpha_{2} \right )
\subseteq SIN_{n}^{<\alpha_{1}} $.\\
\quad \\
\hfill {} \\
{\em 2)}\quad If \  $\alpha _{2}\in SIN_{n}^{<\alpha _{1}}$, then
for every \ $m \leq n$ \
\[
    SIN_{m}^{<\alpha _{1}}\cap \alpha_{2}=SIN_{m}^{<\alpha _{2}} .
\]
\end{lemma}
\textit{Proof.} \ Statements of this lemma one can prove in a
similar
way and so we shall demonstrate it for $1.(i).$
\\
Let \ $\alpha \in SIN_{n}^{<\alpha _{1}}\cap \alpha_{2}$; \ it
should be proved, that for every \ $\overrightarrow{a} < \alpha $
\[
    dom \left ( \mathbf{S}_{n}^{\triangleleft \alpha_{2}}
    (\overrightarrow{a}) \right )  \subseteq \alpha,
\]
so let \ $\beta \in dom \left ( \mathbf{S}_{n}^{\triangleleft
\alpha_{2}} (\overrightarrow{a}) \right ) $. \ Now we have  \
$\beta < \alpha_{2} , \overrightarrow{a} < \alpha_{2} $ \  and
 \  $\alpha_{2} \in SIN_{n-1}^{<\alpha_{1}} $, that is why
for every  \  $\beta^{'} \leq \beta $ \  we can replace the  \
$\vartriangleleft^{\underline{l}} $~-boundary \ $\alpha_{2} $
\ in the proposition  \  $u_{n-1}^{\Pi \triangleleft \alpha_{2}}
\left (F^{\underline{l}} (\beta^{'}), \overrightarrow{a},
\underline{l} \right ) $ \  with the  \
$\vartriangleleft^{\underline{l}} $~-boundary \  $\alpha_{1}
$. \ From here and definition 3.1 it comes  \  $\beta \in dom
\left ( \mathbf{S}_{n}^{\triangleleft \alpha_{1}}
(\overrightarrow{a}) \right ) $ \  and then  \  $ \alpha \in
SIN_{n}^{<\alpha_{1}} $ \ implies  \  $\beta < \alpha $. \ It
means that  \  $dom \left ( \mathbf{S}_{n}^{\triangleleft
\alpha_{2}} (\overrightarrow{a}) \right ) \subseteq \alpha $ \ and
hence  \  $\alpha \in
SIN_{n}^{<\alpha_{2}} $.  \hfill $\dashv$\\
\quad \\

When formulas are equivalently transformed their spectra can
change. It is possible to use this phenomenon for the analysis of
subinaccessible cardinals. To this end we shall introduce the
universal formulas with ordinal spectra containing only
subinaccessible cardinals of smaller level. For more clearness of
constructions formulas without individual constants will be
considered. Let us start with the spectral universal formula for
the class \ $\Sigma _{n}^{\dashv \vdash }$. \  The upper indices \
$^{\Sigma} $, $ ^{\Pi}$ \ will be omitted as usual (if it will
not cause misunderstanding). \\
In what follows it is sufficient to consider bounding ordinals
only from the class \ $SIN_{n-2}$, \ therefore everywhere further
\ $\vartriangleleft^{\underline{l}}$\;- \ or \ $\vartriangleleft
$-bounding ordinals \ $\alpha $ \ are assumed to be \
$SIN_{n-2}$-ordinals or \ $\alpha =k$ \ (if the context does
not mean another case). \\
Thus all such ordinals \ $\alpha \leq k$ \ are cardinals
$\alpha=\omega_{\alpha}$ due to lemma 3.5.
\hfill {} \\

\begin{definition}
\label{3.9.} \hfill {} \\
{\em 1)}\quad We name as the
monotone spectral universal for the class \  $\Sigma _{n}^{\dashv
\vdash }$ \ formula of the level \ $n$ \  the \  $\Sigma
_{n}$-formula
\[
    \widetilde{u}_{n} ( \underline{l} ) =\exists x~\widetilde{u}_{n-1}
    ( x, \underline{l} )
\]
where \  $\widetilde{u}_{n-1} ( \underline{l} ) \in \Pi _{n-1}$
\  and
\[
    \widetilde{u}_{n-1} ( x,\underline{l} ) \dashv \vdash
    \exists x^{\prime}\vartriangleleft
    ^{\underline{l}}x~u_{n-1}^{\Pi}
    ( x^{\prime },\underline{l} ) .
\]
\hfill {} \\
{\em 2)}\quad We name as the
subinaccessibly universal for the class \  $\Sigma _{n}^{\dashv
\vdash }$ \  formula of the level \  $n$ \  the \  $\Sigma
_{n}$-formula
\[
    \widetilde{u}_{n}^{\sin } ( \underline{l} ) =
    \exists x~\widetilde{u}_{n-1}^{\sin} ( x,\underline{l} )
\]
where \ $\widetilde{u}_{n-1}^{\sin }\in \Pi _{n-1} $  \  and
\[
    \widetilde{u}_{n-1}^{\sin } ( x,\underline{l} )
    \dashv \vdash SIN_{n-1} ( x )
    \wedge \widetilde{u}_{n-1} ( x ,
    \underline{l} ) .
\]
\hfill {} \\
The monotone and subinaccessibly universal for the class \
$\Pi_{n}^{\dashv \vdash} $ \ formulas are introduced in the dual
way.
\\
\hfill {} \\
{\em 3)}\quad The Boolean values
\[
    A_{\varphi}^{\vartriangleleft \alpha _{1}} ( \alpha  ), \
    \Delta_{\varphi }^{\vartriangleleft \alpha _{1}} ( \alpha  ),
    \mathrm{\it \ \ and\ the\ spectrum\ \ }
    \mathbf{S}_{\varphi }^{\vartriangleleft \alpha_{1}}
\]
of the formula \ $\varphi = \widetilde{u}_{n}^{\sin }$ \ and its
projections (definition 2.4 where \
$\trianglelefteq^{\underline{l}}$ \ should be replaced with \
$\leq$) are named subinaccessibly universal of the level \ $n$ \
below \ $\alpha_{1}$ \ and are denoted respectively by

\[
    \widetilde{A}_{n}^{\sin \vartriangleleft \alpha _{1}} ( \alpha ) ,
    \quad \widetilde{\Delta }_{n}^{\sin \vartriangleleft \alpha _{1}}
    ( \alpha  ) ,\quad \widetilde{\mathbf{S}}_{n}^{\sin \vartriangleleft
    \alpha _{1}}.
\]
\hfill {} \\
{\em 4)}\quad If \ $( \alpha ,\Delta  ) \in
\widetilde{\mathbf{S}}_{n}^{\sin \vartriangleleft \alpha _{1}}$, \
then \ $\alpha $ \ is named
the jump ordinal of these formula and spectra, while \ $\Delta $ \
is named its Boolean value below \ $\alpha _{1}$.
\hfill {} \\
{\em 5)}\quad The cardinal \ $\alpha_{1}$ \ is named
the carrier of these spectra.
\end{definition}
From this definition and lemmas 3.4~~2) ( for \  $n-1$ \ instead
of \  $n$), 2.5~~3), 2.7 and 3.8 the following easy lemmas come:

\begin{lemma}
\label{3.10.} \quad \\
{\em 1)}\qquad \qquad \qquad \qquad \qquad $u_{n}^{\Sigma} (
\underline{l} ) \dashv \vdash \widetilde{u}_{n}^{\sin } (
\underline{l} )$~.
\\
\\
{\em2)} \qquad   $ dom \ \bigl ( \widetilde{\mathbf{S}}_{n}^{\sin
\vartriangleleft \alpha _{1}} \bigr ) \subseteq
SIN_{n-1}^{<\alpha_{1}} \cap dom  \bigl
(\mathbf{S}_{n}^{\triangleleft \alpha_{1}} \bigr )$~.

\quad \\
{\em 3)} Let
\[
    \alpha \in dom \bigl ( \mathbf{S}_{n}^{\triangleleft \alpha_{1}} \bigr )
    \mathrm{\quad and \quad}
    \alpha^{\prime} = min  \left \{ \alpha^{\prime \prime} \in
    SIN_{n-1}^{< \alpha_{1}} : \alpha^{\prime \prime} > \alpha
    \right \}~,
\]
then
\[
    \alpha^{\prime} \in dom \bigl (
    \widetilde{\mathbf{S}}_{n}^{\sin \triangleleft \alpha_{1}}
    \bigr )~.
\]
{\em 4)} Let
\[
    \alpha \leq \alpha_{1} \mathrm{\it \quad be\ limit\ in \quad}
    SIN_{n-1}^{< \alpha_{1}}~,
\]
then
\[
    \mathrm{sup} dom \bigl ( \widetilde{\mathbf{S}}_{n}^{
    \sin \triangleleft \alpha_{1}} | \alpha \bigr ) =
    \mathrm{sup} dom \bigl ( \mathbf{S}_n^{\triangleleft
    \alpha_{1}} | \alpha \bigr ).
\]
\end{lemma}

\vspace{-6pt}

\begin{lemma} \label{3.11.} \hfill {} \\
\hspace*{1em} Let
\[
    \alpha_{2} \in SIN_{n-2}^{< \alpha_{1}} \quad \mbox{\it and}\quad
    \alpha_{0} =\sup \bigl ( SIN_{n-1}^{<\alpha _{1}}\cap \alpha_{2} \bigr ),
\]
then
\[
    \widetilde{\mathbf{S}}_{n}^{\sin \vartriangleleft
    \alpha _{2}}  | \alpha_{0} =\widetilde{\mathbf{S}}_{n}^{\sin
    \vartriangleleft \alpha _{1}} | \alpha_{0} \quad .
\]
\end{lemma}
\textit{Proof.} \ Let  \  $\alpha<\alpha_{0},(\alpha, \Delta )
\in \widetilde{\mathbf{S}}_{n}^{\sin \triangleleft \alpha_{2} }
$  \  and hence
\[
    \Delta = \widetilde{\Delta}_{n}^{\sin
    \triangleleft \alpha_{2} } (\alpha)>0.
\]
For every  \  $\alpha^{\prime} \leq \alpha $ \  by lemma 3.10~2)
and lemma 3.8 (used for  \  $n-1 $ \  instead of  \  $n $)
\[
    SIN_{n-1}^{< \alpha_{2}} (\alpha^{\prime})
    \longleftrightarrow SIN_{n-1}^{< \alpha_{1}}
    (\alpha^{\prime})
\]
and so the  \  $\vartriangleleft^{\underline{l}}$~-boundary  \
$\alpha_{2} $ \  in the proposition  \ $\widetilde{u}_{n-1}^{
\triangleleft \alpha_{2}} \left (F^{\underline{l}}
(\alpha^{\prime}), \underline{l}) \right )$ \  can be replaced
with the \ $\vartriangleleft^{\underline{l}}$~-boundary \
$\alpha_{1} $. \ It implies
\[
    \Delta= \widetilde{\Delta}_{n}^{\sin
    \triangleleft \alpha_{1} } (\alpha)
    \quad \mbox{\it and} \quad (\alpha, \Delta ) \in
    \widetilde{\mathbf{S}}_{n}^{\sin \triangleleft \alpha_{1} }.
\]
The converse argument completes the proof. \hfill $\dashv$\\

\quad \\

Our aim is to ``compare'' universal spectra with each other on
\textit{different} carriers \ $\alpha_1$ \ disposed cofinally to
\ $k$ \ in order to introduce \textit{monotone} matrix functions.
To this end it is natural to do it by means of using values of
function \ $Od$ \ for such spectra.
\\
Also it is natural to try to find some estimates of
``informational complexity'' of these spectra by means of
estimates of their order types. But from lemmas \ref{2.7.},
\ref{2.8.}, \ref{3.10.} \ it follows that, for instance, spectra
\ $dom \left ( \widetilde{\mathbf{S}}_{n-1}^{\sin \triangleleft
\alpha_{1}} ( \alpha  ) \right ) $ \ accumulate when \  $\alpha $,
\ $\alpha_{1}$ \ are increasing and therefore their order types
increase up to  \  $k $; \ in addition they are closed under \
$\Pi _{n-2}$-functions, etc. Spectra \
$\widetilde{\mathbf{S}}_{n}^{\sin \triangleleft \alpha_{1}}$ \
have analogous properties.
\\
Therefore the required comparison of such spectra can be hardly
carried out in a proper natural way since they \textit{are ``too much
differ''} from each other for arbitrary great carriers \ $\alpha_1$.
\\
So, there is nothing for it but to consider further
\textit{spectra reduced to some fixed cardinal} and, next,
\textit{reduced matrices}.

\newpage
\quad 

\newpage
\section{Reduced Spectra}
\setcounter{equation}{0}

Here we start to form the main material for building matrix
functions -- reduced matrices.
\\
With this end in view first we shall consider the necessary
preliminary constructions -- reduced spectra.
\\

For an ordinal \ $\chi \leq k$ \ let \ $P_{\chi }$ \ denote the
set \  $\{ p\in P:dom ( p ) \subseteq \chi \} $ \ and \  $B_{\chi
}$ \  denote the subalgebra of \  $B$ \  generated by \
$P_{\chi}$ \  in \  $L_{k}$. \  For every \  $A\in B$ \  let us
introduce the set
\[
    A \lceil \chi =\{ p\in P_{\chi }:\exists q~ ( p=\left. q\right|
    \chi \wedge q\leq A ) \}
\]
which is named the value of
\ $A$ \  reduced to \  $\chi $. \ It is known (see~\cite{Jech})
that
\[
    B_{\chi }=\left \{ \sum X:X\subseteq P_{\chi } \right \}
\]
and therefore every \ $A\in B_{\chi }$  \ coincides with \ $ \sum
A\lceil \chi $. \ Therefore let us identify every \ $A\in B_{\chi
}$ \ with its reduced value \ $A\lceil \chi $; \ so, here one
should point out again, that cause of that every value \ $A \in
B_\chi$ \ is the set in \ $L_k$, not class, and \ $B_\chi$ \ is
considered as the set of such values.

\begin{definition} \label{4.1.} \hfill {} \\
\hspace*{1em} Let \  $\chi \leq k$,\quad $\alpha_{1}\leq k$.
\hfill {} \\
{\em 1)} \quad For every \  $\alpha <\alpha _{1}$  \  let us
introduce the Boolean values and the spectrum:
\hfill {} \\
\[
    \mathbf{\ }\widetilde{A}_{n}^{\sin \vartriangleleft \alpha _{1}}
    (\alpha  ) \lceil \chi \mathbf{;~}\widetilde{\Delta }_{n}^{\sin
    \vartriangleleft \alpha _{1}} ( \alpha  ) \overline{\lceil}\chi =
    \widetilde{A}_{n}^{\sin \vartriangleleft \alpha _{1}} ( \alpha  )
    \lceil\chi \mathbf{-}\sum_{\alpha ^{\prime }<\alpha }
    \widetilde{A}_{n}^{\sin \vartriangleleft \alpha _{1}}
    ( \alpha ^{\prime } ) \lceil \chi ;
\]
\[
    \widetilde{\mathbf{S}}_{n}^{\sin \vartriangleleft \alpha _{1}}
    \overline{\overline{\lceil}}\chi =\{
    ( \alpha ;~\widetilde{\Delta }_{n}^{\sin \vartriangleleft \alpha _{1}}
    ( \alpha  ) \overline{\lceil}\chi  ) :\alpha <\alpha _{1}\wedge
    \widetilde{\Delta }_{n}^{\sin \vartriangleleft \alpha _{1}} ( \alpha  )
    \overline{\lceil}\chi >0\}  .
\]
\hfill {} \\
{\em 2)} \quad These values, spectrum and its  first and second
projections are named subinaccessibly universal reduced to
\ $\chi $ \  of the level \ $n$\quad below \  $\alpha _{1}$.
\hfill {} \\
{\em 3)}\quad If \ $( \alpha ,\Delta  ) \in
\widetilde{\mathbf{S}}_{n}^{\sin \vartriangleleft \alpha _{1}}
\overline{\overline{\lceil}} \chi$, \ then \ $\alpha $ \ is
named the jump cardinal of
these spectra, while \ $\Delta $ \ is named
their Boolean value reduced to \ $\chi $ \ below \ $\alpha_{1}$.
\hfill {} \\
{\em 4)}\quad The cardinal \  $\alpha _{1}$ \ is named
the carrier of these spectra.
\end{definition}

In a similar way multi-dimensional reduced spectra can be
introduced. \label{c4}
\endnote{
\ p. \pageref{c4}. \ For this purpose all Boolean values of
multi-di\-men\-sio\-nal spectra should be reduced to certain
cardinal \ $\chi$. \ For example two-dimensional spectrum (see
comment 2) ) transforms to its the following reduced form:
\[
    \mathbf{S}_{\varphi, 2 }^{\vartriangleleft \alpha _{1}}
    ( \overrightarrow{a}) \overline{\overline{\lceil}} \chi = \{
    ( \alpha , \beta, \Delta_{\varphi}^{\vartriangleleft
    \alpha_{1}}  ( \alpha ,\overrightarrow{a} ) \overline{\lceil} \chi,
    \Delta_{\varphi_1^{\prime}}^{\vartriangleleft \alpha_{1}}
    (\alpha , \beta, \overrightarrow{a} ) \overline{\lceil} \chi  ):
\]
\[
    : \alpha,\beta <\alpha_{1} \wedge \ \Delta_{\varphi }^{\vartriangleleft \alpha_1}
    (\alpha , \overrightarrow{a}) \overline{\lceil} \chi > 0 \
    \wedge \
    \Delta_{\varphi_1^{\prime} }^{\vartriangleleft \alpha_1}
    (\alpha, \beta, \overrightarrow{a}) \overline{\lceil} \chi > 0  \}.
\]
\vspace{6pt}
} 

\noindent Further it is always assumed that \  $\chi $ \  is
closed under the pair function; if $\chi =k$, \ all mentionings
about \ $\chi $ \ will be dropped.
\\
It is not difficult to derive the following two lemmas from
definitions \ref{3.9.}, \ref{4.1.} and lemmas \ref{3.10.},
\ref{3.11.}:

\begin{lemma}
\label{4.2.}
\[
    dom \left ( \widetilde{\mathbf{S}}_{n}^{\sin \vartriangleleft
    \alpha _{1}}\overline{\overline{\lceil}}\chi \right )
    \subseteq dom  \left (  \widetilde{
    \mathbf{S}}_{n}^{\sin \vartriangleleft \alpha _{1}} \right )
    \subseteq SIN_{n-1}^{<\alpha _{1}} \cap
    dom \Bigl( \mathbf{S}_{n}^{\vartriangleleft \alpha_1 } \Bigr).
\]
\end{lemma}

\vspace{-18pt}

\begin{lemma} \label{4.3.} \hfill {} \\
\hspace*{1em} Let
\[
    \alpha _{2}\in SIN_{n-2}^{<\alpha _{1}} \mathrm{\quad and \quad}
    \alpha_{0}=\sup  ( SIN_{n-1}^{<\alpha _{1}}\cap \alpha_{2} ),
\]
then
\[
    ( \widetilde{\mathbf{S}}_{n}^{\sin \vartriangleleft
    \alpha _{2}} \overline{\overline{\lceil }}\chi  )
    | \alpha_{0} = ( \widetilde{\mathbf{S}}_{n}^{\sin
    \vartriangleleft \alpha _{1}}\overline{\overline{\lceil }}\chi)
    | \alpha_{0} \quad.
\]
\end{lemma}
The following lemma is analogous to lemma \ref{2.5.} and comes
from definitions:

\begin{lemma} \label{4.4.} \hfill {} \\
\hspace*{1em} Let \ $\alpha <\alpha _{1}$,\quad $ \chi \leq k$, \  then: \\
\hfill {} \\
{\em 1)}\quad $ \sup dom  \left (
\widetilde{\mathbf{S}}_{n}^{\sin \vartriangleleft
\alpha_{1}}\overline{\overline{\lceil }}\chi \right) < k$.
\\
\quad \\
\medskip
{\em 2)}\ $\widetilde{\mathbf{S}}_{n}^{\sin \vartriangleleft
\alpha_{1}}\overline{\overline{\lceil }}\chi $, \ $dom  \left (
\widetilde{ \mathbf{S}}_{n}^{\sin \vartriangleleft \alpha
_{1}}\overline{\overline{\lceil}}\chi  \right ) $ are
$\Delta_{n}$-definable, while
\\
$rng  \left (\widetilde{\mathbf{S}}_{n}^{\sin \vartriangleleft \alpha
_{1}}\overline{\overline{\lceil}}\chi   \right )$ \ is \ $ \Sigma
_{n}$-definable \medskip in  \  $L_{k}$ \  for \  $\alpha _{1}=k$.
\\
For \  $\alpha _{1}<k$ \ these spectra are
$\Delta_{1}$-definable ;
\\
\quad \\
{\em 3)}\quad $\alpha \in dom  \left (
\widetilde{\mathbf{S}}_{n}^{\sin \vartriangleleft \alpha
_{1}}\overline{\overline{\lceil }}\chi  \right ) $ \ iff there
exists an \ $\mathfrak{M}$-generic function
\[
    l\stackrel{\ast }{\in }\widetilde{\Delta }_{n}^{\sin
    \vartriangleleft \alpha _{1}} (\alpha) \overline{\lceil}\chi ;
\]
\hfill {} \\
{\em 4)}\quad \quad \quad \quad \quad \quad \quad  $\left\|
\widetilde{u}_{n}^{\sin \vartriangleleft \alpha_{1}}
( \underline{l} ) \right\| \lceil \chi =\sum rng  \left ( \widetilde{%
\mathbf{S}}_{n}^{\sin \vartriangleleft \alpha _{1}}\overline{\overline{%
\lceil }}\chi  \right ) $.
\end{lemma}
As in lemma \ref{2.5.}~3) statement 3) here makes it possible to
discover jump ordinals \ $\alpha $ \  with the help of generic
functions \  $l$; \ this techniques is used below in the proof of
lemma \ref{4.6.}. But it is possible to get along without it
using instead of \  $l$ \  conditions \  $p\subset l$ \  with a
sufficiently long domain.
\\

Let us turn to the discussion of order spectrum types. If \ $X$
\  is a well ordered set, then its order type is denoted by \ $OT
( X ) $; \  if \  $X$ \  is a function having well ordered
domain, then
we assume \  $OT ( X ) =OT ( dom ( X )  ) $.\\

Rough upper estimate of spectrum types comes from lemma
\ref{1.1.}, \ $\left| P_{\chi }\right| =\left| \chi \right| $ \
and \  $GCH$ \ in \ $L_{k}$~:
\begin{lemma} \label{4.5.}
\[
    OT ( \widetilde{\mathbf{S}}_{n}^{\sin \vartriangleleft
    \alpha_{1}}\overline{\overline{\lceil }}\chi  ) < \chi ^{+}\ .
\]
\end{lemma}

Now let us discuss estimates of such types from below. Here comes
out the lemma essential for the proof of main theorem. It shows,
that as soon as an ordinal \ $ \delta < \chi ^{+}$ \ is defined
through some jump ordinal of the subinaccessibly universal
spectrum reduced to \ $\chi$, \ the order type of this spectrum
exceeds  \ $\delta$ \ under certain natural conditions.

We shall use here and further the method of reasoning that may be
named \textit{splitting method.}
In its simplest version it consists in splitting some value
\[
    \Delta = \widetilde{\mathbf{S}}_{n}^{\sin \triangleleft
    \alpha_1} \overline{\overline{\lceil}} \chi (\alpha)
\]
(or several such values) under consideration in a sequence of its
parts and in assigning them to successive \ $SIN_{n-1}^{<
\alpha_1} $-cardinals respectively after their certain slight
transformation. After that these cardinals become \textit{jump
cardinals} of the spectrum. For this purpose beforehand the
cardinal \  $\alpha $ \ must be \textit{fixed} by some  \  $l $, \
that is an \ $\mathfrak{M} $-generic function \ $l\stackrel{\ast
}{\in } \Delta $ \  must be used.

For some convenience the suitable notation of ordinal intervals
will be used for \ $\alpha _{1}<\alpha _{2}$:
\\
$\left[ \alpha _{1},\alpha _{2}\right[ =\alpha _{2}-\alpha _{1}$;
\quad $\left] \alpha _{1},\alpha _{2}\right[ =\alpha _{2}- (
\alpha _{1}+1 ) $; \quad $\left[ \alpha _{1},\alpha _{2}\right] =
( \alpha _{2}+1 ) -\alpha _{1}$; \quad $\left] \alpha _{1},\alpha
_{2}\right] = ( \alpha _{2}+1 ) - (\alpha _{1} + 1)$ \quad (here \
$\alpha _{1}, \alpha _{2}$ \ are sets of smaller ordinals).

\begin{lemma}
\label{4.6.} {\em (About spectrum type)} \hfill {} \\
\hspace*{1em} Let ordinals \ $ \overline{\delta
}$,~$\overline{\chi }$,~$\overline{\alpha
}_{0}$,~$\overline{\alpha}_1$ \ be such that:
\\
\quad \\
(i) \quad $\overline{\delta }<\overline{\chi
}^{+}<\overline{\alpha}_{0}<\overline{\alpha}_1\leq k$~; \\
\quad \\
(ii) \quad $SIN_{n-2} ( \overline{\alpha}_1 ) \wedge OT (
SIN_{n-1}^{<\overline{\alpha}_1 } ) =\overline{\alpha}_1$~;
\\
\quad \\
(iii) \quad $SIN_{n-1}^{<\overline{\alpha}_1 } ( \overline{\chi }
) \wedge \sup dom  \left ( \widetilde{\mathbf{S}}_{n}^{\sin
\vartriangleleft
\overline{\chi }}  \right ) =\overline{\chi }$~; \\
\quad \\
(iv) \quad $\sum rng  \left ( \widetilde{\mathbf{S}}_{n}^{\sin
\vartriangleleft \overline{\chi }}  \right ) \in
B_{\overline{\chi }}$~;\\
\quad \\
(v) \quad $\overline{\alpha }_{0}\in dom  \left (
\widetilde{\mathbf{S}} _{n}^{\sin \vartriangleleft
\overline{\alpha}_1}\overline{\overline{\lceil }}
\overline{\chi}  \right ) $~;\\
\quad \\
(vi) \quad $\overline{\delta }$ \ is defined in  \ $ L_{k} $ \
through ordinals \  $\overline{\alpha }_{0}$ , $\overline{\chi }$
\  by a formula of the class \  $\Sigma
_{n-2}\cup \Pi _{n-2}$~.\\
\quad \\
Then $\quad \quad \quad \quad \quad \quad \quad \quad \quad
\quad\overline{\delta }<OT ( \widetilde{\mathbf{S}}_{n}^{\sin
\vartriangleleft \overline{\alpha}_1}\overline{\overline{\lceil }}%
\overline{\chi}  ) \quad.$
\end{lemma}
\textit{Proof.} \ Let us introduce the following formulas
describing the essential aspects of the situation below \
$\overline{\alpha}_1$.
\\
By condition $(vi)$ there exists the  \  $\Sigma_{n-2}\cup
\Pi_{n-2}$-formula \  $\psi _{0} ( \alpha _{0},\chi ,\delta  ) $ \
which defines \  $\overline{\delta } $ \ through \  $\alpha _{0}
=\overline{\alpha }_{0}$, $\chi =\overline{\chi }$, \ that is \
$\delta =\overline{\delta }$ \  is the only ordinal satisfying \
$\psi _{0} ( \overline{\alpha }_{0},\overline{\chi },\delta  ) $ \
in \  $L_{k}$.
\\
Due to  $(iii)$ and lemma~3.5 and the minimality of \ $k$ \ \
$\overline{\chi }$ \ is the singular cardinal and then by
lemma~1.3~2) for every \ $\mathfrak{M}$-generic function \ $\left.
l\right| \overline{\chi } $ \ on \ $\overline{ \chi }$ ordinals \
$\overline{\chi }$, \ $\overline{\delta }$ \ are countable in \
$L_{k}\left[ \left. l\right| \overline{\chi }\right] $. \ Let us
denote by\quad $ ( f:\omega _{0}\longrightarrow \delta +1 ) $ \
the formula
\begin{center}
{``$f$ \ maps \ $\omega _{0}$ \  onto \ $ \delta +1$''~~.}
\end{center}
In algebra \ $B_{\overline{\chi }}$ \ it has the Boolean value
\[
    \left\| \exists f\  ( f:\omega _{0}\longrightarrow
    \delta +1 ) \right\|_{\overline{\chi }}=1
\]
and therefore there exists some name \ $\underline{f }\in
L_{k}^{B_{\overline{\chi }}}$ \ for which
\[
    \| \underline{f}:\omega _{0}\longrightarrow
    \delta +1 \| _{\overline{\chi } }=1
\]
in \ $B_{\overline{\chi}}$ \ (see~\cite{Jech}). Hence there is the
ordinal \ $\overline{\beta }$ \ defined through \ $\alpha _{0}=
\overline{\alpha }_{0}$, \ $\chi =\overline{\chi }$ \ by the
following formula which is denoted by \ $\psi_{1} ( \alpha
_{0},\chi ,\beta  ) $:
\\
\quad \\
\hspace*{1em} $ \exists \delta<\alpha _{0}( ^{^{{}}}\psi _{0} (
\alpha _{0},\chi ,\delta ) \wedge$
$$
\qquad \wedge \beta =\min \{ \beta ^{\prime }:F ( \beta ^{\prime }
) \in L^{B_{\chi }}\wedge \left\| F ( \beta ^{\prime } ) :\omega
_{0}\longrightarrow \delta +1\right\| _{\chi }=1\}  ) .
$$
From conditions $(ii)$, $(v)$ and lemmas~\ref{3.8.},~\ref{4.2.} it
follows that \ $\overline{\alpha }_{0} \in SIN_{n-2}$ \  and then
\ $\overline{\beta }<\overline{\alpha}_{0}$. \ Next from $(ii)$,
$(v)$ and the same lemma~\ref{4.2.} it follows that \  $
\overline{\alpha }_{0}\in dom  \left (
\widetilde{\mathbf{S}}_{n}^{\sin \vartriangleleft \overline{\alpha
}_1}  \right ) $ \  and it makes possible to define the jump
cardinal
\[
    \overline{\alpha}_0^\prime = \min  ( \left[ \overline{\chi },
    \overline{\alpha }_{0}\right] \cap dom  \left
    ( \widetilde{\mathbf{S}}_{n}^{\sin \vartriangleleft
    \overline{\alpha}_1}  \right ) ) .
\]
Since \ $\overline{\chi },\overline{\alpha}_0^\prime \in
SIN_{n-1}^{<\overline{\alpha}_1}$ \ from $(iii)$, $(iv)$ it
follows that
\[
    \overline{\alpha}_0^\prime \in dom  \left
    ( \widetilde{\mathbf{S}}_{n}^{\sin \vartriangleleft
    \overline{\alpha}_1} \overline{\overline{\lceil }}\overline{\chi}
    \right )
\]
and the cardinal \  $\overline{\chi }$ \  is definable through \
$\alpha_0^\prime = \overline{\alpha}_0^\prime$ \ by the formula
\[
    \chi =\sup dom  \left ( \widetilde{\mathbf{S}}_{n}^{\sin
    \vartriangleleft \alpha_0^\prime } \right )
\]
which we denote by \ $\psi _{2} ( \alpha_0^\prime , \chi).$
\\
Let us assume that formulas \ $\psi _{i}$,\quad
$i=\overline{0,2}$, \  are transformed in the preuniversal form;
further they will be used in generic extension of \ $L_{k}$ \ and
at that time they will be bounded by the
constructive class and denoted by \  $\psi _{i}^{L}$, $i=\overline{0,2}$.

Now the \textit{splitting method} starts to work. To this end the
notion of direct product of functions is needed. Let \
$l_{0},...,l_{m}\in {}^{k}k$; \ on \ $k$ \ there is defined the
function \ $l=l_{0}\oplus ...\oplus l_{m}$ \ in the following way:
for every ordinals \  $\alpha <k, \quad \alpha =\alpha _{0}+ ( m+1
) i+j$ , \ where the ordinal \ $\alpha _{0}$ \  is limit and \
$i\in \omega _{0}$, \quad $j \leq m $, \ it has the value \ $l (
\alpha  ) =l_{j} ( \alpha _{0}+i ) $ . \ Let us denote by \ $ ( ,
) _{j}^{m}$ \ the operation reconstructing \ $l_{j}$ \ by \ $l$,
that is \  $l_{j}= ( l ) _{j}^{m}$. \ It is known (see Solovay
~\cite{Solovay}) that \ $l$ \ is the \ $\mathfrak{M}$-generic
function on \ $k$ \ \textit{iff} \ $l_{0}$ \  is the \
$\mathfrak{M}$-generic function and \ $l_{j}$ \ is the \
$\mathfrak{M}\left[ l_{0}\oplus ...\oplus l_{j-1}\right] $-generic
function on \  $k$, \ $j=\overline{1,m}.$
\\
We shall take here \ $\mathfrak{M}$-generic \ $l_0$, \ $l_1$ \ to
hold fixed Boolean values of the spectrum on jump cardinals \
$\overline{\alpha}_0$, \ $\overline{\alpha}_0^\prime$, \ that is
take \ $l_0$, \ $l_1$ \ such that
\[
    l_0\stackrel{\ast }{\in }\widetilde{\Delta }_{n}^{\sin
    \vartriangleleft \overline{\alpha}_1} (\overline{\alpha}_0)
    \overline{\lceil}\chi, \quad \quad
    l_1\stackrel{\ast }{\in }\widetilde{\Delta }_{n}^{\sin
    \vartriangleleft \overline{\alpha}_1} (\overline{\alpha}_0^\prime)
    \overline{\lceil}\chi;
\]
there existence comes from lemma~4.4~3).
\\
Now let us consider  the formula \  $\varphi =\exists \alpha \
\varphi _{1} ( \alpha ,\underline{l} ) $ \  collecting all
information about the situation below \ $\overline{\alpha}_1$ \
obtained and specifying the splitting Boolean values on jump
cardinals \ $\overline{\alpha}_0$, \ $\overline{\alpha}_0^\prime$,
\ where \ $\varphi _{1} $ \ is the following formula:
\begin{eqnarray*}
    SIN_{n-1} ( \alpha  ) \wedge \quad \qquad \qquad \qquad \qquad
    \qquad \qquad \qquad \qquad \qquad  \qquad \qquad
\\
    \wedge \exists \alpha_0 <\alpha \exists \alpha_0^\prime
    \leq \alpha _{0}\exists \chi \leq \alpha_0^\prime \exists \beta
    <\alpha _{0}\exists i\in \omega _{0}\exists \delta <\alpha
    _{0}\exists \delta _{i}\leq \delta \exists p\in P_{\chi }
\\
    \left[ SIN_{n-1}^{<\alpha } ( \alpha _{0} ) \wedge
    SIN_{n-1}^{<\alpha } ( \alpha_0^\prime ) \wedge \widetilde{u}%
    _{n-1}^{\vartriangleleft \alpha } ( \alpha _{0,} ( \underline{l})
    _{0}^{3} ) \wedge \widetilde{u}_{n-1}^{\vartriangleleft \alpha}
    ( \alpha_0^\prime, ( \underline{l} ) _{1}^{3} ) \right. \wedge
\\
    \wedge \psi _{0}^{L} ( \alpha _{0},\chi ,\delta  ) \wedge
    \psi_{1}^{L} ( \alpha_0,\chi ,\beta  ) \wedge \psi _{2}^{L}
    (\alpha _{0}^\prime,\chi  ) \wedge
\\
    \wedge ( \underline{l} )
    _{2}^{3} ( \omega _{0} ) =i\wedge p\subseteq  ( \underline{l}%
     ) _{3}^{3}\wedge p\leq \left\| F ( \beta  )  ( i )
     =\delta _{i}\right\| _{\chi }\wedge
\\
    \left . \wedge OT\{ \alpha ^{\prime }<\alpha :SIN_{n-1}^{<\alpha }
    ( \alpha ^{\prime } ) \} = \alpha_{0} + \delta _{i}^{^{{}}} \right].
\end{eqnarray*}
\hfill {} \\
This formula belongs to \  $\Pi _{n-1}$ \  since all its variables
in square brackets are bounded by the  \ $SIN_{n-1}$-variable \ $
\alpha $. \ For every \  $\delta \leq \overline{\delta }$ \  let
us denote by \  $\alpha _{\delta }$ \ the cardinal \ $\alpha <
\overline{\alpha }_1$ \ such that

\vspace{6pt}
\begin{equation}
\label{e4.1} SIN_{n-1}^{<\overline{\alpha}_1 } ( \alpha  ) \wedge
OT\{ \alpha ^{\prime }<\alpha :SIN_{n-1}^{<\alpha } ( \alpha
^{\prime } ) \} = \overline{\alpha}_{0} + \delta ~;
\end{equation}
\vspace{0pt}

\noindent this cardinal does exist due to $(ii)$ and lemma~3.5.
Let us show that every \  $\alpha _{\delta }$ \ is the jump
cardinal of \  $\varphi $ \  below \  $ \overline{\alpha }_1$ \ by
means of \textit{splitting} Boolean values mentioned above. To
this end let us assume \ $\mathfrak{M}$-generic function on \ $k$
\ \ $l=l_{0} \oplus l_{1} \oplus l_{2} \oplus l_{3}$ \ such that

\vspace{6pt}
\begin{equation}
\label{e4.2} L_{k}\left[ l_{0}\right] \vDash
\widetilde{u}_{n-1}^{\vartriangleleft \overline{\alpha}_1} (
\overline{\alpha }_{0,}l_{0} )\quad ;
\end{equation}

\begin{equation}
\label{e4.3}
l_{0}\stackrel{\ast }{\notin }\left\| \exists \alpha <\overline{\alpha }%
_{0} ( SIN_{n-1}^{\vartriangleleft \overline{\alpha}_1} ( \alpha )
\wedge \widetilde{u}_{n-1}^{\vartriangleleft \overline{\alpha}_1 }
( \alpha _{,}\underline{l} )  ) \right\| \lceil \overline{\chi }
\quad;
\end{equation}

\begin{equation}
\label{e4.4} L_{k}\left[ l_{1}\right] \vDash
\widetilde{u}_{n-1}^{\vartriangleleft \overline{\alpha}_1} (
\overline{\alpha}_0^\prime, l_{1} )\quad ;
\end{equation}

\begin{equation}
\label{e4.5} l_{1}\stackrel{\ast }{\notin }\left\| \exists \alpha
<\overline{\alpha}_0^\prime ( SIN_{n-1}^{\vartriangleleft
\overline{\alpha}_1} ( \alpha ) \wedge
\widetilde{u}_{n-1}^{\vartriangleleft \overline{\alpha}_1 } (
\alpha _{,}\underline{l} )  ) \right\| \lceil \overline{\chi
}\quad;
\end{equation}

\begin{equation}
\label{e4.6} l_{2} ( \omega _{0} ) =i_{\delta }\quad;
\end{equation}

\begin{equation}
\label{e4.7} \exists p\in P_{\overline{\chi }} ( p\subset
l_{3}\wedge p\leq \left\| F ( \overline{\beta } )  ( i_{\delta }
) =\delta \right\| _{\overline{\chi }} )
\end{equation}
\vspace{0pt}

\noindent for corresponding \ $i_{\delta } \in \omega_0$. \ Here
the existence of \ $l_{2}$, \ $l_{3}$ \ is obvious. From
(\ref{e4.2}), (\ref{e4.4}), (\ref{e4.6}), (\ref{e4.7}) and
(\ref{e4.1}) it comes

\vspace{6pt}
\begin{equation}
\label{e4.8} l\stackrel{\ast }{\in }\left\| \varphi
_{1}^{\vartriangleleft \overline{\alpha}_1} ( \alpha _{\delta
},\underline{l} ) \right\| \lceil \overline{\chi }\quad.
\end{equation}
\vspace{0pt}

\noindent Then from (\ref{e4.3}), (\ref{e4.5}) it is not difficult
to deduce that

\[
    l\stackrel{\ast }{\notin }\left\| \exists \alpha <
    \alpha_{\delta }\ \varphi_{1}^{\vartriangleleft
    \overline{\alpha}_1} ( \alpha ,\underline{l} )
     \right\| \lceil \overline{\chi }\quad.
\]
\vspace{0pt}

\noindent Along with (\ref{e4.8}) it implies \quad
\(
    \alpha_{\delta } \in dom \left (
    \mathbf{S}_{\varphi }^{\vartriangleleft \overline{\alpha}_1}
    \overline{\overline{\lceil} }\overline{\chi } \right )
\).
\ Then with the help of condition~$(ii)$ it is not hard to see
that for every \
\(
    \delta \leq \overline{\delta }\quad dom \left (
    \widetilde{\mathbf{S}}_{n}^{\sin \vartriangleleft
    \overline{\alpha}_1} \overline{\overline{\lceil }}
    \overline{\chi }\right )
\)
\ contains the cardinal succeeding \ $\alpha _{\delta}$ \ in \
$SIN_{n-1}^{<\overline{\alpha}_1}$. \hfill $\dashv$
\\

Complicating this reasoning insignificantly one can prove this
lemma in the case when  \  $\delta $ \  is defined through  \
$\overline{\chi} $ \  and several jump cardinals
\[
    \overline{\alpha}_{1}, \dots, \overline{\alpha}_{m} \in dom  \left
    (\widetilde{\mathbf{S}}^{\sin \triangleleft \overline{\alpha}_1}
    \overline{\overline{\lceil}} \overline{\chi} \right ) .
\]

\noindent Lemma \ref{4.6.} admits various strengthenings and
versions but they are indifferent for what follows and therefore
are omitted. \label{c5}
\endnote{
\ p. \pageref{c5}. \ For instance, it is possible to require the
definability of \ $\overline{\delta }$ \ not in \  $ L_{k} $, \
but in \ $L_{k}\left[ l\right] $; \ it is possible to weaken
condition (ii) up to the condition

\[
    SIN_{n-2} ( \overline{\alpha}_1 ) \wedge OT ( \;
    ]\overline{\alpha }_0, \overline{\alpha}_1 [ \; \cap \;
    SIN_{n-1}^{< \overline{\alpha}_1 } ) \geq \overline{\chi}^{+} ,
\]

\noindent also the condition imposed on \ $\overline{\chi}$ \ can
be weakened substantially, etc.
\\
\quad \\
} 
\hfill {} \\
\hfill {} \\
In this reasoning the formulas interpreted in \ $L_{k}$ \  were
equipped with the upper index \  $L$ \  when turning to
extensions of \ $L_{k}$. \  Further we shall omit this index if
the corresponding passage is meant by the context.
\hfill {} \\
\hfill {} \\

Lemma~4.6 about spectrum type here is the substantially important
property of reduced spectra along with their informative
properties (see lemma~5.1 below).

\noindent But still there is the following essential
inconvenience: such spectra, taken on their different carriers,
can be hardly compared with each other in view to their basic
properties, because their domains can contain an arbitrary great
cardinals, when these carriers are increasing up to \ $k$. In
order to avoid this obstacle we shall transform them to reduced
matrices.

\newpage
\quad 

\newpage

\section{Reduced Matrices}
\setcounter{equation}{0}

Now we start to form matrix functions. For this purpose, with
reduced spectra in hand, we turn here to their simple
transformation -- reduced matrices, being values of such
functions. These matrices comes from reduced spectra by easy
isomorphic enumeration of their domains:

\begin{definition}
\label{5.1.}
\ \\
{\em 1)}\quad We name as a matrix reduced to an ordinal \ $\chi $
 \  every relation \  $M$ \  satisfying the following
formula denoted by \  $\mu  ( M,\chi  )$ :
\[
    (M \mbox{ is a function})\wedge (dom ( M )
    \mbox{ is an ordinal})\wedge rng ( M ) \subseteq B_{\chi }\quad .
\]

\noindent {\em 2)}\quad Let \  $M$ \  be a matrix and \
$M_{1}\subset k\times B$. We name as a superimposition of \
$M$ \ onto \ $M_{1}$ \  a function \  $f$  \  satisfying the
following formula denoted by \  $f:M\Rightarrow M_{1}$~:
\begin{eqnarray*}
& ( f \ \mbox{ is an order isomorphism } of \ dom (M ) \
{onto}\  \ dom ( M_{1} )  ) \wedge \qquad\qquad
\\
& \wedge \forall \alpha ,\alpha ^{\prime }\forall \Delta ,\Delta
^{\prime } ( f ( \alpha  ) =\alpha ^{\prime }\wedge  ( \alpha
,\Delta  ) \in M\wedge  ( \alpha ^{\prime },\Delta ^{\prime } )
\in M_{1}\longrightarrow \qquad
\\
& \qquad\qquad\qquad\qquad\qquad\qquad\qquad\qquad\qquad
\longrightarrow \Delta =\Delta ^{\prime } ).
\end{eqnarray*}

\noindent If the superimposition exists then we say that \  $M$
 \  is superimposed onto \  $M_{1}$ \  and write \
$M\Rightarrow M_{1} $.

\noindent {\em 3)}\quad If the matrix \  $M$ \  superimposes onto
the spectrum  \ $\widetilde{\mathbf{S}}_{n}^{\sin \vartriangleleft
\alpha }\overline{\overline{\lceil }}\chi $ then \ $M$ \  is
named the matrix of this spectrum on \  $\alpha $, \ or the
subinaccessibly universal matrix of the level \  $n$ \ reduced to \
$\chi$ \ on \ $\alpha $.
\\
{\em 4)}\quad In this case if \ $(\alpha^{\prime} ,\Delta ) \in
\widetilde{\mathbf{S}}_{n}^{\sin \vartriangleleft \alpha
}\overline{\overline{\lceil }}\chi $, \ then  \ $\alpha^{\prime}$
\ is named the jump cardinal
of the matrix \ $M$, \ while \ $\Delta$ \ is named
its Boolean value on \ $\alpha$.
\\
{\em 5)}\quad In this case the cardinal \ $\alpha$ \ is
named the carrier of the matrix \ $M$.
\end{definition}

\noindent By this definition if \ $M\Rightarrow
\widetilde{\mathbf{S}}_{n}^{\sin \vartriangleleft \alpha
}\overline{ \overline{\lceil }}\chi$, \ then
\[
    rng (M) = rng (\widetilde{\mathbf{S}}_{n}^{\sin
    \vartriangleleft \alpha }\overline{ \overline{\lceil }}\chi),
\]
so from here we shall consider the Boolean spectrum on \ $\alpha$
\ also as \ $rng (M)$.
\\
\hfill {} \\
Matrices introduced above can be named one-dimensional; following
this definition it is possible to define multi-dimensional
matrices superimposed onto multi-dimensional reduced spectra of
the same dimension (see comment 4) ) with a view to a finer
analysis of propositions of our language. Matrices of this
multi-dimensional kind were used by the author for a long time as
the main tool of investigations of subinaccessibility.
\\
In what follows carriers \  $\alpha $ \  and reducing cardinals \
$\chi $ \  will be $SIN_{n-2}$-cardinals, \  $\chi <\alpha \leq
k$ \ (if some other case is not considered).
\\
\hfill {} \\
From this definition and lemmas \ref{4.4.}, \ref{4.5.} there come
the following two lemmas:

\begin{lemma}
\label{5.2.} \hfill {} \\
\hspace*{1em} The formulas
\[
    f:M\Rightarrow \widetilde{\mathbf{S}}_{n}^{\sin
    \vartriangleleft \alpha }\overline{\overline{\lceil }}\chi
    \quad , \ M\Rightarrow \widetilde{\mathbf{S}}_{n}^{\sin
    \vartriangleleft \alpha }\overline{\overline{\lceil }}\chi
\]
belong respectively to \  $\Pi_{n},\Sigma _{n+1}$ \  for \
$\alpha =k$ \ and to \ $\Delta _{1}$ \ for \ $\alpha <k$.
\end{lemma}

\vspace{-6pt}

\begin{lemma}
\label{5.3.} \hfill {} \\
\hspace*{1em} Let
\[
    M\Rightarrow \widetilde{\mathbf{S}}_{n}^{\sin
    \vartriangleleft \alpha } \overline{\overline{\lceil }}\chi ,
\]
then \medskip

\noindent {\em 1)} \ $\left\| \widetilde{u}_{n}^{\sin
\vartriangleleft \alpha } ( \underline{l} ) \right\| \lceil \chi
=\sum rng ( M )
$; \\

\noindent {\em 2)} \ $OT ( M ) =dom ( M ) \leq Od ( M ) <\chi
^{+}.$
\end{lemma}

\noindent The last statement shows that now reduced matrices can
be compared (in the sense of G\"{o}del function $Od$)
\textit{within \ $L_{\chi^{+}}$ \ only} and this circumstance will
make possible to define matrix functions with required properties.
\\

The main role further is played by matrices and spectra reduced to
complete cardinals; their existence comes out from lemma
\ref{2.5.}~1) (for \ $\widetilde{u}_n^{\sin}$, \ $k$ \ as \
$\varphi$, \ $\alpha_1$ \ respectively):

\begin{definition}
\label{5.4.} \hfill {} \\
\hspace*{1em} We name as a complete ordinal of level \  $n$
 \  every ordinal \  $\chi $ \  such that
\[
    \exists x~\widetilde{u}_{n-1}^{\sin } ( x,\underline{l} )
    \dashv \vdash \exists x<\chi
    ~\widetilde{u}_{n-1}^{\sin}( x,\underline{l} ) .
\]
The least of these  ordinals is denoted by \  $\chi ^{\ast }$, \
while the value \  $\left\| \widetilde{u}_{n}^{\sin }
(\underline{l} ) \right\| $ \  is denoted by \  $A^{\ast}$~.
\end{definition}

\vspace{-6pt}

\begin{lemma}
\label{5.5.} \

\quad

\noindent {\em 1)}\ $\chi ^{\ast }=\sup dom \left (
\widetilde{\mathbf{S}}_{n}^{\sin } \right ) =\sup dom \left (
\mathbf{S}_{n}\right ) <k$ \quad;

\quad

\noindent  {\em 2)}\ $SIN_{n-1} ( \chi ^{\ast } ) ~ , ~ \chi
^{\ast }=\omega _{\chi ^{\ast }}$\quad;

\quad

\noindent {\em 3)}\ $OT ( \chi ^{\ast }\cap SIN_{n-1} ) =OT (
dom  \left ( \widetilde{ \mathbf{S}}_{n}^{\sin } \right ) ) =OT (
dom \left ( \mathbf{S}_{n} \right )) =\chi ^{\ast }$;

\quad

\noindent {\em 4)}\ $\widetilde{\mathbf{S}}_{n}^{\sin
}=\widetilde{\mathbf{S}}_{n}^{\sin \triangleleft \chi^{\ast}} =
\widetilde{\mathbf{S}}_{n}^{\sin }\overline{\overline{\lceil
}}\chi ^{\ast }; \ \ \mathbf{S}_{n}=
\mathbf{S}_{n}^{\triangleleft \chi^{\ast}} =\mathbf{S}_{n}
\overline{\overline{\lceil }}\chi ^{\ast }; \ \ A^{\ast }=A
^{\ast } \lceil \chi ^{\ast }$;

\quad \\
and similarly for the reduction of these spectra to every \
$SIN_{n-1}$-cardinal \ \ $\geq \chi ^{\ast }$.
\end{lemma}

\noindent \textit{Proof} of 1) follows from lemmas 3.10~2), 2.7,
1.1; statement 2) follows from lemmas~\ref{3.10.}~2),
\ref{3.4.}~1), \ref{3.5.} (for the level \ $n-1$) when \ $\alpha
_{1}=k$. Then 3),~4) are deduced by the splitting method of the
proof of lemma~\ref{4.6.}.~\hfill~$\dashv$
\\
One can show that lemma 5.5 2) is best possible and \
$\chi^{\ast} \notin SIN_{n} $.
\\
\hfill {} \\
From this lemma it follows one more important property of the
complete cardinal reinforcing lemma 3.2 (for  \  $\alpha_{1}=k ,
\alpha=\chi^{\ast} $) because here  \  $\chi^{\ast} $ \  is only
the  \  $ SIN_{n-1} $-cardinal:

\begin{lemma}
\label{5.6.} \hfill {} \\
\hspace*{1em} Let \  $\forall x~\varphi $ \  be a \  $\Pi
_{n}^{\dashv \vdash }$-proposition, \  $\varphi \in \Sigma
_{n-1}^{\dashv \vdash }$, \ then:
\\
\\
{\em1)}\quad if \  $\varphi $ \  contains individual constants
only from \ $L_{\chi ^{\ast }}^{B_{\chi ^{\ast }}}$, then for
every \ $\mathfrak{M}$-generic function \  $l$
\[
    L_{k} [ l] \vDash \left( \forall x\vartriangleleft ^{\underline{l}}
    \chi ^{\ast} ~\varphi \longleftrightarrow \forall
    x~\varphi \right)~;
\]
{\em 2)}\quad if \  $\varphi $ \  contains individual constants
$\vartriangleleft \chi^{\ast} $ \  only from \  $L_{k} $ \ and
does not contain \ $\underline{l}$, \  then
\[
    \quad L_{k}\vDash \left( \forall x\vartriangleleft
    \chi ^{\ast }~\varphi \longleftrightarrow \forall
    x~\varphi \right)~;
\]
{\em 3)}\quad let \  $\omega _{0}^{\ast }=\sup dom \left (
\widetilde{\mathbf{S}} _{n}^{\sin }\overline{\overline{\lceil }}
( \omega _{0}+1 ) \right ) $  \  and \  $\varphi $ \  does not
contain individual constants and \  $\underline{l}$, \ then
\[
    \qquad L_{k}\vDash \left( \forall x\vartriangleleft
    \omega _{0}^{\ast}~\varphi \longleftrightarrow
    \forall x~\varphi \right)~.
\]
\end{lemma}
\textit{Proofs} of 1) - 3) are analogous and come to the fact that
all jump cardinals of the proposition \ $\exists x~\neg \varphi $
\ (if they exist) are less than \ $\chi ^{\ast }$, \  while for 3)
even than \ $\omega _{0}^{\ast }$; \ so, they can be demonstrated
first for~2).
\\
Beforehand the following remark should be done:
for every \ $x \in L_{k}$, \  $m \geq 2 $, \ $\alpha\in SIN_{m}$ \
and  \  $\mathfrak{M} $-generic  \  $l $\quad\quad
$x\vartriangleleft \alpha \longleftrightarrow x
\vartriangleleft^{l} \alpha $.
\\
Therefore the restriction  \  $x \vartriangleleft \alpha $ \ should
be considered as  \  $Od(x)<\alpha $ \  over  \  $L_{k} $ \ and
as  \  $Od^{l}(x)<\alpha $ \  over  \  $L_{k}[l]$.
\\
Let us consider the proposition  \  $\forall x \ \ \varphi
(x,\alpha_{0}) $, \ $\varphi \in \Sigma_{n-1} $, \ having only one
individual ordinal constant  \  $\alpha_{0} \vartriangleleft
\chi^{\ast } $ \ (for more clearness) and let  \  $\alpha_{0} \in
dom \left (\widetilde{\mathbf{S}}_{n}^{\sin} \right ) $ \
(otherwise we can use lemma 5.5.~3)~). First the constant  \
$\alpha_{0} $ \ must be fixed as in the proof of lemma 4.6 by some
\ $\mathfrak{M} $-generic  \  $l \stackrel{\ast }{\in }
\widetilde{\bigtriangleup}_{n}^{\sin} (\alpha_{0}) $. Let us
assume that
\[
    L_{k} \vDash \exists x \neg \varphi (x, \alpha_{0}),
\]
then
\[
    L_{k} [l] \vDash \exists \alpha \ \  \varphi_{1} (\alpha,l)
\]
\\
where  \  $\varphi_{1} (\alpha, \underline{l}) $ \  is the
following  \  $\Pi_{n-1}$-formula:
\[
    SIN_{n-1}(\alpha) \wedge \exists \alpha^{\prime} <
    \alpha  \left ( \underline{l} \stackrel{\ast }{\in }
    \widetilde{\bigtriangleup}_{n}^{\sin \triangleleft \alpha}
    (\alpha^{\prime}) \wedge \exists x \vartriangleleft \alpha
    \neg \varphi^{L} (x,\alpha^{\prime}) \right ) .
\]
The proposition   \  $\varphi_{2}(\underline{l}) = \exists \alpha
\ \ \varphi_{1} (\alpha,\underline{l})  $ \  already has no
individual constants and \  $l \stackrel{\ast }{\in }
\bigtriangleup_{\varphi_{2}}(\alpha) $ \  for some its jump
ordinal \ $\alpha $. \ By lemmas~2.7, \ 5.5~1) \ $ \alpha <
\chi^{\ast} $ \ and for some  \  $\alpha^{\prime} < \alpha $
\[
    L_{k} [l] \vDash  \left ( l \stackrel{\ast }{\in }
    \widetilde{\bigtriangleup}_{n}^{\sin \triangleleft \alpha}
    (\alpha^{\prime}) \wedge \exists x \vartriangleleft \alpha \neg
    \varphi^{L} (x,\alpha^{\prime}) \right ).
\]
Since  \  $\alpha \in SIN_{n-1} $ \  and  \  $\alpha^{\prime} <
\alpha $ \  we can drop here the restriction  \  $
\vartriangleleft^{l} \alpha $ \  and hence  \  $ l \stackrel{\ast
}{\in } \widetilde{\bigtriangleup}_{n}^{\sin} (\alpha^{\prime})
$. \ It implies  \  $\alpha^{\prime} = \alpha_{0} $ \  and
\[
    L_{k} [l] \vDash   \exists x \vartriangleleft \alpha \neg
    \varphi^{L} (x,\alpha_{0})
\]
and, at last,
\[
    L_{k} \vDash   \exists x \triangleleft \chi^{\ast}
    \neg \varphi (x,\alpha_{0}) .
\]
\\
Turning to 3) let us assume that  \  $L_{k} \vDash  \exists x
\neg \varphi (x) $, \ then for every \
\mbox{$\mathfrak{M}$-generic} \ $l $
\[
    L_{k} [l] \vDash \exists \alpha \ \ \varphi_{1}^{L} (\alpha)
\]
where
\[
    \ \varphi_{1} (\alpha) = \exists x \vartriangleleft
    \alpha \neg \varphi(x) .
\]
\\
From lemma~1.4 it is clear that
\[
    \left \| \exists \alpha \ \ \varphi_{1}^{L} (\alpha)  \right \| = 1
\]
and for every \ $\alpha$
\[
    \quad \left \|
    \varphi_{1}^{L} (\alpha) \right \| \in  \left \{ 0;1 \right \} .
\]
Let  \  $\mathfrak{n}$ \  be the G\"{o}del number of the formula
\ $\exists \alpha \ \ \varphi_{1}^{L} (\alpha) $ \  and  \
\mbox{$p_{ \mathfrak{n}}= \left \{ (\omega_{0}, \mathfrak{n})
\right \} $}. \  From here and definitions~2.6, 3.9 one can see
that
\[
    p_{ \mathfrak{n}} \leq  \left \| \exists \alpha \ \
    \widetilde{u}_{n-1} (\alpha,\underline{l}) \right \|
\]
and for every   \  $\alpha $
\[
    p_{\mathfrak{n}} \cdot \left \| \widetilde{u}_{n-1}
    (\alpha,\underline{l}) \right \| = 0
    \mathrm{\quad or \quad} p_{\mathfrak{n}}  \leq \left \| \widetilde{u}_{n-1}
    (\alpha,\underline{l}) \right \|~.
\]

These statements are preserved under reducing to  \  $ \chi =
\omega_{0}+1 $ \  the values  \  $\left \| \widetilde{u}_{n-1}
(\alpha,\underline{l}) \right \| $\ . By definition~4.1 it
implies the existence of some jump cardinal \
$\alpha_{\mathfrak{n}} $ \ such that
\[
    p_{\mathfrak{n}}\cdot \widetilde{\bigtriangleup}_{n}^{\sin}
    (\alpha_{\mathfrak{n}}) \overline{\lceil}(\omega_{0}+1)>0 .
\]
Hence
\[
    \alpha_{\mathfrak{n}} \in dom  \left ( \widetilde{\mathbf{S}}_{n}^{\sin}
    \overline{\overline{\lceil}} (\omega_{0}+1)\right ) \mathrm{\
    and\ } \alpha_{\mathfrak{n}}<\omega_{0}^{\ast} .
\]
After that by definitions~3.9, 2.6 \
\[
    L_{k} [l] \vDash \exists \alpha^{\prime} < \alpha_{\mathfrak{n}} \ \
    \varphi_{1}^{L} (\alpha^{\prime})
\]
and, at last,
\[
    L_{k} \vDash \exists x \vartriangleleft
    \omega_{0}^{\ast} \neg \varphi (x).
\]
\hfill $\dashv$
\\
Using here the arguments from the proof of~2) we can admit
individual constants \ $ \in  dom  \left (
\mathbf{\widetilde{S}}_{n}^{\sin} \overline{\overline{\lceil}}
(\omega_{0}+1) \right ) $ \ in
 \ $\forall x \varphi $,  \ but it is not necessary in what follows.
\\
\hfill {} \\
\hfill {} \\
Further the very special role is played by the so
called \textit{singular matrices}:
\begin{definition}
\label{5.7.} \hfill {} \\
\hspace*{1em} We denote by \  $\sigma  ( \chi ,\alpha  ) $ \ the
conjunction of the following formulas:
\\

{\em 1)}\quad $SIN_{n-2} ( \alpha  ) \wedge  ( \chi\mbox{ is
a limit cardinal }<\alpha )$~; \\

{\em 2)}\quad $OT ( SIN_{n-1}^{<\alpha } ) =\alpha $~;\\

{\em 3)}\quad $\sup dom \left ( \widetilde{\mathbf{S}}_{n}^{\sin
\vartriangleleft \alpha }\overline{\overline{\lceil}}\chi\right )
=\alpha $~.
\\
\\
And let \  $\sigma  ( \chi ,\alpha ,M ) $ \  denote the formula
\[
    \sigma  ( \chi ,\alpha  ) \wedge  ( M\Rightarrow
    \widetilde{\mathbf{S}}_{n}^{\sin \vartriangleleft \alpha }
    \overline{\overline{\lceil }}\chi  ) .
\]
The matrix \  $M$  \  and the spectrum \
$\widetilde{\mathbf{S}}_{n}^{\sin \vartriangleleft \alpha}
\overline{\overline{\lceil }}\chi $ \  reduced to \  $\chi $ \ are
named singular on a carrier \ $\alpha $ \  (on an interval
\ $[\alpha _{1},\alpha _{2}[\ $) \ \textit{iff} \ $\sigma
 ( \chi ,\alpha ,M ) $ \  is fulfilled (for some \
$\alpha \in \lbrack \alpha _{1},\alpha _{2}[\ $).
\\
The symbol \ $S$  \  is used for the common notation of singular
matrices.
\end{definition}

In what follows all matrices \textit{will be reduced to certain
cardinal \ $\chi$ \ and singular on their carriers under
consideration}; all reasoning will be conducted in \ $L_{k}$ \ (
or in \ $\mathfrak{M}$ \ if the context does not point to the
opposite case).
\\
Definition~5.7, lemma~5.2 imply easily

\begin{lemma}
\label{5.8.} \hfill {} \\
{\em 1)}\quad The formulas\quad $\sigma  ( \chi ,\alpha  ) $%
, $\sigma  ( \chi ,\alpha ,M ) $\quad belong to\quad $\Pi _{n-2}$.
\hfill {} \\
{\em 2)}\quad If\quad $\sigma  ( \chi ,\alpha  ) $,\quad then
\[
    \widetilde{A}_{n}^{\sin \vartriangleleft \alpha }
    (\chi  ) \lceil \chi < \sum rng \left
    (\widetilde{\mathbf{S}}_{n}^{\sin \vartriangleleft \alpha }
    \overline{\overline{\lceil}}\chi  \right ) =\left\|
    \widetilde{u}_{n}^{\sin \vartriangleleft \alpha }
    ( \underline{l} ) \right\| \lceil\chi .
\]
\end{lemma}

Due to lemma \ref{5.8.}~2) it is possible to introduce the
following important cardinals: \vspace{-6pt}

\begin{definition}
\label{5.9.} \hfill {} \\
\hspace*{1em} Let \ $\sigma  ( \chi ,\alpha, S ) $ \ fulfills,
then we name as jump cardinal
and prejump cardinal after \ $\chi $ \ of the matrix \ $S$ \ on
the carrier \ $\alpha $, \ or, briefly, of the cardinal \
$\alpha$, \ the following cardinals respectively:
\[
    \alpha _{\chi }^{\downarrow }=\min \{ \alpha ^{\prime}
    \in ]\chi ,\alpha \lbrack ~:~ \widetilde{A}_{n}^{\sin
    \vartriangleleft \alpha } ( \chi  ) \lceil \chi
    <\widetilde{A}_{n}^{\sin \vartriangleleft \alpha } (
    \alpha^{\prime } ) \lceil \chi \wedge
\]
\[
    \qquad \qquad \qquad \qquad \qquad \qquad \qquad \qquad \qquad \qquad
    \wedge SIN_{n-1}^{<\alpha } ( \alpha^{\prime } ) \} ;
\]
\[
    \alpha _{\chi }^{\Downarrow }=\sup \{ \alpha ^{\prime}
    <\alpha _{\chi }^{\downarrow }  ~:~
    \widetilde{A}_{n}^{\sin \vartriangleleft \alpha } ( \chi )
    \lceil \chi =\widetilde{A}_{n}^{\sin \vartriangleleft \alpha }
    (\alpha ^{\prime }) \lceil \chi \wedge
\]
\[
    \qquad \qquad \qquad \qquad \qquad \qquad \qquad \qquad \qquad \qquad
    \wedge SIN_{n-1}^{<\alpha } ( \alpha ^{\prime } ) \}.
\]
\end{definition}

In this notation and everywhere further the index \  $\chi $ \
can be omitted if it is arbitrary or can be restored from the
context.

\begin{lemma}
\label{5.10.} \hfill {} \\
\hspace*{1em} Let \  $\sigma  ( \chi ,\alpha, S ) $ \ fulfills,
then the cardinals \  $\alpha _{\chi }^{\downarrow },$ $\alpha
_{\chi }^{\Downarrow }$ \  do exist and

\quad

{\em 1)}~$~\alpha _{\chi }^{\downarrow }=\min \left\{ \alpha
^{\prime }>\chi :\alpha ^{\prime }\in dom   \left (
\widetilde{\mathbf{S}}_{n}^{\sin \vartriangleleft \alpha }
\overline{\overline{\lceil }}\chi \right ) \right\} \quad;$

\quad

{\em 2)}\quad $\alpha _{\chi }^{\Downarrow }<\alpha _{\chi
}^{\downarrow }<\alpha ;\quad ]\alpha _{\chi }^{\Downarrow
},\alpha _{\chi }^{\downarrow }[\cap SIN_{n-1}^{<\alpha
}=\varnothing \quad;$

\quad

{\em 3)}\quad $\alpha _{\chi }^{\downarrow },$ $\alpha _{\chi
}^{\Downarrow }\in SIN_{n-2} \quad$.
\end{lemma}
\textit{Proof} of 1), 2) follows from the definitions and
lemma~3.4, while 3) -- from lemmas~3.6, \ref{3.8.}~2) (for \
$k,\alpha,n-2 $ \ as \ $\alpha_{1},\alpha_{2},n $ \
respectively). \hfill $\dashv$
\\

The following lemma along with lemma 4.6 about spectrum type
constitutes the main tools of treating matrix functions and
presents the important informative property of reduced matrices
and corresponding spectra.
\\
In the foregoing when considering reduced spectra we have mainly
used in our reasoning their first projections, that is ordinal
spectra. As for their second projections or reduced Boolean
spectra or reduced matrices -- these ones have played an auxiliary
role, serving as an instrument of proving lemmas~\ref{4.6.},
\ref{5.5.}, \ref{5.6.} and others. However, second projections of
Boolean spectra or reduced matrices have the following important
characteristic property: they contain information about parts of
the universe bounded by their jump cardinals on their carriers, so
that when such spectrum or matrix is transposed from one carrier
to another the properties of such parts of the universe are
preserved. More precisely:

\begin{lemma}
\label{5.11.} {\em (About matrix informativeness)} \hfill {} \\
\hspace*{1em} Let \ $S$ \  be a matrix reduced to \ $\chi$ \ on
carriers \  $\alpha _{1},\alpha _{2}>\chi $ \ and superimposed on
spectra
\[
    S \Rightarrow \widetilde{\mathbf{S}} _{n}^{\sin \vartriangleleft
    \alpha_{1}}\overline{\overline{\lceil }} \chi, \quad
    S \Rightarrow \widetilde{\mathbf{S}} _{n}^{\sin \vartriangleleft
    \alpha_{2}}\overline{\overline{\lceil }} \chi
\]
on these carriers and
\[
    \overrightarrow{a_{1}}= ( \alpha _{10},\alpha_{11},...,\alpha_{1m}),
    \quad \overrightarrow{a_{2}}= ( \alpha _{20},\alpha_{21},...,\alpha_{2m})
\]
be trains of cardinals from ordinal spectra
\[
    dom  \left ( \widetilde{\mathbf{S}}_{n}^{\sin \vartriangleleft
    \alpha _{1}} \overline{\overline{\lceil }}\chi  \right ),
    \quad
    dom  \left ( \widetilde{\mathbf{S}}_{n}^{\sin \vartriangleleft
    \alpha _{2}} \overline{\overline{\lceil }}\chi  \right )
\]
on these carriers which correspond to the same Boolean values:
\[
    \widetilde{\mathbf{S}}_{n}^{\sin \vartriangleleft
    \alpha _{1}}\overline{\overline{\lceil }}\chi  (
    \alpha _{1i} ) =\widetilde{\mathbf{S}}_{n}^{\sin
    \vartriangleleft \alpha_{2}}\overline{\overline{\lceil }}
    \chi ( \alpha _{2i} ) , \ \ i=\overline{0,m}\mathit{.}
\]
Let, at last, \  $\psi  ( x_{1},x_{2},...,x_{m},\underline{%
l} )$ \ be an arbitrary formula of an arbitrary level with free
variables \ $ x_{1}$, $~x_{2}$,$...$, $x_{m}$ \  without
individual constants. Then

\[
    \psi ^{\vartriangleleft \alpha _{10}} ( \alpha _{11},...,\alpha _{1m},
    \underline{l} ) \dashv \vdash \psi ^{\vartriangleleft \alpha_{20}}
    ( \alpha _{21},...,\alpha _{2m},\underline{l} ) .
\]
\end{lemma}
\textit{Proof} is conducted by the methods  analogous to those
used in the proof of lemma~\ref{4.6.}. Let us consider for more
transparency the case when the formula  \  $\psi $ \ contains no
free variables and  \ $\underline{l} $ \  and  \
$\overrightarrow{a_{1}}$, $\overrightarrow{a_{2}}$ \  consist only
of the jump cardinals after \  $\chi $ \ on $\alpha_{1}$, $\alpha
_{2}$ \ respectively:
\[
    \alpha _{10}=\alpha _{1\chi }^{\downarrow }, \
    \alpha_{20}=\alpha _{2\chi }^{\downarrow }
\]
because just this case is needed further. Let \ $\alpha _{11}$,
$\alpha _{21}$ \  be the succeeding  \ $\alpha _{10}$, $\alpha
_{20}$ \ jump cardinals in the reduced ordinal spectra considered
respectively on the carriers \ $ \alpha _{1}$, $\alpha _{2}$. \
The idea of the proof is the following:

The satisfiability of the proposition \  $\psi ^{\vartriangleleft
\alpha _{10}}$ \  means that the value
\[
    \Delta =\widetilde{\mathbf{S}}_{n}^{\sin \vartriangleleft
    \alpha_{1}}\overline{\overline{\lceil }}\chi
    ( \alpha _{11} )
\]
contains the corresponding condition \  $ p\in P_{\chi }$ \ in
which it is ``encoded''. When \  $S$ \ is carried over from the
carrier \  $\alpha _{1}$ \  to the carrier \  $\alpha _{2}$ \ the
value
\[
    \Delta =\widetilde{\mathbf{S}}_{n}^{\sin \vartriangleleft
    \alpha _{2}}\overline{\overline{\lceil }}\chi (\alpha_{21} )
\]
still contains \  $p$ \  and it means that \
$\psi^{\vartriangleleft \alpha _{20}}$ \  is fulfilled. But it is
more convenient to use instead of  \  $p $ \  some  \  $
\mathfrak{M} $-generic  function \  $l $ \  including  \  $p $. \
The lower index  \  $\chi $ \  in the notation will be dropped.

So, let  \  $\Delta$ \  be the Boolean value of the spectrum  \
$\widetilde{\mathbf{S}}_{n}^{\sin \triangleleft \alpha_{1}}
\overline{\overline{\lceil}} \chi $ \ corresponding to  \
$\alpha_{10} = \alpha_{1}^{\downarrow} $, \ that is
\[
    \Delta = \widetilde{\Delta}_{n}^{\sin \triangleleft
    \alpha_{1}}(\alpha_{10})\overline{\lceil} \chi,
\]
and let some \  $\mathfrak{M} $-generic function \ $l $ \ holds
this  \ $\Delta$ \ fixed, that is  \ $l\stackrel{\ast }{\in }
\Delta $, \ and hence it holds fixed  \ $\alpha_{10} $ \ when
talking about the carrier \ $\alpha_{1} $ \ and holds fixed  \
$\alpha_{20} $ \ when talking about \ $\alpha_{2} $. \ By
definition 4.1

\begin{equation}
\label{e5.1} l \stackrel{\ast }{\in }  \left \|
\widetilde{u}_{n-1}^{\triangleleft
\alpha_{1}}(\alpha_{10},\underline{l}) \right \| \lceil \chi ~;
\end{equation}

\begin{equation}
\label{e5.2} l \stackrel{\ast }{\notin }  \left \| \exists
\alpha^{\prime} < \alpha_{10}  \left ( SIN_{n-1}^{\triangleleft
\alpha_{1}} (\alpha^{\prime}) \wedge
\widetilde{u}_{n-1}^{\triangleleft \alpha_{1}} (\alpha^{\prime},
\underline{l} ) \right ) \right \| \lceil \chi ~.
\end{equation}
\vspace{0pt}

\noindent Since \ $\alpha_{10} \in SIN_{n-1}^{< \alpha_{1}} $ \
the restriction \ $\vartriangleleft \alpha_{1} $ \ in (5.2) can be
replaced with the restriction \ $\vartriangleleft \alpha_{10} $. \
After dropping the reduction to \  $\chi $ \ in (5.1), (5.2) it
comes

\[
    l_{1} \stackrel{\ast }{\in } \left \|
    \widetilde{u}_{n-1}^{\triangleleft \alpha_{1}}(\alpha_{10},
    \underline{l}) \right \| ~;
\]
\[
    l_{1} \stackrel{\ast }{\notin } \left \|\exists \alpha^{\prime} <
    \alpha_{10}  \left ( SIN_{n-1}^{\triangleleft \alpha_{10}}
    (\alpha^{\prime}) \wedge \widetilde{u}_{n-1}^{\triangleleft
    \alpha_{10}}(\alpha^{\prime}, \underline{l}) \right ) \right \| ~,
\]
\vspace{0pt}

\noindent where  \  $l_{1} $ \  is some \ $\mathfrak{M} $-generic
function coinciding with \ $l $ \ on  \  $\chi $; \ but for some
shortness let us use the present symbol \ $l$. \ Then it means
that

\begin{equation}
\label{e5.3} L_{k}[l] \vDash \widetilde{u}_{n-1}^{\triangleleft
\alpha_{1}}(\alpha_{10},l) \wedge \neg \exists \alpha^{\prime} <
\alpha_{10}  \left ( SIN_{n-1}^{\triangleleft
\alpha_{10}}(\alpha^{\prime}) \wedge
\widetilde{u}_{n-1}^{\triangleleft \alpha_{10}}(\alpha^{\prime},l)
\right )~.
\end{equation}
\vspace{0pt}

\noindent Thus the cardinal  \  $\alpha_{10} $ \  is defined in
 \  $ L_{k}[l] $ \  through  \  $l $ , \ $\alpha_{1}
$ \  and similarly \ $\alpha_{20} $ \ is defined through  \ $l $,
\ $\alpha_{2} $
 \  in  \  $ L_{k}[l] $ \  also. Now let us consider the
proposition  \  $\varphi(\underline{l}) = \exists \alpha \ \
\psi_{1}(\alpha, \underline{l}) $ \  where  \ $\psi_{1}(\alpha,
\underline{l}) $ \  is the following  \  $\Pi_{n-1}$-formula:

\vspace{6pt}
\[
    SIN_{n-1}(\alpha) \wedge \widetilde{u}_{n-1} (\alpha,
    \underline{l}) \wedge \neg \exists \alpha^{\prime} < \alpha
    \left( SIN_{n-1}^{< \alpha} (\alpha^{\prime}) \wedge
    \widetilde{u}_{n-1}^{\triangleleft \alpha}
    (\alpha^{\prime},\underline{l}) \right ) \wedge
    \psi^{\triangleleft \alpha} ~.
\]
\vspace{0pt}

\noindent From (5.3) it is clear that \vspace{6pt}
\begin{equation}
\label{e5.4} L_{k}[l]\vDash \psi^{\triangleleft \alpha_{10}}
\longleftrightarrow L_{k}[l] \vDash \exists \alpha < \alpha_{1} ~
\psi_{1}^{\triangleleft \alpha_{1}} (\alpha, l) ~.
\end{equation}
\vspace{0pt}

After that the argument from the proof of lemma 2.7 should be
repeated: the function  \  $l $ \  must be replaced with the
function \ $l_{0} $ \  as in (2.1) and at the same time the
formula  \ $\varphi $ \  must be transformed to the formula
 \  $\varphi_{2} $ \  just as it was done in this proof, that
is by replacing its subformulas of the kind  \
$\underline{l}(t_{1})=t_{2} $ \ with the subformulas (2.2) and so
on. Let  \  $\mathfrak{n}$ \  be the G\"{o}del number of  \
$\varphi_{2}$ \  and  \  $l_{0}(\omega_{0}) = \mathfrak{n}$. \ As
a result we turn from  \  $\varphi $ \   to the spectral universal
formula and (5.4) implies \vspace{6pt}
\[
    L_{k}[l]\vDash \psi^{\triangleleft \alpha_{10}}
    \longleftrightarrow L_{k}[l_{0}] \vDash \exists x
    \vartriangleleft \alpha_{1} ~ u_{n-1}^{\triangleleft
    \alpha_{1}}(x, l_{0}) ~.
\]
\vspace{0pt}

\noindent Using here  \  $\widetilde{u}_{n-1}^{\sin} $ \ instead
of  \  $u_{n-1}$ \  one can see from definition 5.7 that
\vspace{6pt}
\begin{equation}
\label{e5.5} L_{k}[l]\vDash \psi^{\triangleleft \alpha_{10}}
\longleftrightarrow L_{k}[l_{0}] \vDash \exists x \vartriangleleft
\alpha_{1} \widetilde{u}_{n-1}^{ \sin \triangleleft \alpha_{1}} ~
(x, l_{0}) \ .
\end{equation}
\vspace{0pt}

\noindent Analyzing the construction of the formula \
$\widetilde{u}_{n}^{\sin} $ \  it is not hard to derive from
(5.5), (5.3) that \vspace{6pt}
\[
    L_{k}[l]\vDash \psi^{\triangleleft \alpha_{10}}
    \longleftrightarrow l_{0} \stackrel{\ast }{\in }
    \widetilde{\Delta}_{n}^{\sin \triangleleft
    \alpha_{1}} (\alpha_{11}) ~;
\]
\vspace{0pt}

\noindent let us remind that here  \  $\alpha_{11} $ \  is the
successor of  \  $\alpha_{10} $ \  in the reduced ordinal spectrum
on the carrier  \  $\alpha_{1} $. \ In this reasoning the values \
$l_{0}(\alpha) $ \  for  \  $\alpha \geq \chi $ \  were not used
and that is why \vspace{6pt}
\begin{equation}
\label{e5.6} L_{k}[l]\vDash \psi^{\triangleleft \alpha_{10}}
\longleftrightarrow l_{0} \stackrel{\ast }{\in }
\widetilde{\Delta}_{n}^{\sin \triangleleft \alpha_{1}}
(\alpha_{11}) \overline{\lceil} \chi \quad.
\end{equation}
\vspace{0pt}

\noindent The matrix  \  $S $ \  is superimposed on  \
$\widetilde{\mathbf{S}}_{n}^{\sin \triangleleft \alpha_{2}}
\overline{\overline{\lceil}} \chi $ \  and hence \vspace{6pt}
\begin{equation}
\label{e5.7} L_{k}[l]\vDash \psi^{\triangleleft \alpha_{10}}
\longleftrightarrow l_{0} \stackrel{\ast }{\in }
\widetilde{\Delta}_{n}^{\sin \triangleleft \alpha_{2}}
(\alpha_{21}) \overline{\lceil} \chi ~.
\end{equation}
\vspace{0pt}

\noindent This reasoning up to (5.6) but with  \  $\alpha_{1},
\alpha_{10},\alpha_{11} $ \ replaced respectively with \
$\alpha_{2}, \alpha_{20},\alpha_{21} $ \ establishes that
\vspace{6pt}
\[
    L_{k}[l]\vDash \psi^{\triangleleft \alpha_{20}}
    \longleftrightarrow l_{0} \stackrel{\ast }{\in }
    \widetilde{\Delta}_{n}^{\sin \triangleleft
    \alpha_{2}} (\alpha_{21}) \overline{\lceil} \chi ~.
\]
\vspace{0pt}

\noindent It remains to compare this equivalence with (5.7).
\\
In the case when  \  $\psi $ \  contains free variables  \
$x_{1}...,x_{m} $ \  and  \  $\underline{l} $ \  one should carry
out the analogous reasoning considering several  \  $\mathfrak{M}
$-generic functions
\[
    l_{i} \stackrel{\ast }{\in }
    \widetilde{\Delta}_{n}^{\sin \triangleleft
    \alpha_{1}} (\alpha_{1i}) \overline{\lceil} \chi ,
    \quad i=\overline{0,m} ,
\]
and treating their direct product as in the proof of lemma 4.6.
\hfill $\dashv$
\hfill {} \\
\\
Changing this reasoning slightly it is not hard to conduct it for
the prejump cardinals  \
$\alpha_{10}=\alpha_{1\chi}^{\Downarrow}, \
\alpha_{20}=\alpha_{2\chi}^{\Downarrow}$~.
\\

\noindent The basic instruments of the proof of main theorem are
matrix functions that are sequences of reduced singular matrices
of the special kind. The following lemma makes it possible to
build such functions:

\begin{lemma}
\label{5.12.}
\[
    \forall \chi \forall \alpha _{0} (  ( \chi
    \mbox{\ is a limit cardinal }>\omega _{0} ) \rightarrow
    \exists \alpha _{1}>\alpha _{0}~\sigma  ( \chi ,\alpha _{1} )  ) .
\]
\end{lemma}

\noindent \textit{Proof.} \ This \  $\Pi _{n}$-proposition does
not contain any individual constants or \  $\underline{l}$ \  and
therefore by lemma \ref{5.6.} 3) it is sufficient to prove that
it is fulfilled when the variables \  $\chi $, $ \alpha _{0}$ \
are bounded by the cardinal \  $\omega _{0}^{\ast }$.
\\
So, let \  $\chi $, $\alpha _{0}<\omega _{0}^{\ast }$; \  let us
consider
\[
    \alpha _{1}=\sup dom  \left (
    \widetilde{\mathbf{S}}_{n}^{\sin }\overline{
    \overline{\lceil }}\chi \right ).
\]
One can see that \ $\alpha _{1}\geq \omega _{0}^{\ast }$ \  and
that is why \  $\chi $, $\alpha _{0}<\alpha _{1}$. From lemmas
\ref{4.2.}, \ref{3.4.} (for \  $\alpha _{1}=k$) \  it comes \
$\alpha _{1}\in SIN_{n-1}$ \  and hence conditions 1), 3) of
definition \ref{5.7.} are carried out; condition 2) can be deduced
by the method which is the simplified variant of the splitting
method of the proof of lemma~\ref{4.6.}. Really, suppose that, on
the contrary, \vspace{6pt}
\[
    \alpha_{2} = OT  \left ( SIN_{n-1}^{< \alpha_{1}} \right )
    < \alpha_{1} ~,
\]
\vspace{0pt}

\noindent then there exists the cardinal

\[
    \alpha_{3} \in [ \alpha_{2}, \alpha_{1} [ \cap dom
    \left ( \widetilde{\mathbf{S}}_{n}^{\sin}
    \overline{\overline{\lceil}} \chi \right )
\]

\noindent and some \  $\mathfrak{M} $-generic function  \ $l
\stackrel{\ast }{\in } \widetilde{\Delta}_{n}^{\sin} (\alpha_{3})
\overline{\lceil} \chi $. \ Now let us consider the proposition
\  $\varphi(\underline{l}) = \exists \alpha \ \
\varphi_{1}(\alpha, \underline{l}) $ \  where  \
$\varphi_{1}(\alpha, \underline{l} ) $ \  is the  \ $\Pi_{n-1}
$\- formula:

\begin{eqnarray*}
&   SIN_{n-1}(\alpha)\wedge
\qquad\qquad\qquad\qquad\qquad\qquad\qquad\qquad\qquad\qquad
\\
& \wedge \exists
\alpha^{\prime},\alpha^{\prime\prime}<\alpha \ \ (
\alpha^{\prime}<\alpha^{\prime\prime} \wedge
SIN_{n-1}^{<\alpha}(\alpha^{\prime\prime})\wedge
\widetilde{u}_{n-1}^{\sin \triangleleft \alpha^{\prime\prime}}
(\alpha^{\prime},\underline{l})\wedge
\\
&   \wedge \neg \exists \alpha^{\prime\prime\prime} <
\alpha^{\prime} \ \ \widetilde{u}_{n-1}^{\sin \triangleleft
\alpha^{\prime\prime}}
(\alpha^{\prime\prime\prime},\underline{l}) \wedge
\\
& \qquad\qquad\qquad\qquad \wedge OT \left (
\left \{ \alpha^{\prime\prime\prime}< \alpha^{\prime\prime} :
SIN_{n-1}^{<\alpha^{\prime\prime}}(\alpha^{\prime\prime\prime})
\right \} \right ) = \alpha^{\prime} ) \ .
\end{eqnarray*}
\vspace{0pt}

\noindent It is not hard to see that  \  $ L_{k} [l] \vDash
\exists \alpha \ \ \varphi_{1}(\alpha,l) $ \  and for every  \
$\alpha $

\begin{eqnarray*}
L_{k}[l] \vDash \varphi_{1}(\alpha,l) \longleftrightarrow
SIN_{n-1}(\alpha)\wedge \exists
\alpha^{\prime\prime} < \alpha (
SIN_{n-1}(\alpha^{\prime\prime}) \wedge
\\
\qquad\qquad\qquad\qquad\qquad\qquad\qquad\qquad
\wedge OT \left (\alpha^{\prime\prime} \cap SIN_{n-1} \right ) =
\alpha_{3} ).
\end{eqnarray*}
\vspace{0pt}

\noindent Hence \ $ dom  \left
(\mathbf{S}_{\varphi}\overline{\overline{\lceil}} \chi  \right )$
\ contains the cardinal \ $  > \alpha_{1}  $ \ and with the help
of lemmas~2.7, 3.10 one can see that \  $dom  \left (
\widetilde{\mathbf{S}}_{n}^{\sin}\overline{\overline{\lceil}} \chi
\right ) $ \ contains such cardinals also in contradiction with
the assumption.\hfill $\dashv$
\vspace{12pt} \\

\noindent The building of matrix functions relies on the following
enumeration (in~~$L_{k}$) \  of subinaccessible cardinals:

\begin{definition}
\label{5.13.} \hfill {} \\
\hspace*{1em} Let \  $\alpha _{1}\leq k$.
\\
By the recursion on \  $ \tau <\alpha _{1}$ \  we define the
function \ $\gamma _{f}^{<\alpha _{1}}= ( \gamma _{\tau
}^{<\alpha _{1}} ) _{\tau }$~:

\[
    \gamma _{0}^{<\alpha _{1}}=0 ~; \quad for \quad \tau>0 \qquad
    \qquad\qquad\qquad\qquad\qquad\qquad
\]

\vspace{-6pt}

\[
    \gamma _{\tau }^{<\alpha _{1}}=\min \{ \gamma <\alpha_{1}:
    SIN_{n-1}^{<\alpha _{1}} ( \gamma  ) \wedge \forall
    \tau^{\prime }<\tau \quad
    \gamma_{\tau ^{\prime }}^{<\alpha_{1}}<\gamma \} \ .
\]
\vspace{0pt}

\noindent The inverse function \ $\tau _{f}^{<\alpha _{1}}= (
\tau_{\gamma }^{<\alpha _{1}} ) _{\gamma }$  is defined :

\[
    \tau =\tau _{\gamma }^{<\alpha _{1}}\longleftrightarrow
    \gamma =\gamma _{\tau }^{<\alpha _{1}}.
\]
\end{definition}

The proof of main theorem consists in creation in \ $L_{k}$ of \
the special matrix function possessing inconsistent properties;
this function arises out by a sequential complication of  its
following simplest form:

\begin{definition}
\label{5.14.} \hfill {} \\
\hspace*{1em} We \medskip name as a matrix function of the level \
$n$ \ below $ \alpha _{1}$ reduced to  $\chi $ the following
function \ \ \ $S_{\chi f}^{<\alpha _{1}}= ( S_{\chi \tau
}^{<\alpha _{1}} ) _{\tau }$ \ \ taking values : \vspace{6pt}
\[
    S_{\chi \tau }^{<\alpha _{1}}=\min_{\underline{\lessdot }}\{
    S:\exists \alpha < \alpha_{1}  \left ( \gamma_{\tau}^{<
    \alpha_{1}}< \alpha \wedge \sigma ^{\vartriangleleft
    \alpha _{1}} ( \chi ,\alpha ,S )\right ) \}  .
\]
\end{definition}
\vspace{0pt}

\noindent So, these values are matrices \ $S$ \ reduced to \
$\chi$ \ and singular on these carriers \ $\alpha$.
\\

\noindent As usual, if \  $\alpha _{1}<k$, \ then all functions
introduced are named restricted or relativized to \ $\alpha _{1}$;
\ if  \  $\alpha_{1} = k $, \ then all mentionings about \
$\alpha_{1} $ \ are dropped.
\\
Let us remind that all bounding ordinals  \  $\alpha_{1} $ \ are
assumed to be  \  $SIN_{n-2}$-cardinals or  \ $\alpha_{1}=k $.

\begin{lemma}
\label{5.15.} {\em (About matrix function absoluteness).} \\
\hspace*{1em} Let \ $\chi < \gamma_{\tau+1}^{<\alpha_1} < \alpha_2
< \alpha_1 \leq k$ \ and \ $\alpha_2 \in SIN_{n-2}^{<\alpha_1}$, \
then:
\\

\noindent {\em 1)} functions \ \ $\gamma _{f}^{<\alpha_2}$, \
$\gamma _{f}^{<\alpha_1}$ \ \ coincide on the set \ $\{
\tau^\prime: \gamma _{\tau^\prime} ^{<\alpha_2} \leq
\gamma_{\tau+1} ^{<\alpha_1} \}$:

\[
    \gamma _{\tau^\prime} ^{<\alpha_2} \leq
    \gamma_{\tau+1} ^{<\alpha_1}  \longrightarrow
    \gamma _{\tau^\prime} ^{<\alpha_2} =
    \gamma _{\tau^\prime} ^{<\alpha_1};
\]
\vspace{0pt}

\noindent {\em 2)} functions \ \ $S_{\chi f}^{<\alpha_2}$, \
$S_{\chi f}^{<\alpha _{1}} $ \ \ coincide on the set \ $\{
\tau^\prime: \chi \leq \gamma _{\tau^\prime} ^{<\alpha_2} \leq
\gamma_\tau ^{<\alpha_1} \}$:

\[
    \chi \leq
    \gamma _{\tau^\prime} ^{<\alpha_2} \leq
    \gamma_\tau ^{<\alpha_1}  \longrightarrow
    S _{\chi \tau^\prime} ^{<\alpha_2} =
    S _{\chi \tau^\prime} ^{<\alpha_1}.
\]
\end{lemma}
\textit{Proof} is coming from lemma~\ref{3.8.} 1) (where \  $n$ \
is replaced with \ $n-1$) \ and lemma 5.17~2 $(ii)$ below. \hfill $\dashv$ \\

From lemma~\ref{3.3.} where \ $n$ \ is replaced with \ $n-1$ \ it
comes
\begin{lemma}
\label{5.16.} \hfill {} \\
\hspace*{1em} For \ $\alpha_{1}<k $ \ functions

\[
    \gamma =\gamma _{\tau }^{<\alpha _{1}}, \quad S=S_{\chi
    \tau}^{<\alpha _{1}}
\]
\vspace{0pt}

\noindent are  \  $\Delta _{1}$-definable through \ $\chi $,
$\alpha _{1}$.
\\
For \  $\alpha_{1}=k $ \  these functions are  \ $\Pi_{n-1}
$-definable,  \ $\Delta_{n}$-definable respectively.
\end{lemma}

The following lemma contains the ``embryo'' of all further
reasoning: it establishes that matrix function has the property of
$\underline{\lessdot }$-monotonicity -- and further we shall
modify this function in order to preserve and to exclude this
property simultaneously.
\\
Therefore we shall often return to the idea of this lemma and of
its proof in various forms:

\begin{lemma}
\label{5.17.}
\quad \medskip \\
{\em 1)}\quad The function \  $S_{\chi f}^{<\alpha _{1}}$ \ is
$\underline{\lessdot }$-monotone, that is for every
\\
$\tau _{1}$, $ \tau _{2}\in dom \left( S_{\chi f}^{<\alpha _{1}}
\right) $
\[
    \quad \tau _{1}<\tau _{2}\longrightarrow S_{\chi
    \tau _{1}}^{<\alpha _{1}} \underline{\lessdot }
    S_{\chi \tau _{2}}^{<\alpha _{1}} \ .
\]
{\em 2)}\quad Let \ $ \tau \in dom \left ( S_{\chi f}^{<
\alpha_{1}} \right ) $, \ then:
\\
\begin{itemize}
\item[(i)] \ $ \{ \tau^{\prime}: \chi \leq
\gamma_{\tau^{\prime}}^{< \alpha_{1}} \leq \gamma_{\tau}^{<
\alpha_{1}} \} \subseteq dom \left ( S_{\chi f}^{< \alpha_{1}}
\right ) $~;
\\
\item[(ii)] \ if \ $\gamma_{\tau+1}^{<\alpha_{1}} $ \ and the
matrix \ $S_{\chi \tau}^{< \alpha_{1}}$ \ exist, \ then this
matrix possesses a carrier \ $\alpha \in \;
]\gamma_{\tau}^{<\alpha_{1}}, \gamma_{\tau+1}^{<\alpha_{1}}[ $~.
\end{itemize}
\end{lemma}
\textit{Proof.} \ Statements 1, 2 $(i)$ follow from
definition~5.14 immediately. The proof of 2~$(ii) $ represents
the typical application of lemma~3.2 about restriction. Since
there is the matrix  \  $ S=S_{\chi \tau}^{< \alpha_{1}} $ \  on
some carrier
\[
    \alpha \in ]\gamma_{\tau}^{<\alpha_{1}}, \alpha_{1} [,
\]
the following \  $\Sigma_{n-1} $-proposition  \  $\varphi
(\chi,\gamma_{\tau}^{<\alpha_{1}}, S ) $~:
\[
    \exists \alpha \ \ (\gamma_{\tau}^{< \alpha_{1}}< \alpha
    \wedge \sigma (\chi, \alpha, S))
\]
is fulfilled below  \  $\alpha_{1} $, \ that is fulfilled  \
$\varphi^{\triangleleft
\alpha_{1}}(\chi,\gamma_{\tau}^{<\alpha_{1}}, S) $.

This proposition contains individual constants
\[
    \chi < \gamma_{\tau +1}^{<\alpha _{1}}, \quad
    \gamma _{\tau }^{<\alpha_{1}}<
    \gamma _{\tau +1}^{<\alpha _{1}},\quad S\vartriangleleft
    \gamma _{\tau +1}^{<\alpha _{1}}
\]
due to lemma~\ref{5.3.}~2), and that is why by lemma~\ref{3.2.}
(where \ $n$ \ replaced with \ $n-1$) the \ $ SIN_{n-1}^{<\alpha
_{1}}$-cardinal \ $\gamma _{\tau +1}^{<\alpha_{1}} $ \  restricts
the proposition \ $\varphi $, \ that is fulfilled the formula
\[
    \exists \alpha \in \left[ \gamma _{\tau }^{<\alpha _{1}},
    \gamma_{\tau +1}^{<\alpha _{1}}\right[ ~\sigma ^{\vartriangleleft
    \alpha _{1}} ( \chi ,\alpha ,S ).
\]
\hfill $\dashv$
\\
\quad \\
From here, lemma 5.12 (for \ $\chi = \chi^{\ast}$) and
lemma~\ref{5.3.}~2) \ (for \ $\alpha_1 = k$) \ it comes directly

\begin{lemma}
\label{5.18.} \quad \\
{\em 1)} The unrestricted function \ $S_{\chi^{\ast} f}$ \ is
defined on the final segment of $k$:

\[
    dom \left( S_{\chi ^{\ast } f } \right)=
    \left \{ \tau : \chi^{\ast} \leq \gamma_{\tau} < k \right \};
\]
\vspace{0pt}

\noindent {\em 2)}\quad this function stabilizes, that is there is
an ordinal \ $\tau ^{\ast }>\chi ^{\ast }$ \ such that for every \
$\tau \geq \tau ^{\ast }$ there exists
\[
    S_{\chi ^{\ast }\tau }=S_{\chi ^{\ast }\tau ^{\ast }};
\]
therefore for every \ $\gamma_{\tau} \geq \chi ^{\ast }$ \
\[
    S_{\chi ^{\ast }\tau } \ \underline{\lessdot}
    \ S_{\chi ^{\ast }\tau ^{\ast }}.
\]
\end{lemma}

Now let us make an outline of the first approach to the idea of the main
theorem proof.
\\
We shall try to obtain the required contradiction on the following
way. The lower index \  $\chi ^{\ast } $ \ will be dropped for
some brevity.
\\
Let us consider the matrix function \ $S_f$ in its state of
stabilizing, that is consider \ an arbitrary sufficiently large \
$\tau _{0}>\tau ^{\ast }$, \  the matrix \  $ S_{\tau _{0}}$ \ on
some carrier \ $\alpha _{0}\in \left] \gamma_{\tau _{0}},\gamma
_{\tau _{0}+1}\right[ $, \  the prejump cardinal \
$\alpha^0=\alpha _{0}^{\Downarrow }$ \  and the function \  $
S_{f}^{<\alpha^0}$.
\\
The reasoning forthcoming will be more convenient and transparent,
if the following notion will be introduced:
\\
\quad \\
\textit{We shall say, that one, standing on an ordinal \ $\alpha
> \chi^{\ast}$, can detect (or can see) some information or some
object below \ $\alpha$, \ iff this information or object can be
defined by some formula relativized to this \ $\alpha$, and
defined just as some ordinal \ $< \chi^{\ast +}$.}
\\
\quad \\
So, the idea of the matrix function research in its basic form
will be mainly the following:
\\
\hfill {} \\
\textit{Assume, that some matrix, which is the value of such
function, is considered on its carrier\  $ \alpha _{0}$ .\
Standing on the prejump cardinal  \ $\alpha^0=\alpha
_{0}^{\Downarrow }$ \ after \ $\chi^{\ast}$ \ of this matrix on
this carrier, one should observe the situation below \ $\alpha^0$
and, mainly, the behavior of this very function, but in its
relativized to \ $\alpha^0$ \ form, with a view to detect its
properties}.
\\
\hfill {} \\
The information about this function, that one can detect and
express as some ordinals \ $< \chi^{\ast +}$ -- this information
can provide certain contradictions.
\\
For instance, it can provide the increasing of the order type of
\ $S_{\tau_0}$ \ due to lemma~4.6, and, hence, can destroy the
stabilizing of the whole function \ $S_f$ \ in contraction with
lemma 5.18.
\\
More explicitly in this context:
\\
By lemma~\ref{5.15.}~2) the function \ $S_f^{< \alpha^0}$
coincides with \  $S_{f}$ \ on the ordinal \  $\tau _{0}$ \  and
also \ $\underline{\lessdot }$-monotone. Let us apply
lemma~\ref{4.6.} to this situation, considering
\[
    \overline{\delta }= \mathrm{sup}_{\tau} Od
    (S_{\tau }^{<\alpha^0} ),  \quad \overline{\chi }=
    \chi^{\ast }, \quad \overline{\alpha}_{0}=
    \alpha _{0}^{\downarrow }, \quad \overline{\alpha}_1 =
    \alpha_{0}.
\]
Suppose \ $\overline{\delta }<\chi ^{\ast +}$, \ then one can
detect \ $\overline{\delta}$ \ standing on \ $\alpha^0$. \ But now
from definition~\ref{5.7.} and lemma~\ref{5.5.} it comes that all
conditions of lemma~\ref{4.6.} are fulfilled and so it implies the
contradiction:
\[
    Od ( S_{\tau ^{\ast }} ) \leq \overline{\delta }<OT
    ( S_{\tau_{0}} ) \leq Od ( S_{\tau ^{\ast }} )  .
\]
Hence, in reality \ $\overline{\delta }=\chi ^{\ast +}$ \ and the
function \ $S_{f}^{<\alpha^0}$ \ is \ $\underline{\lessdot
}$-nondecreasing up to \  $\chi ^{\ast +}$.

It means that the function \  $ ( S_{\tau } ) _{\tau <\tau _{0}}$,
\ being \ $\underline{\vartriangleleft }$~-bounded by the ordinal
\  $Od \left ( S_{\tau^{\ast}} \right ) < \chi^{\ast +} $, \ loses
this boundedness after its extension on the set \  $\{ \tau
:\gamma _{\tau}^{<\alpha^0 }<\alpha^0 \} $.
\\
One can see that it happens because of losing the properties of
subinaccessibility of all levels \  $\geq n$ \  by all cardinals \
$\leq \gamma _{\tau _{0}}$ \  after their relativization to \
$\alpha^0 $; \  at the same time there appears some \
$SIN_{n-1}$-cardinals (relativized to \ $\alpha^0$) that are not
such in the universe (Kiselev~\cite{Kiselev3}, \cite{Kiselev4}).
\\

All these conclusions mean that many important properties of lower
levels of the universe do not extend up to the relativizing
cardinals, namely, to jump cardinals of reduced matrices, which
are values of matrix functions.

\textit{In order to prevent this phenomenon we shall introduce
special cardinals named disseminators that are extending such
properties without distortion}.
\\
Extensions of this kind sometimes are produced by subinaccessible
cardinals but for what follows it is needed to widen this notion
significantly.

\newpage
\quad 

\newpage

\section{Disseminators}
\setcounter{equation}{0}

Let us remind that all the reasoning is conducted and formulas are
interpreted in  \  $ L_{k} $ \  (or in \ $\mathfrak{M}$ \ if, of
course, the context does not mean some another case). Also all \
$\vartriangleleft$-bounding ordinals belong to \ $SIN_{n-2}$.
\\
We shall introduce the notion of disseminator only for the
constructive universe for more clearness of reasoning although it
can be done without any loss of generality. Further we shall
consider classes \  $\Sigma _{m}$, \  $\Pi _{m}$ \  of the fixed
arbitrary level \  $m>3$ (unless otherwise specified).
\begin{definition}
\label{6.1.} \hfill {} \\
\hspace*{1em} Let
\[
    0<\alpha <\alpha _{1} \leq k,\quad X\subseteq \alpha, \quad
    X\neq \varnothing .
\]
The ordinal \  $\alpha $ \  is named the disseminator of the level
\  $m$ \  with the data base \  $X$ \ below \  $\alpha _{1}$ \ iff
for every \ $\mathfrak{n} \in \omega_0$ \ and train \
$\overrightarrow{a}$ \  of  ordinals \  $\in X$

\vspace{6pt}
\[
    U_{m}^{ \Sigma \vartriangleleft \alpha _{1}} ( \mathfrak{n},
    \overrightarrow{a} ) \longrightarrow  U_{m}^{\Sigma
    \vartriangleleft \alpha } ( \mathfrak{n}, \overrightarrow{a} ) \quad.
\]
\vspace{0pt}

\noindent The formula defining in \ $L_{k}$ \ the set of all
disseminators with the data base \ $X$ \ of the level \ $m $ \
below \ $\alpha _{1}$ \ is denoted by \ $ SIN_{m}^{<\alpha
_{1}}\left[ X\right] ( \alpha  ) $; \ this set itself is also
denoted by \ $SIN_{m}^{<\alpha _{1}}\left[ X \right] $, \ while
its disseminators -- by the common symbol \  $\delta ^{X}$ \  or,
briefly, by \ $\delta $, \ pointing to \ $\alpha_1$ \ in the
context.
\end{definition}
Let us remind that the symbol \ $U_{m}^{\Sigma}$ \  here denotes
the \  $\Sigma _{m}$-formula which is universal for the class of
\  $\Sigma _{m}$-formulas but without any occurrences of the
constant \ $\underline{l}$ \ (see remark after lemma 2.5).
\\
It is possible to obtain the definition of the disseminator
notion in more wide sense if occurrences of the constant \ $
\underline{l}$ \  are allowed:
\\
\textit{for every train \ $\overrightarrow{a}$ \ of ordinals \
$\in X$}
\[
    \left\| u_{m}^{\Sigma \vartriangleleft \alpha _{1}}
    ( \overrightarrow{a},\underline{l} ) \right\| \leq
    \left\| u_{m}^{\Sigma \vartriangleleft \alpha }
    (\overrightarrow{a},\underline{l} ) \right\| ,
\]
and all the following reasoning can be carried over in this case.
\\
\hfill {} \\
The following two lemmas are quite analogous to lemmas~3.3, 3.2
and can be proved in the same way:

\begin{lemma}
\label{6.2.} \hfill {} \\
\hspace*{1em} The formula \ $SIN_{m}^{<\alpha _{1}}\left[
X\right] ( \alpha  ) $ \ belongs to the class \ $\Pi _{m}$ \  for \  $%
\alpha _{1}=k$ \  and to the class $ \Delta _{1}$ \ for \ $\alpha
_{1}<k$.
\end{lemma}

The following obvious lemma justifies the  term ``disseminator''
since it shows that such an ordinal really extend \
$\Pi_{m}$-properties (containing constants from its base) from
lower levels of the universe up to relativizing cardinals:

\begin{lemma}
\label{6.3.} {\em (About extending)} \hfill {} \\
\hspace*{1em} Let
\[
    X\subseteq \alpha, \quad \alpha <\alpha_{1}, \quad
    \alpha \in SIN_{m}^{<\alpha_{1}}\left[ X\right]
\]
and a proposition \ $\forall x~\varphi ( x,\overrightarrow{a} ) $
\ has a train \ $\overrightarrow{a}$ \ of constants \ $\in X$, \
$\varphi \in \Sigma _{m-1}$, \ then :

\vspace{6pt}
\[
    \forall x\vartriangleleft \alpha ~\varphi ^{\vartriangleleft \alpha }
    (x,\overrightarrow{a} ) \longrightarrow \forall x\vartriangleleft
    \alpha_{1}~\varphi ^{\vartriangleleft \alpha _{1}}
    ( x,\overrightarrow{a} )~.
\]
\vspace{0pt}

\noindent In this case we shall say, as above, that below \
$\alpha _{1} $ \ the ordinal$ \  \alpha $ \  extends or prolongs \
$\forall x~\varphi $ \  up to \  $\alpha _{1}$.
\\
Considering the same in the inverted form for \  $\varphi \in \Pi
_{m-1}$:

\vspace{6pt}
\[
    \exists x\vartriangleleft \alpha _{1}~\varphi ^{\vartriangleleft
    \alpha _{1}} ( x,\overrightarrow{a} ) \longrightarrow \exists
    x\vartriangleleft \alpha ~\varphi ^{\vartriangleleft \alpha }
    ( x,\overrightarrow{a} )~,
\]
\vspace{0pt}

\noindent we shall say that below \  $\alpha _{1}$ \  the ordinal
\ $\alpha $ \  restricts or relativizes the proposition \ $\exists
x~\varphi $. \hspace*{\fill} $\dashv$
\end{lemma}
It is not hard to see that the class of disseminators \ $\alpha
\in SIN_{m}^{<\alpha _{1}} \left[ X\right]$ \  т \textit ({in wide
sense}) below \ $\alpha _{1}$ \ of level \ $m$ \ and with base \
$X=\alpha $ \ coincides with \ $ SIN_{m}^{<\alpha _{1}}$.
\\
 The properties of disseminators become more transparent if, in
addition, they possess the  properties of subinaccessibility of
lower levels:
\begin{lemma}
\label{6.4.} \hfill {} \\
\hspace*{1em} Let
\[
    X\subseteq \alpha, \quad \alpha <\alpha _{1}, \quad
    \alpha \in SIN_{m-1}^{<\alpha _{1}},
\]
then the following statements are equivalent:
\\
\\
{\em 1)}\quad $\alpha \in SIN_{m}^{<\alpha _{1}}\left[ X\right]$~;
\\
\\
{\em 2)}\quad for every \ $\mathfrak{n} \in \omega_0$ \ and every
train \ $\overrightarrow{a}$ \ of constants \ $\in X $
\[
    L_{k} \vDash \exists x\vartriangleleft
    \alpha _{1}~U_{m-1}^{\Pi \vartriangleleft \alpha _{1}}
    (\mathfrak{n}, x,\overrightarrow{a} ) \longleftrightarrow \exists
    x\vartriangleleft \alpha ~U_{m-1}^{\Pi \vartriangleleft
    \alpha _{1}} (\mathfrak{n}, x,\overrightarrow{a} )~;
\]

\noindent {\em 3)}\quad for every \  $\Pi_{m-1}$-formula \  $
\varphi ( x,\overrightarrow{a} ) $ \  and every train \
$\overrightarrow{a} $ \ of  constants~~$\in X$
\[
    L_{k} \vDash \exists x\vartriangleleft \alpha _{1}~
    \varphi^{\vartriangleleft \alpha _{1}} ( x,\overrightarrow{a} )
    \longleftrightarrow \exists x\vartriangleleft \alpha ~
    \varphi^{\vartriangleleft \alpha _{1}} ( x,\overrightarrow{a} )~.
\]
\hspace*{\fill} $\dashv$
\end{lemma}

This lemma advances the following important

\begin{definition}
\label{6.5.} \hfill {} \\
{\em 1)}\quad The minimal disseminator of the class
\[
    SIN_{m}^{<\alpha _{1}}\left[ X\right] \cap
    SIN_{m-1}^{<\alpha_{1}}
\]
is named the generating disseminator with the data base \ $X$ \
below \ $\alpha_1$ \ and is denoted by the general symbol \
$\check{\delta}^{X}$, \ or \ $\check{\delta}$ \ or \ $\delta$;

\noindent {\em 2)}\quad without this condition of minimality,
disseminators of this class are named floating disseminators below
\ $\alpha_1$ \ and are denoted by the general symbol \
$\widetilde{\delta }^{X}$ \ or, briefly, by \ $\widetilde{\delta}$
\ or \ $\delta$, \ pointing to \ $\alpha_1$ \ and to their other
attributes in the context. \hspace*{\fill} $\dashv$
\end{definition}

As usual, indices \ $m$, \ $\alpha _{1}$, $X$ \ are omitted if
they are arbitrary or pointed out by the context.

\noindent These terms ``generating disseminator'' and ``floating
disseminator'' are justified by the following circumstance: the
\textit{generating disseminator} is \textit{uniquely} defined
through its base \ $X$ \ below \ $\alpha_1$ \ as minimal, while
the \textit{floating disseminator} may be not, and its value is
not specified (it is floating); moreover, the ``generating
disseminator'' term is justified by the following

\begin{lemma}
\label{6.6.} \hfill {} \\
\hspace*{1em} Let
\[
    X\subseteq \alpha _{0}, \quad  \alpha _{0}<\check{\delta}<\alpha _{1},
\]
\ $\check{\delta}$ \  \smallskip be a generating disseminator \
$\in SIN_{m}^{<\alpha _{1}}\left[ X\right] $ \ and \  $\varphi  (
\alpha ,\overrightarrow{a} ) $ \ be a \ $\Sigma _{m-1}$-formula
\smallskip with the ordinal variable \ $\alpha $ \ and the train
\  $\overrightarrow{a}$ \ of \smallskip constants~~$\in X$.
\\
Suppose
\[
    \forall \alpha \in \left] \alpha _{0} ,\check{\delta}
    \right[~\varphi ^{\vartriangleleft \check{\delta}}
    ( \alpha ,\overrightarrow{a} ),
\]

\noindent then there exists some \  $\alpha _{0}^{\prime }\in
\left[ \alpha _{0}^{{}},\check{\delta}\right[ $ \ such that

\[
    \forall \alpha \in \left] \alpha _{0}^{\prime },\alpha _{1}
    \right[ ~~\varphi^{\vartriangleleft \alpha _{1}}
    ( \alpha ,\overrightarrow{a} ).
\]
\end{lemma}
\textit{Proof.} \ One should point out here that the variable \
$\alpha$ \ under the quantor \ $\forall$ \ runs not through all \
$\check{\delta}$ ,\ as it was in lemma~6.3, but only through its
some final segment \ $\left[ \alpha
_{0}^{{}},\check{\delta}\right[ $ and nevertheless this \
$\check{\delta}$ \ \emph{generates} the extension of \ $\varphi (
\alpha ,\overrightarrow{a} )$\ up to \ $\alpha_{1}$.\
\\
Beforehand the following remark should be done: for arbitrary
ordinals \ $\alpha \in SIN_{e}$, \ $e>1$, \ $\beta $: \qquad
\qquad $ \beta < \alpha \longleftrightarrow \beta \vartriangleleft
\alpha$.
\\
It is sufficient to consider first the case when  \  $\alpha_{0}
\in SIN_{m-1}^{<\alpha_{1}} $. \ From the minimality of  \
$\check{\delta}$ \ it follows that \ $\alpha_{0} \notin
SIN_{m}^{<\alpha_{1}}[X] $ \ and by lemma 6.4 it means that there
is  \ $\Sigma_{m-1}$-formula \ $\varphi_{1}(\alpha,
\overrightarrow{a_{1}}) $ \  with the train \
$\overrightarrow{a_{1}} $ \  of constants $\in X $ \ such that

\begin{equation}
\label{e6.1} \forall \alpha < \alpha_{0} \ \
\varphi_{1}^{\triangleleft \alpha_{0}}(\alpha,
\overrightarrow{a_{1}})~,
\end{equation}
\begin{equation}
\label{e6.2} \exists \alpha^{\prime} < \alpha_{1} \neg
\varphi_{1}^{\triangleleft \alpha_{1}}(\alpha^{\prime},
\overrightarrow{a_{1}})~.
\end{equation}
\vspace{0pt}

\noindent Since  \  $\alpha_{0}, \check{\delta}\in
SIN_{m-1}^{<\alpha_{1}} $, \ statement (6.1) is equivalent to

\vspace{6pt}
\begin{equation}
\label{e6.3} \forall \alpha < \alpha_{0} \ \
\varphi_{1}^{\triangleleft  \check{\delta}}(\alpha,
\overrightarrow{a_{1}})~.
\end{equation}
\vspace{0pt}

\noindent The cardinal  \ $\check{\delta}$ \  is the disseminator
and that is why (6.2) implies

\vspace{6pt}
\[
    \exists \alpha^{\prime} < \check{\delta} \neg
    \varphi_{1}^{\triangleleft  \check{\delta}}(\alpha^{\prime},
    \overrightarrow{a_{1}})~.
\]
\vspace{0pt}

\noindent From here and (6.3) it comes that there is the ordinal
 \  $\alpha_{0}^{\prime} \in [ \alpha_{0}, \delta [ $ \  such
that

\begin{equation}
\label{e6.4} \forall \alpha < \alpha_{0}^{\prime} \ \
\varphi_{1}^{\triangleleft \check{\delta}}(\alpha,
\overrightarrow{a}_{1})\wedge \neg \varphi_{1}^{\triangleleft
 \check{\delta}}(\alpha_{0}^{\prime}, \overrightarrow{a}_{1})~.
\end{equation}
\vspace{0pt}

\noindent Now from the condition of the lemma it comes the
proposition

\vspace{6pt}
\[
    \forall \alpha^{\prime}, \alpha < \check{\delta} \left ( \neg
    \varphi_{1}^{\triangleleft \check{\delta}}(\alpha^{\prime},
    \overrightarrow{a}_{1}) \wedge \alpha > \alpha^{\prime}
    \longrightarrow  \varphi^{\triangleleft \check{\delta}}(\alpha,
    \overrightarrow{a}) \right )~.
\]
\vspace{0pt}

\noindent It contains constants $\in X $ \  so by lemma~6.4~3) the
disseminator  \ $\check{\delta}$ \ extends it up to  \  $
\alpha_{1} $:

\begin{equation}
\label{e6.5} \forall \alpha^{\prime}, \alpha < \alpha_{1}   \left
( \neg \varphi_{1}^{\triangleleft \alpha_{1}}(\alpha^{\prime},
\overrightarrow{a}_{1}) \wedge \alpha > \alpha^{\prime}
\longrightarrow  \varphi^{\triangleleft \alpha_{1}}(\alpha,
\overrightarrow{a}) \right )~.
\end{equation}
\vspace{0pt}

\noindent Since   \ $\check{\delta}\in SIN_{m-1}^{<\alpha_{1}} $ \
from (6.4) it comes  \  $\neg \varphi_{1}^{\triangleleft
\alpha_{1}}(\alpha_{0}^{\prime}, \overrightarrow{a}_{1}) $ \ and
\smallskip therefore (6.5) implies  \  $\forall \alpha \in [
\alpha_{0}^{\prime}, \alpha_{1}[ ~~ \varphi^{\triangleleft
\alpha_{1}}(\alpha, \overrightarrow{a}) $.\hfill $\dashv$
\\

Using this lemma we shall say, as above, that \ $\check{\delta}$ \
\textit{extends or prolongs} the proposition  \ $\forall \alpha \
\ \varphi $ \  up to  \  $\alpha_{1} $, \ or, in the inverted
form, that \ $\check{\delta}$ \  \textit{restricts or relativizes}
\ $\exists \alpha \neg \varphi $ \  below  \ $\alpha_{1} $, \
pointing out ordinals  \  $\alpha_{0}, \alpha_{0}^{\prime} $ \ by
the context (if it will be necessary).

\begin{corollary}
\label{6.7.} \hfill {} \\
\hspace*{1em} Let  \ $\check{\delta}\in SIN_{m}^{<\alpha
_{1}}\left[ X\right] $ \  be a generating disseminator below  \
$\alpha_{1} $ ,\ then
\[
    \sup  ( \check{\delta}\cap SIN_{m-1}^{<\alpha _{1}} ) = \check{\delta}.
\]
\end{corollary}
\textit{Proof.} \ Let, on the contrary,
\[
    \alpha_{0} = sup  \left ( \check{\delta}\cap SIN_{m-1}^{<
    \alpha_{1}} \right ) < \check{\delta},
\]
then \ $\alpha_{0}\in SIN_{m-1}^{<\alpha_{1}} $ \ and the argument
from the previous proof can be applied. The disseminator \
$\check{\delta}$ \ restricts the proposition
\[
    \exists \alpha > \alpha_{0}  \ \ SIN_{m-1}(\alpha)
\]
that is fulfilled below \ $\alpha_{1} $ \  since \
$\check{\delta}\in SIN_{m-1}^{<\alpha_{1}} $. \ So there exists
\[
    \alpha \in \; ] \alpha_{0}, \check{\delta} [ \; \cap SIN_{m-1}^{< \alpha_{1}}
\]
contrary to the assumption.\hfill $\dashv$
\\
Now it is clear that in the proof of lemma 6.6 the condition
$\alpha_{0}\in SIN_{m-1}^{<\alpha_{1}} $ \  should be dropped.
\\
\quad \\
In what follows we shall investigate the mutual disposition of
disseminators of one and the same level \ $m$ \ below one and the
same \ $\alpha _{1}$ \ but with different data bases. How these
bases influence  their location? The discussion of such questions
is more clear when\ the  data base \ $X$ \ is some ordinal \
$\rho$ \ (that is the set of all smaller ordinals). In order to
point out this case we shall write
\[
    SIN_{m}^{<\alpha _{1}}\left[ <\rho \right]
    \mathrm{\it \quad instead\ of\quad} SIN_{m}^{<\alpha _{1}}
    \left[ X\right].
\]
Further all bases \ $\rho $ \  are \textit{limit} ordinals (if
some another case is not meant).
\\
\hfill {} \\
The following lemma, although obvious, represents, nevertheless,
several important technical tools of disseminator analysis:

\begin{lemma}
\label{6.8.}
\quad \\
{\em 1)}\quad Let
\[
    \delta < \alpha _{1}, \quad \delta \in SIN_{m}^{<\alpha _{1}}\left[ X\right]
    \cap SIN_{m-1}^{<\alpha_{1}}
\]
\vspace{0pt}

\noindent then every \ $SIN_{m-1}^{<\alpha_{1}}$-cardinal in \ $[
\delta, \alpha_1 [$ \ is again the disseminator with the same base
below \ $\alpha_1$:

\[
    \alpha \in SIN_{m-1}^{<\alpha_{1}} \wedge \alpha > \delta
    \longrightarrow \alpha \in SIN_{m}^{<\alpha_{1}}\left[
    X\right];
\]
\vspace{0pt}

\noindent {\em 2)}\quad the increasing of the base \ $\rho$ \
implies the nondecreasing of the generating disseminator below \
$\alpha_1$:

\[
    \rho _{1}<\rho _{2}\longrightarrow \check{\delta}^{\rho
    _{1}}\leq \check{\delta}^{\rho _{2}};\quad \rho _{1}^{+}\leq
    \rho_{2}\longrightarrow \check{\delta}^{\rho_{1}}<
    \check{\delta}^{\rho_{2}};
\]
\vspace{0pt}

\noindent {\em 3)}\quad the passage to the limit of bases
involves the passage to the limit of generating disseminators
below \ $\alpha_1$:
\[
    \lim_{i} \rho _{i}=\rho \longrightarrow \lim_i \check{\delta}^{\rho
    _i}=\check{\delta}^{\rho }.
\]
\end{lemma}
\textit{Proof.} 1). The upper index \ $\alpha_{1} $ \ will be
dropped for shortness. Let us turn to definition~6.1 (or to
lemma~6.4) and consider the proposition
\[
    \exists x \ \
    U_{m-1}^{\Pi}( \mathfrak{n},x,\overrightarrow{a})
\]
for an arbitrary \ $\mathfrak{n} \in \omega_0$ \ and a train \
$\overrightarrow{a} $ \ of constants~~$\in X $.
\\
Suppose, that
\[
    \exists x \vartriangleleft \alpha_1
    ~U_{m-1}^{\Pi}(\mathfrak{n},x,\overrightarrow{a}).
\]
Since \ $\delta $ \ is the disseminator it implies
\[
    \exists x \vartriangleleft \delta ~
    U_{m-1}^{\Pi}(\mathfrak{n},x,\overrightarrow{a})
\]
and, thus,
\[
    \exists x \vartriangleleft \alpha \ \
    U_{m-1}^{\Pi}(\mathfrak{n},x,\overrightarrow{a}) .
\]
Here  \ $U_{m-1}^{\Pi}$ \ should be  \ $\vartriangleleft$-bounded
by  \ $\alpha $ \ since  \ $\alpha \in  SIN_{m-1} $ \  and so \
$\alpha \in SIN_{m}[X] $. \ Statement 2) comes from definitions
and
along with 1) implies 3). \hfill $\dashv$\\
\quad \\
To obtain more detailed information about disposition of
disseminators it is natural and reasonable to consider matrix
carriers.
\\
Further by \  $\widehat{ \rho }$ \  is denoted the closure of \
$\rho $ \  under the pair function.

\begin{definition}
\label{6.9.} \hfill {} \\
\hspace*{1em} Let \ $\chi<\alpha<\alpha_1$ \ and \ $S$ \ be a
matrix which is reduced to a cardinal \ $\chi$ \ and is singular
on a carrier \ $\alpha $.
\\
{\em 1)}\quad We name as disseminator for \ $S$ \  on \ $\alpha $
\ (or as disseminator for this carrier) of the level \  $m$ \ with
a data base \  $X$ \  every disseminator \ $\delta $ \ below the
prejump cardinal \ $\alpha_{\chi}^{\Downarrow }$ \ with the same
parameters, that is every
\[
    \delta \in SIN_{m}^{<\alpha_{\chi}^{\Downarrow }}\left[ X\right]
    \cap SIN_{m-1}^{<\alpha_{\chi}^{\Downarrow }}.
\]
In this case we say also that the matrix \  $S$ \ on \ $\alpha$ \
(or its carrier \ $\alpha $) \ possesses this disseminator (with
these parameters).
\\
\quad \\
{\em 2)}\quad We call as the {\sl eigendisseminator} of the matrix
\ $S$ \ on \ $\alpha$ \ of the level \ $m$ \ every its
disseminator of this level with the data base \ $\rho =
\widehat{\rho_1}$, $\rho_1 = Od(S)$, \ and denote it through the
general symbol \ $\widetilde{\delta}^S$ \ or \ $\delta^S$, \ and
its data base \ $\rho$ --- through \ $\rho^S$; \ if \
$\widetilde{\delta}^S$ \ is minimal with this base \ $\rho^S$, \
then it is called the {\sl generating eigendisseminator} for \ $S$
\ on \ $\alpha$ \ and is denoted through \ $\check{\delta}^S$.
\\
\quad \\
{\em 3)}\quad The matrix \ $S$ \ is called the disseminator matrix
or, more briefly, the \ $\delta $-matrix of the level \ $m$ \
admissible on the carrier \ $\alpha$ for \ $\gamma =
\gamma_{\tau}^{<\alpha_1}$ \ below \ $\alpha_1$ \ iff it possesses
some disseminator \ $ \delta < \gamma$ \ of the level \ $m$ \ with
some base \ $\rho \le \chi^{\ast +}$ \ on this carrier such that \
$S \vartriangleleft \rho, \ \rho = \widehat{\rho }$; \ such \
$\rho$ is also called admissible on  \ $\alpha$ for \ $\gamma =
\gamma_{\tau}^{<\alpha_1}$ \ (also below \ $\alpha_1$).
\label{c6}
\endnote{
\ p. \pageref{c6}. \ Here one should point out the following
interesting notion of base nonexcessiveness: \textit{the base \
$\rho$ \ of the matrix S admissible on the carrier \ $\alpha$}\ is
named nonexcessive on this carrier, iff the decreasing of \ $\rho$
\ implies the decreasing of generating disseminator for \ $S$ \ on
the same carrier \ $\alpha$:
\[
    \forall \rho^{\prime} \left ( \rho^{\prime} < \rho
    \longrightarrow \check{\delta}^{\rho^{\prime}} <
    \check{\delta}^{\rho} \right ).
\]

\noindent The point is that among bases of the matrix \ $S$ \
disseminators on the carrier \ $\alpha$ \ there are possible bases
\ $\rho$ \ which possess certain ``excessive'' information in the
sense that for some \ $\rho^{\prime} < \rho$ \ still
\[
    \check{\delta}^{\rho^{\prime}} = \check{\delta}^{\rho}
\]
on \ $\alpha$ \ and, so, instead of \ $\rho$ \ it \textit{can be
used the smaller base \ $\rho^{\prime} < \rho $ \ without any
loss} for the location of the generating disseminator (and, hence,
for its action). It is natural to use bases free of such
excessiveness, that are nonexcessive bases. One should note here
that the using of some smaller base, may be, is possible,
\textit{but only} when the corresponding disseminator decreases.

\noindent It is not hard to prove the lemma which establishes the
remarkable property: not only such base defines the location of
generating disseminator, but, inversely, such disseminator defines
its nonexcessive base even for different matrices on different
carriers:
\\

\noindent {\bf Lemma}
\\
\hspace *{1em} \em \noindent Let \ $S_1$, \ $S_2$ \ be matrices on
carriers \ $\alpha_1$, \ $\alpha_2$ \ which possess generating
disseminators \ $\check{\delta}_1$, \ $\check{\delta}_2$ \ of the
level \ $m$ \ with nonexcessive bases \ $\rho_1$, \ $\rho_2$ \ on
these carriers respectively, then:
\[
    if \quad \check{\delta}_1 = \check{\delta}_2 \in
    SIN_{m-1}^{<\alpha_1^{\Downarrow}} \cap
    SIN_{m-1}^{<\alpha_2^{\Downarrow}}
    \quad then \quad \rho_1 = \rho_2.
\]
\em 
\\
} 
\\

\noindent {\em 4)}\quad In what follows by \ $Lj^{<\alpha } ( \chi
) $ \ is denoted the \ $\Delta _{1}$-formula:
\[
    \chi <\alpha \wedge SIN_{n-1}^{<\alpha } ( \chi  ) \wedge
    \sum rng \left ( \widetilde{\mathbf{S}}_{n}^{\sin
    \vartriangleleft \chi } \right ) \in B_{\chi }\wedge
    \sup dom \left (\widetilde{\mathbf{S}}_{n}^{\sin
    \vartriangleleft \chi } \right )=\chi ~
\]
and the cardinal \ $\chi$ \ with this property is called saturated
below \ $\alpha$.
\\
Now let us accept the following important agreement:
\\
everywhere further we shall consider singular matrices \ $S$ \ on
carriers \ \mbox{$\alpha \in SIN_{n-2}$} \ reduced to some
saturated cardinal \ $\chi$ \ which has the property
\\
\hspace*{9em} $\chi^\ast \le \chi \wedge Lj^{<\alpha } ( \chi )$.
\hspace*{\fill} $\dashv$
\end{definition}

It should be noted that there exist cardinals of this kind, for
instance, the cardinal \  $\chi = \chi ^{\ast }$ \ for every \
$\alpha >\chi ^{\ast }$, $\alpha \in SIN_{n-2}$, \ as one can see
it from lemmas~\ref{5.5.}, \ref{3.8.} (for \ $n-1$ \ as \ $n$).
\\
This agreement makes it possible to apply lemma~\ref{4.6.} about
the spectrum type to these  cardinals.
\\
Let us \textit{fix this cardinal} \ $ \chi$ \ up to a special
remark, so the lower index  \ $\chi $ \  will be dropped as a
rule; for example, the prejump cardinal  \
$\alpha_{\chi}^{\Downarrow} $ \ will be denoted simply by \
$\alpha^{\Downarrow} $ \  and so on.
\\

\noindent Next lemma analyzes the set of matrix disseminators with
bases exceeding this matrix order type itself.

\begin{lemma}
\label{6.10.} \hfill {} \\
\hspace*{1em} Let a matrix \ $S$ \ on a carrier \ $\alpha$ \ be
reduced to the cardinal \ $\chi$, \ and let it possess some
disseminator of the level \ $m$ \ with a base \ $\in \left[ OT ( S
) ,\chi ^{+}\right[ $\;.
\\
Then \ $S$ \ possesses all disseminators of this level with all
bases $<\chi ^{+}$ \ on this carrier; analogously for generating
and floating disseminators.
\end{lemma}

\noindent \textit{Proof.} \ Suppose it is wrong, then there is
some
\[
    \rho_{0} \in [ OT(S), \chi^{+} [
\]
such that for every  \  $\rho \in [\rho_{0}, \chi^{+} [ $ \  the
matrix  \  $S $ \  does not possess any disseminator with the base
\ $\rho $. \ In this case standing on \  $\alpha^{\downarrow} $
\  one can detect the ordinal \  $\rho_{1} \in  [ OT(S), \chi^{+}
[ $ \ which can be defined through  \ $\chi $ \ by the formula:
\[
    \exists \delta < \alpha^{\Downarrow} \ SIN_{m}^{<
    \alpha^{\Downarrow} } [<\rho_{1}] (\delta) \wedge \forall \rho
    \in \;] \rho_{1}, \chi^{+} [ \ \forall \delta < \alpha^{\Downarrow}
    \neg SIN_{m}^{< \alpha^{\Downarrow} } [< \rho] (\delta)~.
\]

\noindent But then by lemma 4.6 \ $OT(S) \leq \rho_{1} < OT(S) $.
\hfill $\dashv$
\\

\noindent By the similar argument one can easily prove

\begin{lemma}
\label{6.11.} \hfill {} \\
\hspace*{1em} Let a matrix \ $S$ \ on a carrier \ $\alpha$ \ be
reduced to the cardinal \ $\chi$, \ and let it possess generating
disseminators \ $ \check{\delta}^{\rho _{0}}$, \
$\check{\delta}^{\rho _{1}}$ \ of the level \ $m$ \ with bases \
$\rho_{0}$, \ $\rho_{1}$ \ respectively on this carrier and

\[
    \rho_{0} < OT(S) \leq \rho_{1} <\chi ^{+}.
\]

\noindent Then \qquad \qquad\qquad
$\check{\delta}^{\rho_{0}}<\check{\delta}^{\rho _{1}}$ \quad or
\quad $\check{\delta}^{\rho _{0}}= \check{\delta}^{\rho
_{1}}=\check{\delta}^{\chi ^{+}}$.

\end{lemma}

\noindent Let us describe some methods of producing \ $\delta
$-matrices that are used in what follows. One should pay very
special attention to the following important condition
\[
    A_n^{\vartriangleleft \alpha_1}(\chi^{\ast}) =
    \left\| u_{n}^{\vartriangleleft \alpha _{1}}
    ( \underline{l}) \right\|~,
\]
where, recall, the formula  \ $u_n(\underline{l})$ \ is the
spectral universal \ $\Sigma_n$-formula\ $u_{n}^{\Sigma }$ \
without individual constants, but containing occurrences  of \
$\underline{l}$ \ (see definition 2.6);  this condition plays the
very essential role everywhere further and means that no \
$\Sigma_n$-proposition \ $\varphi(\underline{l})$ \ has jump
ordinals after \ $\chi^{\ast}$ \ below \ $\alpha_1$.

\begin{lemma}
\label{6.12.} {\em (About singular matrices producing)}
\\
\hspace*{1em} Let cardinals \ $\chi <\alpha _{0}<\alpha _{1}$ \
fulfill conditions:
\\

(i)\quad $ \alpha _{0}\in SIN_{n-1}^{<\alpha _{1}}$~; \medskip

(ii)\quad $ A_{n}^{\vartriangleleft \alpha _{1}} ( \chi^{\ast}  )
=\left\| u_{n}^{\vartriangleleft \alpha _{1}} ( \underline{l} )
\right\| $~.
\\
\quad \\
Then there exists the matrix \ $S_{0}$ \ which is reduced to the
cardinal \ $\chi$ \ and singular on the carrier \ $\alpha
_{0}^{\prime }\in \left] \alpha _{0},\alpha _{1} \right[ $ \  such
that \ $\alpha _{0}=\alpha _{0\chi }^{\prime \Downarrow }$.
\\
In this case the matrix \  $S_{0}$ \ is named produced by the
cardinal \  $\alpha _{0}$ \ below \ $\alpha _{1}$.
\end{lemma}
\textit{Proof.} \ Let us apply lemma 5.12. By this lemma the \
$\Pi_n$-proposition \ $\varphi$ 

\[
    \forall \chi^\prime \forall \alpha _{0}^\prime \bigl(  ( \chi^\prime
    \mbox{\ is a limit cardinal }>\omega _{0} ) \longrightarrow
    \exists \alpha _{1}^\prime > \alpha _{0}^\prime ~
    \sigma  ( \chi^\prime ,\alpha _{1}^\prime )  \bigr) .
\]
\vspace{0pt}

\noindent holds; hence it holds below every \ $\Pi_{n-1}$-cardinal
due to lemma~\ref{3.2.} (for \ $n-1$ \ as \ $n$); \ but, moreover,
it holds also below any cardinal \ $\alpha_1$ \ which satisfies
conditions of this lemma.
\\
Justification of this fact is provided by condition $(ii)$, which
implies the preservation of this proposition \ $\varphi$ \ below \
$\alpha_1$. \ It is not hard to see, that in the opposite case the
\ $\Sigma_n$-statement \ $\neg \varphi$ \ receives the jump
ordinal \textit{after \ $\chi^\ast$ \ below \ $\alpha_1$} --
namely, the \textit{minimal} pair \ $(\chi^\prime,
\alpha_0^\prime)$ \ such that the statement

\[
    \bigl(  ( \chi^\prime
    \mbox{\ is a limit cardinal }>\omega _{0} ) \ \wedge \
    \forall \alpha _{1}^\prime > \alpha _{0}^\prime ~
    \neg \sigma  ( \chi^\prime ,\alpha _{1}^\prime )  \bigr)
\]
\vspace{0pt}

\noindent holds below \ $\alpha_1$; \ therefore the \
$\Sigma_n$-universal statement \ $u_n(\underline{l})$ \ receives
this jump ordinal as well due to lemma~\ref{2.7.}. It
implies the following violation of condition \ $(ii)$:
\[
    A_n^{\vartriangleleft \alpha _{1}} ( \chi^{\ast} ) < \|
    u_n^{\vartriangleleft \alpha _{1} } (\underline{l}) \|~.
\]
So, \ $\varphi^{< \alpha_1}$ \ holds; from here it comes the
existence of the cardinal

\[
    \alpha _{0}^{\prime }=\min \{ \alpha \in
    \left] \alpha _{0},\alpha _{1}\right[ :\exists S~\sigma
    ( \chi,\alpha ,S ) \}
\]
\vspace{0pt}

\noindent and of the matrix \  $S_{0}$ \  on the carrier \ $\alpha
_{0}^{\prime }$; \ this minimality provides \ \mbox{$\alpha
_{0}=\alpha _{0\chi }^{\prime \Downarrow }$}. \label{c7}
\endnote{
\ p. \pageref{c7}. \ Here it is quite necessary to pay attention
to the fact, that condition $(ii)$ of this lemma
\[
A_{n}^{\vartriangleleft \alpha _{1}} ( \chi^{\ast}  ) =\left\|
u_{n}^{\vartriangleleft \alpha _{1}} ( \underline{l} ) \right\|
\]
justifies the above-stated proof strongly, and all other
forthcoming reasoning as well. To illustrate in what can result
misunderstanding of this condition the author would like to
present (for an example) the report of one referee in 2000:
\\

\noindent {\em Referee}:  ``Lemma 6.12 is false. Here is a
counterexample. Let \ $\chi$ \ be the ordinal \ $\chi^{\ast}$ \
defined on page 25 [here definition \ref{5.4.} on page
\pageref{5.4.} --- author]. Then for any \ $\alpha_1$, \ condition
$(ii)$ of 6.12 holds. Let \ $\alpha_1$ \ be the \textit{least}
ordinal \ $\alpha$ \ satisfying the following three properties:
\\
\\
(a) \ $\alpha > \chi$;
\\
(b) \ $\alpha \in SIN_{n-2}$;
\\
(c) \ $\sup \left\{ \beta < \alpha: \beta \in SIN_{n-1}^{< \alpha}
\right\} = \alpha$.
\\
\\
Clearly such an \ $\alpha$ \ exists, and any such \ $\alpha$ \
satisfies condition $(i)$ of 6.12. Let \ $\alpha_0$ \ be any
member of \ $SIN_{n-1}^{<\alpha_1}$ \ greater than \ $\chi$. \
Thus the hypothesis of 6.12 is satisfied. But the conclusion of
6.12 is violated. For if \ $S_0$ \ and \ $\alpha_0^{\prime}$ \
satisfy the conclusion, then by definition of ``singular matrix''
and ``carrier'', \ $\alpha_0^{\prime}$ \ satisfies (a)-(c);
but \ $\alpha_0^{\prime} < \alpha_1$ \ contradicting the minimality
of \ $\alpha_1$.
\\
\\
The mistake in the author's proof of 6.12 involves the use of
Lemma 5.12. Now 5.12 is, indeed, true in the model \ $L_k$, \ but
it is not true in the model \ $L_{\alpha_1}$. \ To rephrase this
in the author's notation, if  \ $\varphi$ \  is the statement of
5.12, then \ $\varphi$ \ is true but \ $\varphi^{<\alpha_1}$ \ is
not true. \ $\varphi$ \ is --- as the author points out in the
proof of 5.12
--- a \ $\Pi_n$ \ proposition, and \ $\alpha_1$ \ is merely a
level \ $n-2$, \ not a level \ $n-1$, subinaccessible; so there is
no justification for claiming that \ $\varphi$ \ can be restricted
to \ $\alpha_1$. \ The counterexample shows that it cannot be.''
\\
\\
{\em Commentary of the author in his answer in 2000}: ``No, this
counterexample is false: not for every \ $\alpha_1$ \ condition
$(ii)$ holds; for example, for \ $\alpha_1$ \ used by the referee!
Really, the referee is right that in this case the statement \
$\varphi$ \ of 5.12, used by him, is true in \ $L_k$ \ (and
therefore preserves below \ $\chi^{\ast}$) \ and is not true below
\ $\alpha_1$ \ but he makes the wrong conclusion. Vice versa, the
opposite situation holds: the statement \ $\neg \varphi$ \
receives the jump ordinal after \ $\chi^{\ast}$ \ below \
$\alpha_1$; therefore the \ $\Sigma_n$-universal statement \
$u_n(\underline{l})$ \ receives such ordinal as well (moreover,
infinitely many such ordinals in conditions of the
counterexample). It means the violation of condition (ii), that is
in reality
\[
    A_{n}^{\vartriangleleft \alpha_{1}} ( \chi^{\ast}  ) < \left\|
    u_{n}^{\vartriangleleft \alpha_{1}} ( \underline{l} ) \right\|.
\]
Yes, indeed, \ ``$\varphi$ \ is a $\Pi_n$ \ proposition, and \
$\alpha_1$ \ is merely a level \ $n-2$, \ not a level \ $n-1$'',
but still \ $\varphi$ \ preserves below \ $\alpha_1$ if condition
$(ii)$ fulfills.''
\\
\quad \\
} 
\hfill $\dashv$ 

\begin{lemma}
\label{6.13.} {\em (About $\delta $-matrices producing).}
\\
\hspace*{1em} Let cardinals \ $\chi <\alpha _{0}<\alpha _{1}$ \
fulfill conditions:

\begin{itemize}
\item[(i)] $\alpha _{0}\in SIN_{n-1}^{<\alpha _{1}}$~; \medskip

\item[(ii)] $A_{n}^{\vartriangleleft \alpha _{1}} ( \chi^{\ast} )
=\left\| u_{n}^{\vartriangleleft \alpha _{1}} ( \underline{l} )
\right\| $~; \medskip

\item[(iii)] $\alpha _{0}= \sup (\alpha_0 \cap \sup
SIN_{m-1}^{<\alpha _{1}})$~; \medskip

\item[(iv)]some formula \ $\psi
^{\vartriangleleft \alpha _{0}} ( \beta ,\gamma) $\ defines  on a
final segment of \ $\alpha _{0}$ \ some function \  $f ( \beta  )
=\gamma <\chi ^{+}$ \ which is nondecreasing up to \ $\chi ^{+}$.
\end{itemize}
Then there exists the matrix \  $S_{0}$ \ produced by \ $\alpha
_{0} $ \ below \ $\alpha _{1}$, \ which is reduced to the cardinal
\ $\chi$, \ singular on the carrier \ $\alpha_0^\prime$ \ and
which possesses all generating disseminators on this carrier
\[
    \check{\delta}^{\rho }\in SIN_{m-1}^{<\alpha_{0}}\quad,
    \quad \check{\delta}^{\rho} < \alpha_{0}
\]
of the level \  $m$ \  with all bases \ $\rho <\chi ^{+}$.
\\
In this case we say that the \ $\delta $-matrix \  $S_{0}$ \  is
produced by \ $\alpha _{0}$ \ below \ $\alpha _{1}$.
\\
If in addition \  $\psi \in \Sigma _{m}$, \  then the set of
these disseminators is cofinal to \  $\alpha _{0}$.
\end{lemma}
\textit{Proof.} \ Recall, that by lemma~\ref{6.12.} \ $\alpha_0$ \
is the prejump cardinal \ $\alpha_0^{\prime \Downarrow}$; \ due to
lemma~\ref{6.10.} it is enough to establish the existence of the
disseminator  \ $\delta ^{\rho _{0}}\in SIN_{m-1}^{<\alpha _{0}}$
\ with the base $\rho _{0}=OT ( S_{0} ) $. \ From \textit{(iii)},
\textit{\ (iv)} it comes that there exists a cardinal
\[
    \beta ^{S}\in SIN_{m-1}^{< \alpha_{0}},~~ \beta^S > \chi,
    \mathrm{\quad for\ which\quad} \gamma^{S}=f ( \beta ^{S} )
    >\rho _{0};
\]
from lemmas~\ref{6.4.}, \ref{4.6.} it follows that the ordinal
\[
    \beta^{S}\in SIN_{m}^{<\alpha _{0}} \left[ <\rho _{0}\right].
\]
is the required disseminator.
\\
Really, suppose it is wrong, then
there exists some \ $\Sigma_{m-1}$-formula
 \  $\varphi(\alpha, \overrightarrow{a}) $ \  having the
train of constants \ $< \rho_{0} $ \ such that

\vspace{6pt}
\begin{equation}
\label{e6.6} \forall \alpha < \beta^{S} \varphi^{\triangleleft
\alpha_{0}} (\alpha, \overrightarrow{a}) \wedge \exists \alpha
\in [\beta^{S}, \alpha_{0} [ \ \ \neg \varphi^{\triangleleft
\alpha_{0}} (\alpha, \overrightarrow{a}) ~.
\end{equation}
\vspace{0pt}

\noindent Let us take for some brevity \ $\overrightarrow{a} $ \
consisting  of the one \ $\alpha_{1}^{\prime} < \rho_{0} $ \ only.
It should be used advantage of the jump cardinal
\[
    \alpha_{2} \in dom \left ( \widetilde{\mathbf{S}}_{n}^{\sin
    \triangleleft \alpha_{0}^{\prime}} \overline{\overline{\lceil}}
    \chi \right ),
\]
where, recall, \ $\alpha_{0}^{\prime}$ \ is the minimal carrier $
> \alpha_{0} $ \ of \ $S_{0} $ \ and

\vspace{6pt}
\begin{equation}
\label{e6.7} OT  \left ( \alpha_{2} \cap dom  \left (
\widetilde{\mathbf{S}}_{n}^{\sin \triangleleft
\alpha_{0}^{\prime}} \overline{\overline{\lceil}} \chi \right )
\right ) = \alpha_{1}^{\prime} ~;
\end{equation}
\vspace{0pt}

\noindent the existence of such  \  $\alpha_{2} $ \  is evident
from definition 5.7. After that using (6.6), (6.7) one should
define the ordinal  \  $\gamma \in [ \gamma^{S}, \chi^{+} [ $ \
through  \  $\chi, \alpha_{0}, \alpha_{2} $ \  by the following
formula
\[
    \exists \alpha^{\prime}, \beta < \alpha_{0} \ ( OT
    \left ( dom \left ( \widetilde{\mathbf{S}}_{n}^{\sin
    \triangleleft \alpha_{2}} \overline{\overline{\lceil}}
    \chi \right ) \right ) = \alpha^{\prime} \wedge \forall
    \alpha < \beta \ \ \varphi^{\triangleleft \alpha_{0}}
    (\alpha,\alpha^{\prime}) \wedge
\]
\[
    \qquad\qquad\qquad\qquad\qquad\qquad\qquad\qquad
    \wedge \neg \varphi^{\triangleleft \alpha_{0}} (\beta,
    \alpha^{\prime} ) \wedge  \psi^{\triangleleft \alpha_{0}}
    (\beta,\gamma) ) ~.
\]
\vspace{0pt}

\noindent But then lemma 4.6 implies the contradiction:

\[
    OT(S_{0}) < \gamma < OT(S_{0}) ~.
\]
\vspace{0pt}

\noindent So, \ $\beta^S$ \ is the disseminator with the base \
$\rho_0$ \ and, hence, \ $S_0$ \ possesses generating
disseminators with all bases \ $\rho < \chi^+$.
\\
Now suppose that \ $\psi \in \Sigma_{m} $ \ but, on the contrary,
the set of disseminators is not cofinal to \ $\alpha_{0} $. \ By
lemma~6.10, 6.8~3) in this case there exists the disseminator

\[
    \beta_{1} \in SIN_{m}^{< \alpha_{0}} [< \chi^{+}]
\]
\vspace{0pt}

\noindent having the base  \ $\chi^{+} $. \ For  \
$\gamma_{1}=f(\beta_{1}) $ \ the cardinal \ $\beta_{1} $ \
extends up to  \  $\alpha_{0} $ \ the  \ $\Pi_{m}$-proposition

\[
    \forall \beta, \gamma \quad (\psi(\beta, \gamma)
    \longrightarrow \gamma < \gamma_{1} )
\]
\vspace{0pt}

\noindent and the function  \  $f $ \  becomes bounded by  \
$\gamma_{1} < \chi^{+} $ \  contrary to  \  $(iv) $.\hfill
$\dashv$
\\
\quad \\
The similar reasoning establishes the following important

\begin{lemma}
\label{6.14.} {\em (About \ $\delta$-matrix functions producing)}
\\
\hspace*{1em} Let cardinals \ $\chi <\alpha _{1}$ \ fulfill the
conditions: \medskip
\begin{itemize}
\item[(i)] $A_{n}^{\vartriangleleft \alpha _{1}} ( \chi^{\ast} ) =\left\|
u_{n}^{\vartriangleleft \alpha _{1}} ( \underline{l} )
\right\|$~; \medskip

\item[(ii)] $\alpha _{1}=\sup SIN_{m-1}^{<\alpha _{1}}=\sup
SIN_{n-1}^{<\alpha _{1}}$~; \medskip

\item[(iii)] $cf ( \alpha_{1} ) >\chi ^{+}$~; \medskip

\item[(iv)] the cardinal \ $\chi \ge \chi^{\ast}$ \ is definable through \
$\chi^\ast$ \ by a formula of the class \ $\Sigma_{n-2} \cup
\Pi_{n-2}$~.
\end{itemize}
Then there exists the cardinal \  $\gamma _{\tau _{0}}^{<\alpha
_{1}}$ \ such that on the set
\[
    T=\{ \tau : \gamma _{\tau_{0}}^{<\alpha _{1}}
    \leq \gamma _{\tau }^{<\alpha _{1}}\}
\]
there is defined the function \ $\mathfrak{A}$ \ which possesses
the following properties for every  \ $\tau \in T$ :\
\\

{\em 1)}\quad\quad $\mathfrak{A} ( \tau  ) $ \ is the \ $\delta
$-matrix \ $S$ \ of the level \ $m$ \ reduced to \ $\chi$ \ and
admissible on some carrier \ $\alpha \in \left] \gamma _{\tau
}^{<\alpha _{1}}, \gamma _{\tau + 1}^{<\alpha _{1}} \right[ $ \
for \ $\gamma _{\tau }^{<\alpha _{1}}$ \ below \ $\alpha_1$;
\\

{\em 2)}\quad the set of all generating disseminators of \ $S$ \
on \ $\alpha$ \  of this level is cofinal to \ $\alpha _{\chi
}^{\Downarrow }$;
\\

{\em 3)}\quad $cf ( \alpha _{\chi }^{\Downarrow } ) =\chi ^{+}$;
\\

{\em 4)}\quad $( \gamma _{\tau }^{<\alpha _{1}}+1 ) \cap
SIN_{m-1}^{<\alpha _{1}}\subseteq SIN_{m-1}^{<\alpha _{\chi
}^{\Downarrow }}$~.
\\
\quad \\
In such case we say that the function  \  $\mathfrak{A}$ \  is
produced by the cardinal \  $ \alpha _{1}$ \ and name it the
matrix \ $\delta$-function.
\end{lemma}
\textit{Proof} consists in the repetition of the function  \
$S_{f}^{< \alpha_{1}} $ \  construction from definition
\ref{5.14.} slightly changed. The upper index \  $< \alpha_{1} $ \
will be omitted for shortness. Let us apply lemma~6.12 to every
cardinal \ $\gamma_{\tau_{1}} > \chi $, \ $\gamma_{\tau_{1}} \in
SIN_{m-1} $ \  of the cofinality \ $\neq \chi^{+} $; \ the
existence of such cardinals follows from \ $(iii) $ \  and
lemma~3.4~1) (for \ $n=m-1$). Evidently the matrix \  $S $ \  on
the carrier  \ $\alpha \in \; ]
\gamma_{\tau_{1}},\gamma_{\tau_{1+1}} [ $ \ produced by  \
$\gamma_{\tau_{1}} $ \  possesses the following property for \
$\tau = \tau_{1} $:
\begin{equation}
\label{e6.8}cf(\alpha^{\Downarrow}) \neq \chi^{+} \wedge \forall
\gamma \leq \gamma_{\tau}  \left ( SIN_{m-1}(\gamma)
\longrightarrow SIN_{m-1}^{< \alpha^{\Downarrow}} (\gamma) \right
)~.
\end{equation}

By lemma~3.2 about restriction (for  \  $n-1 $ \  instead of
 \  $n $) the carriers  \  $\alpha $ \  of this kind
exist in every  \  $ [ \gamma_{\tau},\gamma_{\tau_{+1}} [ $ \ for
 \  $\chi < \gamma_{\tau_{+1}} \leq  \gamma_{\tau_{1}}
$ \  (see the proof of lemma~5.17~2)~(ii)~). Really, let us
consider the \  $\Pi_{n-2}$-formula \ $\varphi(\alpha, \chi,
\gamma^{m}_{\tau}, \gamma_{\tau}, S)$:

\vspace{6pt}
\[
    \gamma_{\tau} < \alpha  \wedge cf(\alpha^{\Downarrow})
    \neq \chi^{+} \wedge \sigma (\chi, \alpha, S) \wedge
    SIN_{m-1}^{< \alpha^{\Downarrow}} (\gamma_{\tau}^{m})
\]

\noindent where

\[
    \gamma_{\tau}^{m} = \sup  \left ( (\gamma_{\tau}+1)
    \cap SIN_{m-1} \right ).
\]
\vspace{0pt}

\noindent The \ $ \Sigma_{n-1} $-proposition \ $ \exists \alpha \
\ \varphi $ \  has constants
\[
    \chi, \quad \gamma_{\tau}^{m}, \quad \gamma_{\tau} < \gamma_{\tau+1},
    \quad S \vartriangleleft \gamma_{\tau+1},
\]
so the  \ $SIN_{n-1} $-cardinal  \  $ \gamma_{\tau+1} $ \
restricts this proposition and the matrix  \  $S $ \  appears on
the carrier  \  $ \alpha \in [ \chi_{\tau}, \chi_{\tau+1} [ $ \
that fulfills (6.8).
\\
After that let us repeat definition 5.14 of the function  \  $
S_{f} $ \  below  \  $\alpha_{1} $ \  but imposing the additional
requirement (6.8) upon the carrier  \  $\alpha $; \ as a result we
obtain the function  \  $S_{f}^{m} = (S_{\tau}^{m})_{\tau} $ \
taking values:

\vspace{6pt}
\[
    S_{\tau}^{m}=\min_{\underline{\lessdot }}\{ S:\exists
    \alpha < \alpha_{1}(\gamma_{\tau}<\alpha \wedge
    cf(\alpha^{\Downarrow}) \neq \chi^{+} \wedge
    \sigma ( \chi ,\alpha ,S ) \wedge
\]

\vspace{-6pt}

\[
    \qquad \qquad \qquad \qquad \qquad
    \wedge \forall \gamma \leq \gamma_{\tau} \ \
    (SIN_{m-1}(\gamma) \longrightarrow SIN_{m-1}^{<
    \alpha^{\Downarrow}}(\gamma))) \} ~.
\]
\vspace{0pt}

\noindent This function  \  $\underline{\lessdot } $-nondecreases
just as the function  \  $S_{f} $ \  in lemma 5.17, so due to
\textit{(iii)} there exists an ordinal \ $\tau^{m} \in dom \left (
S_{f}^{m} \right )$ \ such that
\[
    \gamma_{\tau^{m}} \in SIN_{m-1} \mathrm{\quad and\ for\ every \quad}
    \tau \geq \tau^{m} \quad S_{\tau}^{m} = S_{\tau^{m}}^{m}~.
\]
Let us consider an arbitrary  \  $\tau
> \tau^{m} $ \  and the carrier  \  $\alpha_{\tau}^{m} \in [
\gamma_{\tau}, \gamma_{\tau+1}[ $ \  of the matrix  \
$S_{\tau}^{m} $. \ The definition of  \  $S_{f}^{m} $ \  should be
repeated below  \  $\alpha_{\tau}^{1} = \alpha_{\tau}^{m
\Downarrow} $ \  and the resulting function  \
$S_{f}^{m<\alpha_{\tau}^{1}} $ \  is still \ $\underline{\lessdot
} $-nondecreases and coincides with \ $S_{f}^{m} $ \  on  \ $\tau
$. \ The reasoning from the end of the previous section shows with
the help of condition $(iv)$, that it \ $\underline{\lessdot }
$-nondecreases up to \ $\chi^{+} $ \ and, so, appears the ordinal

\[
    \alpha_{0 \tau} = sup \left \{ \gamma_{\tau^{\prime}}^{<
    \alpha_{\tau}^{1}} : \tau^{\prime} \in dom \left (
    S_{f}^{m < \alpha_{\tau}^{1}} \right ) \right \}
\]
such that
\[
    \alpha_{0 \tau} = sup  \left ( \alpha_{0 \tau} \cap
    SIN_{m-1}^{< \alpha_{\tau}^{1}} \right ).
\]
\vspace{0pt}

\noindent It is obvious that \ $ cf(\alpha_{0 \tau}) = \chi^{+} $,
\ hence condition \ $ cf(\alpha_{\tau}^{1}) \neq \chi^{+} $ \
implies
\[
    \gamma_{\tau}< \alpha_{0 \tau} < \alpha_{\tau}^{1}.
\]
Now we can apply lemma~6.13 to  \  $\alpha_{0 \tau}$, \
$\alpha_{\tau}^{1} $ \ as \ $\alpha_{0}$, \ $\alpha_{1} $ \
respectively and to the function
\[
    f(\beta)= sup  \left \{ Od  \left ( S_{\tau^{\prime}}^{m <
    \alpha_{\tau}^{1}} \right ) : \gamma_{\tau^{\prime}}^{<
    \alpha_{\tau}^{1}} < \beta \right \} ~,
\]
which meets all requirements of this lemma.

\noindent The matrix \ $S = S_{0 \tau} $ \ produced by \
$\alpha_{0 \tau} $ \ below  \ $\alpha_{\tau}^{1} $ \ on the
carrier \ $\alpha_{0 \tau}^{\prime} $ \ possesses generating
disseminators of the level  \  $m $ \  with arbitrary bases  \  $<
\chi^{+} $. \ It is not hard to see that the set of these
disseminators is cofinal to  \ $\alpha_{0 \tau} = \alpha_{0
\tau}^{\prime \Downarrow} $ \  and that it has some generating
disseminator \ $\check{\delta} < \gamma_{\tau} $ \ admissible for
\ $\gamma_\tau$.
\\
It comes out that the function  \ $\mathfrak{A}(\tau) = S_{0 \tau}
$ \ fulfills all conclusions of this lemma completely. \hfill
$\dashv$
\\
\hfill {} \\

Now the basic theory is developed sufficiently for the further analysis of
matrix functions, which will be undertaken in Part II of this
work forthcoming. \label{c8}
\endnote{
\ p. \pageref{c8}. \ This section can be finished by one more
comment on the disposition of disseminators. In
lemmas~\ref{6.13.},~\ref{6.14.} \ $\delta $-matrices produced on
the carriers \  $\alpha $ \ have generating disseminators that are
disposed cofinally to \ $\alpha _{\chi }^{\Downarrow }$ \  and \
$cf ( \alpha _{\chi }^{\Downarrow } ) =\chi ^{+}$. \ Using methods
of reasoning from proofs of lemmas 6.13, 6.14 one can prove that
for \ $\underline{\lessdot }$-minimal \  $\delta $-matrices and \
$m \geq n+1$ \ this is unavoidable:
\\

\noindent {\bf Lemma} \\
\hspace*{1em} \em Let cardinals \ $\chi <\gamma
_{\tau}^{<\alpha_{1}}$ \  and a matrix \  $S $ \  fulfill the
conditions:
\\

(i)\quad $ A_{n}^{\vartriangleleft \alpha _{1}} ( \chi^{\ast} )
=\left \| u_{n}^{\vartriangleleft \alpha _{1}} ( \underline{l} )
\right\| $ ~;
\\

(ii)\quad $S$ \ is the \ $\delta $-matrix of the level \ $m \geq
n+1$ \ admissible on some carrier \ $\alpha \in \left] \gamma
_{\tau }^{<\alpha _{1}},\alpha _{1}\right[ $ \ for \ $\gamma
_{\tau }^{<\alpha _{1}}$ \ below \ $\alpha_1$;
\\

(iii)\quad $S$ \ is \ $\underline{\lessdot }$-minimal of all \
$\delta $-matrices with this property.
\\
\quad \\
Then \ $S$ \ on the carrier \  $\alpha $ \ possesses generating
disseminators of the level \  $m$  \  with all bases \  $\rho
<\chi ^{+}$ \  disposed cofinally to \  $\alpha _{\chi
}^{\Downarrow }$ \  and \  $ cf ( \alpha _{\chi }^{\Downarrow } )
=\chi ^{+}$. \em
\\

\noindent This lemma is not used in this work and so the proof is
omitted.
\\
\quad \\
} 

\theendnotes

\newpage

\thispagestyle{empty}

\begin{center}
\quad \\
\quad \\
\quad \\
\quad \\
\quad \\
\quad \\
\quad \\
\quad \\
\quad \\
{Scientific edition} \\
\quad \\
{\bf Kiselev} Alexander\\
\quad \\
{\bf INACCESSIBILITY \\
AND \\
SUBINACCESSIBILITY}\\
\quad \\
{In two parts}\\
{Part I}\\
\quad \\
{\small Responsible for release {\em T. E. Yanchuk}} \\
\quad \\
{\footnotesize Signed for publication on 14.04.2008. Format
60$\times$84/16. Paper offset. Typeface Times. Risograph.
Conditional quires 6,51. Registration-publishing sheets 5,85.
Circulation is 100 ex. Order \selectlanguage{russian} \No
\selectlanguage{english} 359.} \\
\quad \\
{\footnotesize Republican unitary enterprise} \\
{\footnotesize ``Publishing center of Belarusian State University''} \\
{\footnotesize \selectlanguage{russian} ╦╚ \No
\selectlanguage{english}
02330/0131748 on 01.04.2004.} \\
{\footnotesize 220030, Minsk, Krasnoarmeyskaya Str., 6.}\\
\quad \\
{\footnotesize Printed from the author's layout in Republican
unitary
enterprise} \\
{\footnotesize ``Publishing center of Belarusian State University''} \\
{\footnotesize \selectlanguage{russian} ╦╧ \No
\selectlanguage{english} 02330/0056850 on 30.04.2004.}
\end{center}
\label{end}

\newpage

\thispagestyle{empty}

\end{document}